\documentclass[11pt]{amsart}
\usepackage{amsmath}
\usepackage{amssymb}
\usepackage{bm}
\usepackage{esint}
\usepackage{graphicx}
\usepackage{verbatim}
\usepackage{enumitem}
\usepackage{stackengine}
\usepackage{algorithm}
\usepackage{algpseudocode}
\stackMath
\setlist[enumerate,1]{leftmargin=*}
\setlist[itemize,1]{leftmargin=*}
\usepackage{mathrsfs}
\usepackage[dvipsnames]{color}
\usepackage{todonotes}
\usepackage[pdftex,pdfstartview=FitH,pdfborderstyle={/S/B/W 1},%
colorlinks=true, linkcolor=blue, urlcolor=blue, citecolor=blue,%
pagebackref=true]{hyperref}
\usepackage[nobysame,alphabetic,initials,msc-links]{amsrefs}
\usepackage{imakeidx}
\usepackage{tikz}
\usepackage{multicol}
\usetikzlibrary{decorations.markings}
\usepackage{mathtools}
\usepackage{mathrsfs}
\usepackage{comment}
\usepackage{bbm}

\newtheorem{Lemma}{Lemma}[section]
\newtheorem{Theorem}[Lemma]{Theorem}
\newtheorem{Proposition}[Lemma]{Proposition}
\newtheorem{Corollary}[Lemma]{Corollary}

\newtheorem{Claim}[Lemma]{Claim}

\newtheorem*{lem:longHammingSep}{Lemma  \ref{lem:longHammingSep}}
\newtheorem*{prop:cauchySeq}{Proposition  \ref{prop:cauchySeq}}
\newtheorem*{cor:sl2ergodic}{Corollary  \ref{cor:sl2ergodic}}

\theoremstyle{definition}

\newtheorem{Definition}[Lemma]{Definition}

\newtheorem{Remark}[Lemma]{Remark}

\numberwithin{equation}{section}

\newcommand {\N} {{\mathbb N}}
\newcommand {\Q} {{\mathbb Q}}
\newcommand {\R} {{\mathbb R}}

\newcommand {\T} {{\mathbb T}}

\newcommand {\Z} {{\mathbb Z}}
\newcommand {\absolute}[1] {\left| {#1} \right|}
\newcommand {\norm}[1] {\left\| {#1} \right\|}
\newcommand {\innproduct}[1] {\left\langle {#1} \right\rangle}

\DeclareMathOperator {\PSL} {PSL}
\DeclareMathOperator {\SL} {SL}
\DeclareMathOperator {\slfrak} {\mathfrak{sl}}
\DeclareMathOperator {\SO} {SO}

\DeclareMathOperator {\SU} {SU}

\DeclareMathOperator {\Kak} {Kak}
\DeclareMathOperator {\id} {id}
\DeclareMathOperator {\ad} {ad}
\DeclareMathOperator {\Ad} {Ad}

\DeclareMathOperator {\GL} {GL}
\DeclareMathOperator {\inj} {inj}

\newcommand {\bu} {\mathbf{u}}
\newcommand{\bou} {\hat{\mathbf{u}}}
\newcommand {\ba} {\mathbf{a}}
\newcommand {\bx} {\mathbf{x}}

\newcommand {\tr} {\mathrm{tr}}
\DeclareMathOperator {\Lie} {\operatorname{Lie}}

\DeclareMathSymbol{\shortminus}{\mathbin}{AMSa}{"39}

\newcounter{c}

\newcounter{C}
\newcommand{\Cnew}{%
\refstepcounter{C}%
\ensuremath{C_{\theC}}}
\newcommand{\Cold}[1]{\ensuremath{C_{\ref{#1}}}}

\newcounter{e}
\newcommand{\enew}{%
\refstepcounter{e}%
\ensuremath{\epsilon_{\thee}}}
\newcommand{\eold}[1]{\ensuremath{\epsilon_{\ref{#1}}}}

\newcounter{de}
\newcommand{\denew}{%
\refstepcounter{de}%
\ensuremath{\delta_{\thede}}}
\newcommand{\deold}[1]{\ensuremath{\delta_{\ref{#1}}}}

\newcounter{d}

\newcounter{R}
\newcommand{\Rnew}{%
\refstepcounter{R}%
\ensuremath{R_{\theR}}}
\newcommand{\Rold}[1]{\ensuremath{R_{\ref{#1}}}}

\newcounter{r}

\newcounter{K}
\newcommand{\Knew}{%
\refstepcounter{K}%
\ensuremath{{\mathcal K}_{\theK}}}
\newcommand{\Kold}[1]{\ensuremath{{\mathcal K}_{\ref{#1}}}}

\newcounter{O}

\newcounter{Z}
\newcommand{\Znew}{%
\refstepcounter{Z}%
\ensuremath{{\mathcal Z}_{\theZ}}}
\newcommand{\Zold}[1]{\ensuremath{{\mathcal Z}_{\ref{#1}}}}

\newcounter{k}

\newcounter{w}

\newcounter{ps}

\newcounter{th}

\newcounter{n}

\newcounter{omega}

\makeindex[title=Index of notation,columns=2]

\newcommand{\Uh}{\widehat{U}}

\newcommand{\xt}{\Tilde{x}}
\newcommand{\yt}{\Tilde{y}}

\newcommand{\sol}{\overline{s}}

\newcommand{\ebf}{\mathbf{e}}

\newcommand{\Gbf}{\mathbf{G}}
\newcommand{\Hbf}{\mathbf{H}}
\newcommand{\Nbf}{\mathbf{N}}

\newcommand{\Lbf}{\mathbf{L}}
\newcommand{\Sbf}{\mathbf{S}}
\newcommand{\Mbf}{\mathbf{M}}
\newcommand{\Zbf}{\mathbf{Z}}
\newcommand{\Dbf}{\mathbf{D}}
\newcommand{\Ebf}{\mathbf{E}}
\newcommand{\Pbf}{\mathbf{P}}
\newcommand{\Qbf}{\mathbf{Q}}
\newcommand{\Rbf}{\mathbf{R}}
\newcommand{\Tbf}{\mathbf{T}}
\newcommand{\Wbf}{\mathbf{W}}

\newcommand{\gfrak}{\mathfrak{g}}
\newcommand{\hfrak}{\mathfrak{h}}
\newcommand{\zfrak}{\mathfrak{z}}
\newcommand{\lfrak}{\mathfrak{l}}
\newcommand{\nfrak}{\mathfrak{n}}
\newcommand{\mfrak}{\mathfrak{m}}
\newcommand{\sfrak}{\mathfrak{s}}

\newcommand{\ffrak}{\mathfrak{f}}
\newcommand{\pfrak}{\mathfrak{p}}
\newcommand{\qfrak}{\mathfrak{q}}

\newcommand{\Ocal}{\mathcal{O}}
\newcommand{\Hcal}{\mathcal{H}}

\newcommand{\Core}{\operatorname{Core}}

\newcommand{\Sur}{\operatorname{Sur}}
\newcommand{\nc}{\operatorname{nc}}

\begin{document}

\title[Time change rigidity]{Rigidity of time changes of $\slfrak_2$-induced unipotent flows on finite volume quotients}
\author[H.-H. Nguyen]{Huu-Hoan Nguyen}
\thanks{H.-H.N. acknowledges President's Graduate Fellowship (PGF) Scholarship from NUS}
\address{H.-H.N.: Department of Mathematics, National University of Singapore, 119076, Singapore}
\email{hoancats196@gmail.com, e1353503@u.nus.edu}

\author[D. Wei]{Daren Wei}
\thanks{D.W. acknowledges support by NUS startup grants A-0009806-00-00 and A-0009806-01-00}
\address{D.W.: Department of Mathematics, National University of Singapore, 119076, Singapore}
\email{darenwei@nus.edu.sg}

\setlength{\marginparwidth }{2 cm}

\date{\today}

\begin{abstract}
A unipotent flow $U_t=\exp(t\mathbf{u})$ on a finite volume quotient $\mathbf{G}/\Gamma$ of the identity component $\mathbf{G}$ of a real linear algebraic group is $\mathfrak{sl}_2$-induced if $0\neq \mathbf{u}\in\mathfrak{h}\cong\mathfrak{sl}_2(\mathbb{R})$ for some Lie subalgebra $\mathfrak{h}\subseteq\operatorname{Lie}(\mathbf{G})$, and splitting if $\mathfrak{h}$ is an ideal in $\operatorname{Lie}(\mathbf{G})$. We prove that every measurable factor map between time changed ergodic nonsplitting $\mathfrak{sl}_2$-induced unipotent flows is affine up to flow translations. Hence an ergodic $\mathfrak{sl}_2$-induced unipotent flow is loosely Kronecker if and only if it is splitting, and two ergodic $\mathfrak{sl}_2$-induced unipotent flows are Kakutani equivalent if and only if they are isomorphic or both splitting.

Here $\mathbf{G}$ need not to be semisimple, and $L_0^2(\mathbf{G}/\Gamma)$ need not have a spectral gap. Our proof relies on a decomposition of $\mathbf{G}$ and $\mathbf{G}/\Gamma$, strong relative property (T), and removal of flow translations along compact and central directions via solving cohomological equation.
\end{abstract}

\maketitle

\setcounter{tocdepth}{1}
\tableofcontents

\section{Introduction}\label{Section:Introduction}
A landmark phenomenon of homogeneous dynamics is the rigidity of unipotent flows on finite volume quotients, which is established in seminal works of Dani, Ratner and Witte \cites{MR578655,MR629475,MR657240,MR717825,MR796188,MR910358,MR1075042,MR1135878}. It follows from these celebrated works that orbit closures, ergodic measures, and joinings of unipotent flows have a rigid algebraic structure. 

One of the natural questions is how much of this rigidity survives under measurable time changes. Equivalently, one asks whether rigidity is preserved under \textbf{Kakutani equivalence} \cites{kakutani1943131,katok1977monotone,MR653094}, which preserves the orbit monotonicity while distorting the speed along orbits. Already on semisimple quotients, Ratner proved that horocycle flows on finite volume quotients of $\SL_2(\R)$ are \textbf{loosely Kronecker} \cite{MR516248} and hence all mutually Kakutani equivalent, whereas products of horocycle flows exhibit rigidity under Kakutani equivalence \cites{MR571396,ratner1981some}. Under  additional regularity assumptions on the time changes, rigidity holds even for horocycle flows \cites{MR822329,MR880955}, with latter extensions to broader semisimple settings in \cites{MR4385813,MR4651586,artigiani2024rigiditypropertiestimechangesunipotent}. Without regularity assumptions, Kanigowski, Vinhage and the second author classified all loosely Kronecker unipotent flows on semisimple quotients \cite{MR4255058}, and Lindenstrauss and the second author recently established full time change rigidity on quotients of semisimple Lie groups without compact factors \cite{lindenstrauss2025timechangerigidityunipotent}.

\medskip

Witte \cites{MR796188,MR910358} and Ratner \cite{MR1135878} proved that for classical unipotent flows, the rigidity phenomena hold on finite volume quotients $\mathbf{G}/\Gamma$ of any connected Lie group $\Gbf$. On nilmanifolds $N/\Lambda$, it follows from \cites{katok1977monotone,MR653094} that every aperiodic ergodic nilflows is loosely Kronecker, which shows that time change rigidity may fail without semisimple structure. Even if our unipotent flow $U_t=\exp(t\bu)$ on $\Gbf/\Gamma$ is $\mathfrak{sl}_2$-induced, with $0\neq\bu\in\hfrak\cong\slfrak_2(\R)$, the base group $\Gbf$ may still contain a nontrivial unipotent radical, compact semisimple factors and a positive dimensional center $Z(\Gbf)$. Here $H=\langle\exp(\hfrak)\rangle$ commutes with the compact semisimple factors and $Z(\Gbf)$,  and $L_0^2(\Gbf/\Gamma)$ need not have a spectral gap. The purpose of this paper is to establish the time change rigidity and a Kakutani equivalence dichotomy for all ergodic $\slfrak_2$-induced unipotent flows on the quotients $\Gbf/\Gamma$ of the identity component $\Gbf$ of any real linear algebraic group.

\medskip

To state our results precisely, let $\Gbf_i$ be a connected Lie group, $\Gamma_i<\Gbf_i$ a lattice, and $\mu_i$ the Haar probability measure on $\Gbf_i/\Gamma_i$, $i=1,2$. Given a subgroup $H<\Gbf_i$, its \textbf{normal core} $\Core_{\Gbf_i}(H)$ is the largest normal subgroup of $\Gbf_i$ that lies in $H$. For $i=1,2$, let
\[
\Gbf'_i=\Gbf_i/\Core_{\Gbf_i}(\Gamma_i),\quad \Gamma'_i=\Gamma_i/\Core_{\Gbf_i}(\Gamma_i),\quad \vartheta_i(g\Gamma_i)=(g\Core_{\Gbf_i}(\Gamma_i))\Gamma'_i.
\]
A measurable map $\phi:\Gbf_1/\Gamma_1\to\Gbf_2/\Gamma_2$ is called \textbf{affine} if there exists a pair $(g_0,\alpha)\in\Gbf'_2\times \Sur(\Gbf'_1,\Gbf'_2)$, where $\Sur(\Gbf'_1,\Gbf'_2)$ denotes the set of continuous surjective homomorphisms $\Gbf'_1\to\Gbf'_2$, such that 
\[
\alpha(\Gamma'_1)\subseteq g_0\Gamma'_2g_0^{-1}
\]
and for $\mu'_1=(\vartheta_1)_*\mu_1$-a.e. $g\Gamma'_1\in\Gbf'_1/\Gamma'_1$,
\[
\vartheta_2\circ\phi\circ\vartheta_1^{-1}(g\Gamma'_1)=\alpha(g)g_0\Gamma'_2. 
\]
Such a pair $(g_0,\alpha)$ is called an \textbf{affine pair associated with $\phi$}.

With these notions, the isomorphism rigidity results of Ratner and Witte \cites{MR796188,MR910358,MR1075042} can be stated as follows: a measurable factor map between ergodic unipotent flows on quotients of general connected Lie groups is affine. Moreover, if $\phi$ is a measurable isomorphism, then any affine pair $(g_0,\alpha)$ associated with $\phi$ satisfies 
\[
\alpha(\Gamma'_1)=g_0\Gamma'_2g_0^{-1}
\]
and $\alpha:\Gbf'_1\to\Gbf'_2$ is a Lie group isomorphism.

\medskip

Let $(X,\mu, T_t)$ be a probability measure preserving flow. The \textbf{time changed flow} of $T_t$ under a \textbf{time change} $\tau\in L^1_+(X,\mu)$, denoted by $T^{\tau}_t$, is given as
\[
T^\tau_tx=T_{\rho(x,t)}x,\qquad\mu\text{-a.e. }x\in X,\quad\forall t\in\R,
\]
where $\rho(x,t)\in\R$ satisfies $\int_0^{\rho(x,t)}\tau(T_sx)ds=t$. The flow $T^{\tau}_t$ preserves the probability measure $d\mu_\tau=\frac{\tau d\mu}{\int_X\tau d\mu}$. Moreover, if $(X,\mu,T_t)$ is ergodic, then so is $(X,\mu_\tau,T^{\tau}_t)$.

For probability measure preserving flows $(X,\mu,T_t)$ and  $(Y,\nu,S_t)$, a measurable map $\psi:X\to Y$ is called a \textbf{Kakutani equivalence} between two flows if there is a time change $\tau\in L^1_+(X,\mu)$ such that $\psi$ is a measurable isomorphism between $(X,\mu_{\tau},T^{\tau}_t)$ and $(Y,\nu,S_t)$. In this case, $(X,\mu,T_t)$ and $(Y,\nu,S_t)$ are said to be \textbf{Kakutani equivalent}. The notion of \textbf{Kakutani factor map} is defined in the same way, where \emph{measurable isomorphism} is replaced by \emph{factor map}. For more details, we refer to \cites{kakutani1943131,katok1977monotone, MR653094}.

Since time change preserves the monotone orbit structure but may change the speed of the flow, rigidity up to time changes naturally allows translations along target orbits. Let $\Gbf_i$ be a connected Lie group and $\Gamma_i<\Gbf_i$ a lattice for $i=1,2$. Suppose that $U_t$ is a volume preserving flow on $\Gbf_2$, a measurable map $\psi:\Gbf_1/\Gamma_1\to\Gbf_2/\Gamma_2$ is called \textbf{affine up to $U_t$-translations} if for some measurable function $\sigma:\Gbf_1/\Gamma_1\to\R$, the composition map $U_{\sigma(\cdot)}\psi(\cdot)$ is affine. An affine pair $(g_0,\alpha)\in\Gbf_2'\times\Sur(\Gbf_1',\Gbf_2')$ associated with $U_{\sigma(\cdot)}\psi(\cdot)$ is called an \textbf{affine pair associated with $\psi$ up to $U_t$-translations}.

\medskip

Recall that a flow is \textbf{loosely Kronecker} if it is Kakutani equivalent to an irrational linear flow on the $2$-torus $\T^2$. Following foundational rigidity results for products of horocycle flows, Lorentz quotients, regular time changes and classification of loosely Kronecker unipotent flows on semisimple quotients \cites{MR822329,MR880955,MR4385813,MR4651586,MR5152412}, full measurable time change rigidity on the semisimple quotients without compact factors was established in \cite{lindenstrauss2025timechangerigidityunipotent}: if $\psi:\mathbf{G}_1/\Gamma_1\to\mathbf{G}_2/\Gamma_2$ is a Kakutani factor map between non-loosely Kronecker ergodic unipotent flows $(\mathbf{G}_i/\Gamma_i,\mu_i,U^{(i)}_t),i=1,2$, on quotients of semisimple groups and $\Gbf_1$ does not have compact factors, then $\psi$ is affine up to $U^{(2)}_t$-translations. Moreover, if $\psi$ is a Kakutani equivalence, then additionally, an affine pair $(g_0,\alpha)\in\Gbf_2'\times\Sur(\Gbf_1',\Gbf_2')$ associated with $\psi$ up to $U^{(2)}_t$-translations satisfies $\alpha(\Gamma_1')=g_0\Gamma_2'g_0^{-1}$ and $\alpha$ is an isomorphism.

\medskip

On general algebraic quotients, the natural setting for our rigidity theorems is the class of ergodic \textbf{$\slfrak_2$-induced} unipotent flows. Let $\Gbf$ be the identity component of a real linear algebraic group, $\gfrak=\Lie(\Gbf)$ the corresponding Lie algebra, $\Gamma<\Gbf$ a lattice, and $\mu$ the Haar probability measure on $\Gbf/\Gamma$. A unipotent flow $U_t=\exp(t\bu)$ on $\Gbf/\Gamma$ is  \textbf{$\slfrak_2$-induced} if there exists a Lie subalgebra $\hfrak\subseteq\gfrak$ satisfying $\hfrak\cong\slfrak_2(\R)$ and $\bu\in\hfrak\setminus\{0\}$; any such $\hfrak$ is called an $\slfrak_2(\R)$-subalgebra corresponding to $U_t$. By the \textit{Jacobson-Morozov theorem} (cf. \cite{MR1920389}*{Theorem 10.3}), every nontrivial unipotent flow on a finite volume quotient of a connected semisimple Lie group is $\slfrak_2$-induced.

When $\Gbf$ is not semisimple, a unipotent flow on $\Gbf/\Gamma$ need not be $\slfrak_2$-induced. On any nilmanifold $N/\Lambda$ with $\dim N\geq 2$, every ergodic nilflow is loosely Kronecker and hence no rigidity up to measurable time changes \cites{katok1977monotone,MR653094}. Non $\slfrak_2$-induced generators also happen when $\Gbf$ has a noncompact Levi factor. For $\bu_0=\big(\big(\begin{smallmatrix}
    0 & 1 \\ 0 & 0
\end{smallmatrix}\big),(0,1)^{\mathsf{T}}\big)\in\slfrak_2(\R)\ltimes\R^2$, a direct check gives $[\ba,\bu_0]\neq2\bu_0$ for all $\ba\in\slfrak_2(\R)\ltimes\R^2$, so $\bu_0$ does not belong to any $\slfrak_2(\R)$-subalgebra. Assuming $U_t=\exp(t\bu)$ is $\slfrak_2$-induced is a natural assumption that provides an $\slfrak_2$-triple so that the renormalization and matching mechanism can be activated.

Given an $\slfrak_2$-induced unipotent flows $U_t$ on $\Gbf/\Gamma$ with corresponding subalgebra $\hfrak\cong\slfrak_2(\R)$, we call $U_t$ \textbf{splitting} if $\hfrak$ is an ideal of $\gfrak$, and \textbf{nonsplitting} otherwise. Since $\slfrak_2(\R)$ is simple, $\hfrak$ is an ideal of $\gfrak$ if and only if $\gfrak=\hfrak\oplus\mathfrak{g}'$ as a Lie algebra direct sum with ideal $\mathfrak{g}'=C_{\gfrak}(\hfrak)$, in which case $\hfrak$ is the unique $\slfrak_2(\R)$ subalgebra of $\gfrak$ containing $\bu$. It follows from \cites{MR516248,ratner1981some,MR4255058} that on semisimple quotients, an ergodic unipotent flow is splitting if and only if it is loosely Kronecker. As noted in \cites{MR822329,MR880955,MR5152412}, for such flows, the rigidity phenomenon fails completely unless additional regularity conditions are imposed on the time changes.

\medskip

To state our results concisely, we introduce the following terminology.
\begin{Definition}[Time change rigidity]\label{def:TimeChangeRigid}
    An ergodic nonsplitting $\slfrak_2$-induced unipotent flow $(\Gbf_1/\Gamma_1,\mu_1,U^{(1)}_t)$ is \textbf{time change rigid} if every Kakutani factor map $\psi:\Gbf_1/\Gamma_1\to\Gbf_2/\Gamma_2$ to an ergodic nonsplitting $\slfrak_2$-induced unipotent flow $(\Gbf_2/\Gamma_2,\mu_2,U^{(2)}_t)$ is affine up to $U^{(2)}_t$-translations. Moreover, if $\psi$ is a Kakutani equivalence, an affine pair $(g_0,\alpha)\in\Gbf_2'\times\Sur(\Gbf_1',\Gbf_2')$ associated with $\psi$ up to $U^{(2)}_t$-translations satisfies $\alpha(\Gamma_1')=g_0\Gamma_2'g_0^{-1}$ and $\alpha$ is an isomorphism.
\end{Definition}

The following is our main theorem:
\begin{Theorem}\label{thm:main}
    Every ergodic nonsplitting $\slfrak_2$-induced unipotent flow is time change rigid.
\end{Theorem}

\textit{Theorem~\ref{thm:main}} removes both the semisimplicity hypothesis and the assumption in \cite{lindenstrauss2025timechangerigidityunipotent} that base group does not have compact factors. In particular, it includes three geometric phenomena beyond the scope of previous results: irreducible lattices coupling noncompact and compact semisimple factors, semidirect products with non central unipotent radicals, and twist extensions combining various fibers, which are demonstrated as follows: 
\begin{itemize}
    \item $\Gbf=\SL_2(\R)\times\SL_2(\R)\times\SU(2)$, $\Gamma=\SL_2(\Z)\times\Lambda$, where $\Lambda<\SL_2(\R)\times\SU(2)$ is an irreducible lattice, and $U_t=\left(\left(\begin{smallmatrix}
        1 & t\\
        0 & 1
    \end{smallmatrix}\right),\left(\begin{smallmatrix}
        1 & t\\
        0 & 1
    \end{smallmatrix}\right),I_2\right)$. This, in particular, is an example for the case $\Gbf$ is semisimple with compact factors.
    \item $\Gbf=\SL_d(\R)\ltimes\R^d$, $\Gamma<\Gbf$ a lattice, $U_t=\exp(t\bu)$ for some nontrivial nilpotent $\bu\in\slfrak_d(\R)\ltimes\{0\}$. This, in particular, is an example for the case $\Gbf=\Sbf\ltimes\Nbf$, where $\Sbf$ is the identity component of a real semisimple linear algebraic group and $\Nbf$ is a simply connected, connected, nilpotent Lie group.
    \item $\Gbf$ is the subgroup of $\SL_9(\R)$ consisting of matrices of the form
    \[
    \left(\begin{smallmatrix}
        A &  &  &  & x\\
         & B &  &  & 0\\
         &  & C &  & 0\\
         & & & 1 & y\\
         & & & & 1
    \end{smallmatrix}\right),
    \]
    where $A\in\SL_2(\R)$, $B\in\SO(3)$, $C\in\SO(2)$, $x\in\R^2$, and $y\in\R$; $\Gamma<\Gbf$ is the lattice consisting of matrices
    \[
    \left(\begin{smallmatrix}
        A &  &  &  & x\\
         & \rho_1(A) &  &  & 0\\
         &  & \rho_2(A) &  & 0\\
         & & & 1 & n+\rho_3(A)\\
         & & & & 1
    \end{smallmatrix}\right),
    \]
    where $A\in\Gamma_0=\langle\left(\begin{smallmatrix}
        1 & 2\\
        0 & 1
    \end{smallmatrix}\right),\left(\begin{smallmatrix}
        1 & 0\\
        2 & 1
    \end{smallmatrix}\right)\rangle$, $x\in\Z^2$, $n\in\Z$, and
    \[
    \rho_1:\Gamma_0\to\SO(3),\qquad\rho_2:\Gamma_0\to\SO(2),\qquad\rho_3:\Gamma_0\to\R
    \]
    are homomorphisms such that
    \[
    \big\{\big(\rho_1(A),\rho_2(A),\rho_3(A)+\Z\big):A\in\Gamma_0\big\}
    \]
    is dense in $\SO(3)\times\SO(2)\times(\R/\Z)$; and $U_t$ is the flow that has $\left(\begin{smallmatrix}
        1 & t\\
        0 & 1
    \end{smallmatrix}\right)$ as the $\SL_2(\R)$-factor, and the other factors are trivial. This, in particular, is not covered by the previous two classes of examples.
\end{itemize}

As an application of \textit{Theorem~\ref{thm:main}}, we obtain a complete algebraic characterization of loosely Kronecker flows among ergodic $\slfrak_2$-induced unipotent flows on general algebraic quotients, extending \cite{MR4255058}*{Corollary 1.2}:
\begin{Corollary}\label{Coroll:splitting<=>looselyKronecker}
    Let $(\Gbf/\Gamma,\mu,U_t)$ be an ergodic $\slfrak_2$-induced unipotent flow. Then $U_t$ is loosely Kronecker if and only if $U_t$ is splitting.
\end{Corollary}

\textit{Corollary~\ref{Coroll:splitting<=>looselyKronecker}} in particular shows that the nonsplitting assumption in \textit{Theorem~\ref{thm:main}} is sharp. It is also worth mentioning that the proof of \textit{Corollary~\ref{Coroll:splitting<=>looselyKronecker}} does not depend on the computation of Kakutani invariant and hence is different from the semisimple quotient case in \cite{MR4255058}. A more detailed discussion on this is provided at the end of \textit{Section~\ref{Section:ProofOfSplitting<=>LooselyKronecker}}.

\medskip

With \textit{Theorem~\ref{thm:main}} and \textit{Corollary~\ref{Coroll:splitting<=>looselyKronecker}}, we immediately obtain the following dichotomy for ergodic $\slfrak_2$-induced unipotent flows with respect to Kakutani equivalence, as in \cite{lindenstrauss2025timechangerigidityunipotent}:
\begin{Corollary}\label{Coroll:Dichotomy}
    Two ergodic $\slfrak_2$-induced unipotent flows are Kakutani equivalent if and only if at least one of the following holds:
    \begin{enumerate}[label=(\arabic*)]
        \item\label{item:Dichotomy:1} Both flows are splitting.
        \item\label{item:Dichotomy:2} The two flows are measurably isomorphic.
    \end{enumerate}
\end{Corollary}

\subsection*{Main difficulties and overview of the arguments}
The three examples following the statement of \textit{Theorem~\ref{thm:main}} above illustrate new obstacles in our general case beyond \cites{lindenstrauss2025timechangerigidityunipotent}. In particular, the presence of nontrivial (nil)radical and/or compact simple factors can make the left action of $\Hbf_1=\langle\exp(\hfrak_1)\rangle$ on $L^2_0(\Gbf_1/\Gamma_1,\mu_1)$ have no spectral gap, where $\hfrak_1\subseteq\gfrak_1$ is an $\slfrak_2$-subalgebra corresponding to an ergodic $\slfrak_2$-induced unipotent flow $(\Gbf_1/\Gamma_1,\mu_1,U^{(1)}_t)$. To overcome this difficulty, we decompose the base group $\Gbf_1$, and consequently, the quotient space $\Gbf/\Gamma$. Then, base on this decomposition we show that the $\Hbf_1$-action on the subspace of $\Wbf$-invariant functions in $L^2_0(\Gbf_1/\Gamma_1)$ has a spectral gap, where $\Wbf$ is the product of the compact semisimple factor and an abelian factor of $\Gbf_1$ ``compatible" with $\Gamma_1$. And to compensate for the restriction of $\Wbf$-invariance on $L^2_0(\Gbf_1/\Gamma_1)$, we control the behavior of a Kakutani factor map $\psi:\Gbf_1/\Gamma_1\to\Gbf_2/\Gamma_2$ between ergodic $\slfrak_2$-induced unipotent flows $(\Gbf_i/\Gamma_i,\mu_i,U^{(i)}_t)$, $i=1,2$, along $\Wbf$-directions, where $U^{(2)}_t$ is assumed to be nonsplitting. We further explain our ideas in the following.

\medskip

Given an ergodic $\slfrak_2$-induced unipotent flow $(\Gbf_1/\Gamma_1,\mu_1,U^{(1)}_t)$, we obtain a decomposition
\begin{equation}\label{eq:Introduction:DecompositionsOfG_1}
    \Gbf_1=\Sbf_1\Pbf_1\Qbf_1,
\end{equation}
where $\Sbf_1$ is the semisimple factor of $\Gbf_1$, $\Pbf_1$ is a connected normal subgroup of the unipotent radical of $\Gbf_1$ such that $[\Sbf_1,\Pbf_1]=\Pbf_1$, $\Qbf_1$ is a connected subgroup centralizing $\Gbf_1$ such that $\Gamma_{\Qbf_1}=\Gamma_1\cap\Qbf_1$ is a lattice in $\Qbf_1$, and $\Gamma_{1,\nc}=\Gamma_1\cap\Qbf_{1,\nc}$ is a lattice in $\Qbf_{1,\nc}$. Here, $\Qbf_{1,\nc}$ is the noncompact factor of $\Qbf_1$, and $\Qbf_1$ is exactly the abelian factor ``compatible" with $\Gamma_1$ mentioned previously. This decomposition follows from the results in \cites{MR910358} and the theory of linear algebraic groups and their lattices. And with it, we decompose the homogeneous space $\Gbf_1/\Gamma_1$ as a skew-product
\begin{equation}\label{eq:Introduction:DecompositionsOfG_1/Gamma_1}
    \Gbf_1/\Gamma_1=\Sbf_1\Pbf_1/\Gamma_{\Sbf_1\Pbf_1}\ltimes\Qbf_1/\Gamma_{\Qbf_1},    
\end{equation}
where $\Gamma_{\Sbf_1,\Pbf_1}$ is a lattice in $\Sbf_1\Pbf_1$. 

\medskip

Denote the compact semisimple factor in $\Sbf_{1}$ by $\Sbf_{1,\operatorname{c}}$ and $\Wbf=\Sbf_{1,\operatorname{c}}\Qbf_1$, as mentioned previously. With the decomposition \eqref{eq:Introduction:DecompositionsOfG_1/Gamma_1}, we can identify the subspace $L^2_0(\Gbf_1/\Gamma_1)^\Wbf$ of $\Wbf$-invariant functions in $L^2_0(\Gbf_1/\Gamma_1)$ with $L^2_0(\Sbf_1\Pbf_1/\Gamma_{\Sbf_1\Pbf_1})^{\Sbf_{1,\operatorname{c}}}$. Using a decomposition similar to \eqref{eq:Introduction:DecompositionsOfG_1/Gamma_1} for $\Sbf_1\Pbf_1/\Gamma_{\Sbf_1\Pbf_1}$, we further identify this subspace as the direct sum of $\Hcal_2\cong L^2_0(\Sbf_{1,\nc}/\Gamma_{\Sbf_{1,\nc}})$ and the orthogonal complement $\Hcal_2^\perp$. Here, $\Sbf_{1,\nc}$ is the noncompact semisimple factor in $\Sbf_1$ and $\Gamma_{\Sbf_{1,\nc}}<\Sbf_{1,\nc}$ is a lattice. Notice that as $\Hbf_1\subseteq\Sbf_{1,\nc}$ contains no compact simple factors, the $\Hbf_1$-action on $\Hcal_1\cong L^2_0(\Sbf_{1,\nc}/\Gamma_{1,\nc})$ has a spectral gap. And as $\Pbf_1$ is normal in $\Gbf_1$, $\Hcal_2^\perp$ is a subspace of $\Hcal^\perp$-the orthogonal complement of $\Hcal$ in $L^2_0(\Sbf_1\Pbf_1/\Gamma_{\Sbf_1\Pbf_1})$, where $\Hcal$ is the subspace of $\Pbf_1$-invariant functions in $L^2_0(\Sbf_1\Pbf_1/\Gamma_{\Sbf_1\Pbf_1})$. Hence, to show $\Hbf_1$ acts on $L^2_0(\Gbf_1/\Gamma_1)^\Wbf$ with a spectral gap, we only have to show $\Hbf_1$ acts on $\Hcal^\perp$ with a spectral gap. To do so, we argue based on the results in \cites{MR3054613,MR1779896}. More specifically, we show that $(\Dbf,\Ebf)$ has strong relative property (T) via \cites{MR3054613}, where $\Dbf$ is the normal closure of $\Hbf_1$ and $\Ebf=[\Dbf,\Pbf_1]$. With this, we deduce that the trivial representation of $\Dbf$ is not weakly contained in $\Hcal^\perp\vert_\Dbf$, i.e. $\mathbbm{1}_\Dbf\nprec\Hcal^\perp\vert_\Dbf$, by showing $\Hcal^\perp$ has no $\Ebf$-invariant vectors. This is based on the fact that $[\Sbf_1,\Pbf_1]=\Pbf_1$. We then obtain a spectral gap for the $\Hbf_1$-action on $\Hcal^\perp$ via \cites{MR1779896}.

\medskip

Next, we study the behavior of $\psi:\Gbf_1/\Gamma_1\to\Gbf_2/\Gamma_2$ a Kakutani factor map between ergodic $\slfrak_2$-induced unipotent flows $(\Gbf_i/\Gamma_i,\mu_i,U^{(i)}_t)$, $i=1,2$, with respect to $\Wbf$, assuming that $U^{(2)}_t$ is nonsplitting. Recall the decomposition \eqref{eq:Introduction:DecompositionsOfG_1} and the behavior of $\psi$ acting on close points and along $A^{(1)}_t$-orbits established in \cites{lindenstrauss2025timechangerigidityunipotent}, where $\psi$ is adjusted so that its time change satisfy the regularity condition of $\epsilon$-good (cf. \textit{Definition~\ref{Definition:epsilon-goodTimeChange}}) and $\psi(A^{(1)}_1x)$, $A^{(2)}_1\psi(x)$ lie on the same $U^{(2)}_t$-orbit, $A^{(i)}_t=\exp(t\ba^{(i)})$, and $\{\bu^{(i)},\ba^{(i)},\bou^{(i)}\}$ is an $\slfrak_2$-triple corresponding to $U^{(i)}_t$, $i=1,2$. With this, we apply an inductive argument to deduce that
\[
\psi(gx)=U^{(2)}_{f(x,g)}\varpi(g)\psi(x),\qquad g\in\Wbf.
\]
Here, $\varpi$ is a continuous homomorphism from $\Wbf$ to $\Zbf_2=C_{\Gbf_2}(\Hbf_2)$, the centralizer of $\Hbf_2=\big\langle\exp\big(\R\bu^{(2)}\oplus\R\ba^{(2)}\oplus\R\bou^{(2)}\big)\big\rangle$, such that
\[
\varpi(\gamma)\in\Gamma_2\cap Z(\Gbf_2),\qquad\forall\gamma\in\Gamma_{1,\nc}
\]
(technically, we need to lift to the universal cover $\widetilde{\Wbf}$ of $\Wbf$, but we skip this detail in this overview). And $f$ can be considered as a measurable cocycle $\Gbf_1/\Gamma_1\times(\Wbf/\Gamma_{1,\nc})\to\R$. Note that this is not the case for the flow $A^{(1)}_t$, due to the fact that the flow $A^{(1)}_t$ is not periodic, while
\[
\gamma x=x,\qquad\forall x\in\Gbf_1/\Gamma_1,\quad\gamma\in\Gamma_{1,\nc}.
\]
This difference allows us to solve in this case, based on \cite{MR1784202}, the cohomological equation
\begin{equation}\label{eq:Introduction:CohomologicalEquation}
    b(x)-b(gx)=f(x,g),
\end{equation}
upon which we deduce that $\Tilde{\psi}(\cdot)=U^{(2)}_{b(\cdot)}\psi(\cdot)$ is a Kakutani factor map between $(\Gbf_i/\Gamma_i,\mu_i,U^{(i)}_t)$, $i=1,2$, satisfying
\begin{equation}\label{eq:Introduction:Controlw.r.t.Wbf}
    \Tilde{\psi}(gx)=\varpi(g)\Tilde{\psi}(x).
\end{equation}
To make sure that the corresponding time change $\Tilde{\tau}$ of $\Tilde{\psi}$ is still $\epsilon$-good, we use the convolution technique twice: once in the process of solving \eqref{eq:Introduction:CohomologicalEquation} and once after obtaining the solution of \eqref{eq:Introduction:CohomologicalEquation}.

\medskip

Now, with both the existence of a spectral gap for the $\Hbf_1$-action on $L^2_0(\Gbf_1/\Gamma_1)^\Wbf$ and the behavior of $\psi$ with respect to $\Wbf$ via \eqref{eq:Introduction:Controlw.r.t.Wbf}, we move to the renormalization procedure, establishing the limit $\phi$ of $\psi_n=A^{(2)}_{-n}\psi A^{(1)}_n$ as $n$ runs on an increasing subsequence in $\N$ of large proportion. Here, we construct certain ``good" sets via the \textit{Pointwise ergodic theorem} and \textit{Lusin's theorem}, and fill their $\Wbf$-orbits, making them $\Wbf$-invariant. This allows us to apply an $\SL_2(\R)$-ergodic theorem (cf. \textit{Proposition~\ref{Prop:SpectralGapCondition:SL_2-ergodicTheorem}}), which stems from the fact that $\Hbf_1$-acts on $L^2_0(\Gbf_1/\Gamma_1)^\Wbf$ with a spectral gap, while the behavior \eqref{eq:Introduction:Controlw.r.t.Wbf} of $\psi$ with respect to $\Wbf$ ensures the ``good" sets, after $\Wbf$-filling, are still ``good".

\medskip

In the end, we finish the proof of \textit{Theorem~\ref{thm:main}} via \textit{Ratner's Measure Classification Theorem} and the rigidity of unipotent flows under measurable factor maps.

\subsection*{Organization of the paper} The structure of our paper is as follows.

\begin{itemize}
    \item In \textit{Section~\ref{Section:ProofOfSplitting<=>LooselyKronecker}}, assuming \textit{Theorem~\ref{thm:main}}, we prove \textit{Corollary~\ref{Coroll:splitting<=>looselyKronecker}} by applying classical results in \cites{katok1977monotone,MR653094} and \textit{Theorem~\ref{thm:main}} on unipotent flows $(\SL_{2n+1}(\R)/\Gamma_1,\mu_1,U^{(1)}_t)$. We then deduce \textit{Corollary~\ref{Coroll:Dichotomy}}.
    \item In \textit{Section~\ref{Section:Preliminaries}}, we recall some necessary notions and properties of ergodic (nonsplitting) $\slfrak_2$-induced unipotent flows from \cite{lindenstrauss2025timechangerigidityunipotent}*{Sections~3-8}. 
    \item In \textit{Section~\ref{Section:StructureOfErgodicUnipotentFlow}}, given an ergodic $\slfrak_2$-induced unipotent flow $(\Gbf/\Gamma,\mu,U_t)$, based on \cite{MR910358}, we deduce a decomposition of the base group $\Gbf$, and consequently a decomposition of the quotient $\Gbf/\Gamma$.
    \item In \textit{Section~\ref{Section:SpectralGapCondition}}, with the settings in \textit{Section~\ref{Section:StructureOfErgodicUnipotentFlow}}, we extend the \textit{$\SL_2(\R)$-ergodic theorem} \cite{lindenstrauss2025timechangerigidityunipotent}*{Corollary 3.15} to our case by showing that it holds for $\Sbf_{\operatorname{c}}\Qbf$-invariant functions, where $\Sbf_{\operatorname{c}}$ is the compact factor of $\Sbf$ (cf. \textit{Proposition} \ref{Prop:SpectralGapCondition:SL_2-ergodicTheorem}). 
    \item In \textit{Section~\ref{Section:BehaviorOfpsiOnClosePointsAndAlongA_t}}, we recall from \cite{lindenstrauss2025timechangerigidityunipotent}*{Sections~8 and~10} the behavior of a Kakutani factor map $\psi:\Gbf_1/\Gamma_1\to\Gbf_2/\Gamma_2$ between ergodic $\slfrak_2$-induced unipotent flows $(\Gbf_i/\Gamma_i,\mu_i,U^{(i)}_t)$, $i=1,2$, acting on close points and along the flow $A^{(i)}_t$, $i=1,2$, where $U^{(2)}_t$ is assumed to be nonsplitting. 
    \item  In \textit{Section~\ref{Section:BehaviorOfpsiW_t}}, with the Kakutani factor map $\psi$ obtained from \textit{Section~\ref{Section:BehaviorOfpsiOnClosePointsAndAlongA_t}}, we establish the behavior of $\psi$ with respect to $\Wbf=\Sbf_{1,\operatorname{c}}\Qbf_1$, where $\Sbf_{1,\operatorname{c}}$ is the compact semisimple factor and $\Qbf_1$ is the abelian factor $\Qbf_1$ constructed in \textit{Section~\ref{Section:StructureOfErgodicUnipotentFlow}}. More specifically, we show that up to $U^{(2)}_t$-translations, $\psi$ is $\Wbf$-equivariant via a continuous homomorphism $\Tilde{\varpi}:\widetilde{\Wbf}\to\Zbf_2=C_{\Gbf_2}(\Hbf_2)$, where $\widetilde{\Wbf}$ is a universal cover of $\Wbf$.
    \item In \textit{Section~\ref{Section:Renormalization}}, with the Kakutani factor map $\psi$ obtained from \textit{Sections~\ref{Section:BehaviorOfpsiOnClosePointsAndAlongA_t}} and \textit{\ref{Section:BehaviorOfpsiW_t}}, we sketch the renormalization arguments in \cite{lindenstrauss2025timechangerigidityunipotent}*{Section~11}, and identify the technical details we need to adjust in our settings. This results in \textit{Theorem~\ref{Theorem:Renormalization:phi=lim_(n)psi_n}}, which establishes the limit $\phi$ of $A^{(2)}_{-n}\psi A^{(1)}_n$ as $n$ runs on a ``large" subset of $\N$, and shows $\phi U^{(1)}_t=U^{(2)}_t\phi$ for all $t\in\Q$.
    \item In \textit{Section~\ref{Section:ProofOfTheorem:MainTheorem}}, we conclude with the proof of \textit{Theorem~\ref{thm:main}} using results in previous sections, \textit{Ratner's Measure Classification Theorem}, and the rigidity of unipotent flows under measurable factor maps.
\end{itemize}

\medskip

The following graph explains the relations between \textit{Sections~\ref{Section:StructureOfErgodicUnipotentFlow}-\ref{Section:ProofOfTheorem:MainTheorem}}, while \textit{Section~\ref{Section:Preliminaries}} is used throughout \textit{Sections~\ref{Section:StructureOfErgodicUnipotentFlow}-\ref{Section:ProofOfTheorem:MainTheorem}}.
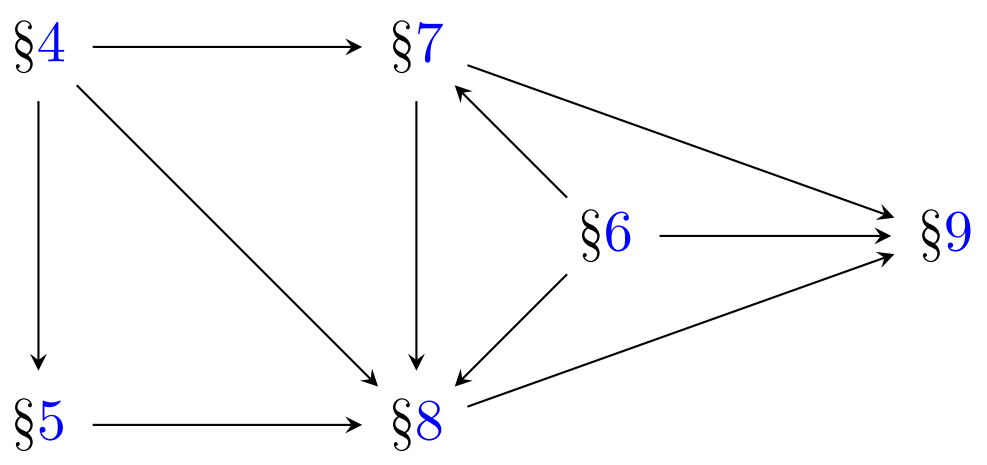
\begin{figure}[H]
    \centering
    \scalebox{1}
    {
    \begin{tikzpicture}[scale=5]
	\tikzstyle{vertex}=[circle,minimum size=20pt,inner sep=0pt]
	\tikzstyle{selected vertex} = [vertex, fill=red!24]
	\tikzstyle{edge1} = [draw,line width=5pt,-,red!50]
    \tikzstyle{edge2} = [draw,line width=5pt,-,green!50]
    \tikzstyle{edge3} = [draw,line width=5pt,-,blue!50]
    \tikzstyle{edge4} = [draw,line width=5pt,-,brown!50]

    \tikzstyle{edge} = [draw,thick,-,black]
        \node[vertex] (v4) at (0,0.5) {\S\ref{Section:StructureOfErgodicUnipotentFlow}};

        \node[vertex] (v5) at (0,0)
        {\S\ref{Section:SpectralGapCondition}};

        \node[vertex] (v6) at (0.75,0.25)
        {\S\ref{Section:BehaviorOfpsiOnClosePointsAndAlongA_t}};

        \node[vertex] (v7) at (0.5,0.5)
        {\S\ref{Section:BehaviorOfpsiW_t}};

        \node[vertex] (v8) at (0.5,0)
        {\S\ref{Section:Renormalization}};

        \node[vertex] (v9) at (1.2,0.25)
        {\S\ref{Section:ProofOfTheorem:MainTheorem}};

        \draw[-stealth] (v4)--(v5);
        \draw[-stealth] (v4)--(v7);
        \draw[-stealth] (v4)--(v8);

        \draw[-stealth] (v5)--(v8);

        \draw[-stealth] (v6)--(v7);
        \draw[-stealth] (v6)--(v8);
        \draw[-stealth] (v6)--(v9);

        \draw[-stealth] (v7)--(v8);
        \draw[-stealth] (v7)--(v9);

        \draw[-stealth] (v8)--(v9);
    \end{tikzpicture}
    }
\caption{Structure of paper}\label{fig:StructureMap}
\end{figure}

\subsection*{Convention} On the notations in this paper, we use $\epsilon_*$, $\delta_*$, $\omega_*$, $C_*$, $R_*$, $\mathcal{K}_*$, $\mathcal{Z}_*$, $*\in\N$, to denote constants and sets that are defined uniquely for the entire paper. Any other notations for constants and sets are local to the proof in which they appear.

\subsection*{Acknowledgement} The authors are grateful for helpful discussions with Adam Kanigowski and Elon Lindenstrauss.

\subsection*{AI disclosure} AI tools were used for literature searches. All mathematical ideas and contributions are due to the authors.

\section{Proof of \textit{Corollary~\ref{Coroll:splitting<=>looselyKronecker}} and \textit{Corollary~\ref{Coroll:Dichotomy}}}\label{Section:ProofOfSplitting<=>LooselyKronecker}
In this section, we first prove \textit{Corollary~\ref{Coroll:splitting<=>looselyKronecker}} and \textit{Corollary~\ref{Coroll:Dichotomy}} assuming \textit{Theorem~\ref{thm:main}}, and end this section with some additional discussion on \textit{Definition~\ref{def:TimeChangeRigid}}.

\subsection{Proof of Corollary~\ref{Coroll:splitting<=>looselyKronecker}}
\subsubsection{Splitting implies loosely Kronecker}\label{subSection:Splitting<=>LooselyKronecker:=>}
Let $(\Gbf/\Gamma,\mu,U_t)$ be an ergodic splitting $\slfrak_2$-induced unipotent flow. Following the proof of the fact that horocycle flows are loosely Kronecker by Ratner \cite{MR516248} with suitable modifications, we deduce that $U_t$ is loosely Kronecker. For more details, we refer to \cite{kanigowski2018horocycleflownegativevariable}.

Alternatively, we can also argue using the characterization of loosely Kronecker unipotent flows on semisimple quotients \cite{MR4255058}*{Corollary 1.2} and the fact that a compact extension of a loosely Kronecker flow (that is ergodic) is also loosely Kronecker \cite{MR653094}*{Section 7, Theorem 7.3}. We will return to this point in \textit{Remark \ref{Remark:StructuralProperties:SplittingImpliesLooselyKronecker}}, after providing in \textit{Section~\ref{Section:StructureOfErgodicUnipotentFlow}} some structural properties of ergodic $\slfrak_2$-induced unipotent flows.

\subsubsection{Loosely Kronecker implies splitting}\label{subSection:Splitting<=>LooselyKronecker:<=}
Let $(\Gbf/\Gamma,\mu,U_t)$ be a loosely Kronecker $\slfrak_2$-induced unipotent flow. Assume that $U_t$ is nonsplitting, we use \textit{Theorem~\ref{thm:main}} to reach a contradiction, hence deducing that $U_t$ is splitting.

Consider a (nontrivial) unipotent flow $U^{(1)}_t$ on $(\SL_{2n+1}(\R)/\Gamma_1,\mu_1)$, where $n\in\N$, $\Gamma_1<\SL_{2n+1}(\R)$ is a lattice, and $\mu_1$ is the Haar probability measure on $\SL_{2n+1}(\R)/\Gamma_1$. It easily follows from the \textit{Jacobson-Morozov theorem} (cf. \cite{MR1920389}*{Theorem 10.3}), \textit{Mautner phenomenon} (cf. \cite{MR586875}*{Theorem 1.1}), and the fact that $\slfrak_{2n+1}(\R)$ is simple, that $U^{(1)}_t$ is an ergodic non-splitting $\slfrak_2$-induced unipotent flow on $(\SL_{2n+1}(\R)/\Gamma_1,\mu_1)$. Hence, \textit{Theorem~\ref{thm:main}} implies that $U^{(1)}_t$ is time change rigid.

Since $(\Gbf/\Gamma,\mu,U_t)$ is loosely Kronecker and $(\SL_{2n+1}(\R)/\Gamma_1,\mu_1,U^{(1)}_t)$ is ergodic, it follows from \cite{katok1977monotone}*{Theorem 1}, \cite{MR653094}*{Theorems 4.1, 6.7} that there exists a Kakutani factor map $\psi:\SL_{2n+1}(\R)/\Gamma_1\to\Gbf/\Gamma$ between the two flows. Since $U^{(1)}_t$ is time change rigid and $(\Gbf/\Gamma,\mu,U_t)$ is an ergodic nonsplitting $\slfrak_2$-induced unipotent flow, it follows from \textit{Definition \ref{def:TimeChangeRigid}} that there exists  a continuous surjective homomorphism,
\[
\alpha:\SL_{2n+1}(\R)\to\Gbf/\Core_\Gbf(\Gamma).
\] 
But as $\SL_{2n+1}(\R)$ has no proper nontrivial normal subgroup, $\alpha$ is in fact an isomorphism, hence
\[
\SL_{2n+1}(\R)\cong\Gbf/\Core_\Gbf(\Gamma).
\]
This gives us a contradiction, since $n\geq1$ is arbitrary. 

\medskip

Thus, we finish the proof of \textit{Corollary~\ref{Coroll:splitting<=>looselyKronecker}}.\qed

\subsection{Proof of Corollary~\ref{Coroll:Dichotomy}}
If \ref{item:Dichotomy:1} holds, then both ergodic $\slfrak_2$-induced unipotent flows are splitting. Hence, both are loosely Kronecker by \textit{Corollary~\ref{Coroll:splitting<=>looselyKronecker}}, and thus Kakutani equivalent. If \ref{item:Dichotomy:2} holds then obviously the two flows are Kakutani equivalent. Conversely, given two Kakutani equivalent ergodic $\slfrak_2$-induced unipotent flows, if one of them is splitting, then it is loosely Kronecker by \textit{Corollary~\ref{Coroll:splitting<=>looselyKronecker}}. Hence, so is the other one, and thus \ref{item:Dichotomy:1} holds. And if both flows are nonsplitting, then by \textit{Theorem~\ref{thm:main}} and \textit{Definition~\ref{def:TimeChangeRigid}}, \ref{item:Dichotomy:2} holds.\qed

\subsection{Remarks on time change rigidity and the (non)splitting condition}\label{subSection:RemarkOnTimeChangeRigidityAndSplitting}
We end this section with two observations based on \textit{Corollary~\ref{Coroll:splitting<=>looselyKronecker}}, in addition to the comments following our definition of (non)splitting $\slfrak_2$-induced unipotent flow in \textit{Section~\ref{Section:Introduction}}, on why we assume the $\slfrak_2$-induced unipotent flows in \textit{Definition~\ref{def:TimeChangeRigid}} to be nonsplitting.

\medskip

In the following, $(\Gbf_i/\Gamma_i,\mu_i,U^{(i)}_t)$, $i=1,2$, are ergodic $\slfrak_2$-induced unipotent flows, and \textit{Corollary~\ref{Coroll:splitting<=>looselyKronecker}} is applied without mention.

\medskip

In \textit{Definition~\ref{def:TimeChangeRigid}}, given $(\Gbf_1/\Gamma_1,\mu_1,U^{(1)}_t)$ nonsplitting, i.e. not loosely Kronecker, we do not consider Kakutani factor maps $\psi:\Gbf_1/\Gamma_1\to\Gbf_2/\Gamma_2$ where $(\Gbf_2/\Gamma_2,\mu_2,U^{(2)}_t)$ is splitting, i.e. loosely Kronecker. This is because then $\Gbf_2$ may be taken to have an arbitrary large dimension, which forbids $\psi$ to be affine. Note that such $\psi$ always exists, as explained in \textit{Subsection~\ref{subSection:Splitting<=>LooselyKronecker:<=}}.

\medskip
    
We do not consider in \textit{Definition \ref{def:TimeChangeRigid}} the notion of time change rigidity for $(\Gbf_1/\Gamma_1,\mu_1,U^{(1)}_t)$ splitting, i.e. loosely Kronecker. This is because classical results \cite{katok1977monotone}*{Theorem 1}, \cite{MR653094}*{Corollary 6.8} imply that for any $(\Gbf_2/\Gamma_2,\mu_2,U^{(2)}_t)$ nonsplitting, i.e. not loosely Kronecker, there is no Kakutani factor map $\psi:\Gbf_1/\Gamma_1\to\Gbf_2/\Gamma_2$ between $(\Gbf_i/\Gamma_i,\mu_i,U^{(i)}_t)$, $i=1,2$.

\section{Basic notions and properties of \texorpdfstring{$\slfrak_2$}{}-induced unipotent flows}\label{Section:Preliminaries}
In this section, we recall from \cite{lindenstrauss2025timechangerigidityunipotent} the basic notations and necessary preliminaries, along with some auxiliary lemmas.

\subsection{$(\delta,\epsilon,R)$-two sides matching}\label{subSection:Preliminaries:2SidesMatchings}
We begin by recalling the notion of two sides matching in \cite{lindenstrauss2025timechangerigidityunipotent}*{Definition 3.1} for describing the relation between orbits of two points along a flow:
\begin{Definition}\label{Def:Preliminaries:TwoSidesMatching}
    Let $(X,d)$ be a metric space and $T_t$ a flow on $X.$ Given $\delta,\epsilon\in(0,1)$ and $R>1$, $x,y\in X$ are called $(\delta,\epsilon,R)$-two sides matchable if there are closed subsets $A,A'$ of $[-R,R]$ with $\absolute{A},\absolute{A'}>(1-\epsilon)2R$ and an increasing absolutely continuous map $h:A\to A'$ such that $0\in A,h(0)=0$, and for $t\in A,d(T_{h(t)}x,T_ty)<\delta$ and $\absolute{h'(t)-1}<\epsilon.$
\end{Definition}
Here, and for the entire paper, we denote by $\absolute{A}$ the Lebesgue measure of $A\subseteq\R$ whenever $A$ is Lebesgue-measurable. It follows from \cite{lindenstrauss2025timechangerigidityunipotent}*{Lemma 3.2} that the function $h$ above can be taken to be a $C^1$-smooth function $h:[-R,R]\to[-R,R]$ with $\absolute{h'(t)-1}<\epsilon$ for all $\absolute{t}\leq R$, given that $x,y$ are $(\delta',\epsilon',R)$-two sides matchable for some sufficiently small $\delta',\epsilon'>0$ depending on $\delta,\epsilon$. Hence, we will consider the matching function $h$ in \textit{Definition~\ref{Def:Preliminaries:TwoSidesMatching}} to be $C^1$ from $[-R,R]$ to $[-R,R]$.

\subsection{$\slfrak_2$-induced unipotent flow}\label{subSection:Preliminaries:HorocyclicUnipotentFlow}
Consider an $\slfrak_2$-induced unipotent flow $(\Gbf/\Gamma,\mu,U_t)$. In the Lie algebra $\hfrak\cong\slfrak_2(\R)$ corresponding to $U_t=\exp(t\bu)$, let $\ba,\bou\in\hfrak$ be such that $\{\bu,\ba,\bou\}$ is an $\slfrak_2$-triple corresponding to $U_t$: 
\[
[\ba,\bu]=2\bu,\qquad [\ba,\bou]=-2\bou,\qquad [\bu,\bou]=\ba.
\]
For $t\in\R$, denote $A_t=\exp(t\ba)$ and $\Uh_t=\exp(t\bou)$.

Since $\hfrak\cong\slfrak_2(\R)$ and $\Gbf$ is the identity component of a real linear algebraic group, $\mathbf{H}=\langle\exp(\hfrak)\rangle$ is a connected closed subgroup of $\Gbf$ isomorphic to either $\SL_2(\R)$ or $\PSL_2(\R)$, and the homomorphism $\SL_2(\R)\to\mathbf{H}$ given by
\begin{equation}\label{eq:Preliminaries:IsomorphismSL_2toH}
    \Big(\begin{smallmatrix}
    1 & t\\
    0 & 1
\end{smallmatrix}\Big)\mapsto U_t,\quad\Big(\begin{smallmatrix}
    e^t & 0\\
    0 & e^{-t}
\end{smallmatrix}\Big)\mapsto A_t,\quad\Big(\begin{smallmatrix}
    1 & 0\\
    t & 1
\end{smallmatrix}\Big)\mapsto\Uh_t,\quad t\in\R,
\end{equation}
is either an isomorphism or a $2$-to-$1$ covering. Hence, one can directly check the following lemma:
\begin{Lemma}\label{Lemma:Preliminaries:SL_2Computation}

    For $z,p\in\R$, we have $ U_{\phi(t)}\Uh_zA_pU_{-t}=\Uh_{z(t)}A_{p(t)}$, where 
    \[
    \phi(t)=\frac{te^p}{e^{-p}-ze^pt},\quad z(t)=ze^p(e^{-p}-ze^pt),\quad p(t)=-\log(e^{-p}-ze^pt),
    \]
    for all $t\in\R$ such that $e^{-p}-ze^pt>0.$ Moreover, if $\absolute{p}<10^{-2}$ then for all $\absolute{t}\leq10^{-2}\absolute{z}^{-1}$,
    \[
    \absolute{z(t)}\leq2\absolute{z},\qquad\absolute{p(t)}\leq\absolute{p}+2\absolute{t}\absolute{z}.
    \]
    And if $\absolute{p}<\epsilon$ for some $\epsilon\in(0,10^{-2}]$ then $\absolute{\phi'(t)-1}<5\epsilon$ for all $\absolute{t}\leq\epsilon\absolute{z}^{-1}$.
    
\end{Lemma}

The following lemma, which is used at several points in the paper without mentioning, follows from the \textit{Mautner phenomenon} (cf. \cite{MR586875}*{Theorem~1.1}) and the ergodicity of $U_t$:
\begin{Lemma}\label{Lemma:Preliminaries:ErgodicityOfU_t,A_t}

   The flow $A_t$ on $(\mathbf{G}/\Gamma,\mu)$ is ergodic. Moreover, for every $t_0\in\R\setminus\{0\}$, the maps $U_{t_0},A_{t_0}$ on $(\mathbf{G}/\Gamma,\mu)$ are also ergodic.
    
\end{Lemma}

\subsection{Chain basis and metric}\label{subSection:Preliminaries:ChainBasisAndMetric}
Let $\{\bu,\ba,\bou\}$ be an $\slfrak_2$-triple in $\hfrak\subseteq\gfrak$ as before. Consider the representation of $\hfrak\cong\slfrak_2(\R)$ on $\gfrak$ via the adjoint representation, it follows that $\gfrak$ decomposes as a direct sum of finite dimensional irreducible representations of $\hfrak\cong\slfrak_2(\R)$, which gives us a basis for $\gfrak$:
\begin{equation}\label{eq:Preliminaries:ChainBasis}
\{\bu,\ba,\bou,\bx_{0,1},\dots,\bx_{q_1,1},\dots,\bx_{0,n},\dots,\bx_{q_n,n}\},
\end{equation}
where for $1\leq j\leq n$ and $0\leq i\leq q_j$, by setting $\bx_{-1,j}=\bx_{q_j+1,j}=0$, we have
\begin{equation*}
    \begin{aligned}
        &\ad(\bu)\bx_{i,j}=\bx_{i-1,j},\\
        &\ad(\ba)\bx_{i,j}=(q_j-2i)\bx_{i,j},\\
        &\ad(\bou)\bx_{i,j}=(i+1)(q_j-i)\bx_{i+1,j}.
    \end{aligned}
\end{equation*}
We call this \textbf{a chain basis for $\gfrak$ with respect to $\{\bu,\ba,\bou\}$,} or just a \textbf{chain basis} when it is clear what the $\slfrak_2$-triple is.

With a fixed chain basis as in (\ref{eq:Preliminaries:ChainBasis}), we can define a right-invariant Riemannian metric $d_{\mathbf{G}}$ on $\mathbf{G}$ by equipping on $\gfrak$ the inner product $\innproduct{\cdot,\cdot}$ making this basis orthonormal. For $W\subseteq \mathbf{G}$ and $r>0$, we denote 
\[
B^W_r=\left\{g\in W:d_{\mathbf{G}}(g,\ebf_\Gbf)<r\right\}.
\]

The right-invariant metric $d_{\mathbf{G}}$ induces a metric on $\mathbf{G}/\Gamma$ by 
\[
\textstyle{d_{\mathbf{G}/\Gamma}(g\Gamma,h\Gamma)=\inf_{\gamma\in\Gamma}d_{\mathbf{G}}(g,h\gamma).}
\]

For $x\in \mathbf{G}/\Gamma$ and $r>0$, denote
\[
B^{\mathbf{G}/\Gamma}_r(x)=\{y\in \mathbf{G}/\Gamma:d_{\mathbf{G}/\Gamma}(x,y)<r\}.
\]

With $r>0$ small enough, we have an isometry $B^{\mathbf{G}}_r\to B^{\mathbf{G}/\Gamma}_r(x):g\mapsto gx.$ Moreover, if $K\subset \mathbf{G}/\Gamma$ is a compact set, then there exists $r>0$ such that we have the above isometry for all $x\in K.$ Define the injective radius of $K$, denote by $\inj(K)$, to be half the supremum of such $r>0.$ 

As mentioned in \textit{Section~\ref{Section:Introduction}}, we always consider $\Gbf$ as the identity component of a real linear algebraic group. We fix a closed embedding of $\Gbf$ in $\SL_N(\R)$ for some $N\in\N$, and for $g\in \mathbf{G}$, we denote $\norm{g}$ as the Hilbert-Schmidt norm of $g$.

\subsection{Kakutani-Bowen balls}\label{subSection:Preliminaries:KakBall}
With the notations in \textit{Subsections~\ref{subSection:Preliminaries:HorocyclicUnipotentFlow},~\ref{subSection:Preliminaries:ChainBasisAndMetric}}, we denote the centralizer of $\mathbf{H}=\langle\exp(\hfrak)\rangle$ in $\Gbf$ as
\[
\mathbf{Z}=C_{\Gbf}(\mathbf{H}).
\]
Its Lie algebra can be characterized by the chain basis defined in \eqref{eq:Preliminaries:ChainBasis}:
\[
\textstyle{\zfrak=\bigoplus_{j\in J}\R\bx_{0,j}},\qquad J=\{1\leq j\leq n:q_j=0\}.
\]
To simplify the notation, we denote the complementary part of $\zfrak\oplus\hfrak$ as
\[
\textstyle{\mathfrak{tr}=\bigoplus_{j\notin J}\bigoplus_{i=0}^{q_j}\R\bx_{i,j}.}
\]

For $g\in \mathbf{G}$ sufficiently close to $\ebf_\Gbf$, we have the following decomposition:
\begin{equation}\label{eq:Preliminaries:FormulaOfElementsCloseToIdentity}
    \begin{aligned}
        g=g^{\hfrak}g^{\zfrak}g^{\tr},
 \end{aligned}
\end{equation}
where $g^{\hfrak}\in\mathbf{H}$, $g^{\zfrak}\in\exp(\zfrak)$ and $g^{\tr}\in\exp(\mathfrak{tr})$ so that
\begin{equation*}
    \begin{aligned}
        g^{\hfrak}=U_{c_{\bu}}\Uh_{c_{\bou}}A_{c_{\ba}},\,
        g^{\zfrak}=\exp\left(\textstyle{\sum_{j\in J}c_{0,j}\bx_{0,j}}\right),\, g^{\tr}=\exp\left(\textstyle{\sum_{j\notin J}\sum_{i=0}^{q_j}c_{i,j}\bx_{i,j}}\right).
    \end{aligned}
\end{equation*}

With this, a special type of neighborhoods of a point in $\mathbf{G}$ or $\mathbf{G}/\Gamma$, called Kakutani-Bowen balls, was defined in \cite{lindenstrauss2025timechangerigidityunipotent}*{Definition 3.9}. We recall it in the following:  
\begin{Definition}[Kakutani-Bowen balls]\label{Def:Preliminaries:Kakutani-BowenBalls}

    For $\delta\in(0,1),1<R<\infty,\xt\in \mathbf{G}$, denote by $\Kak(R,\delta,\xt)$ the set of all $\yt\in \mathbf{G}$ such that $g=\yt\xt^{-1}$ can be written as in (\ref{eq:Preliminaries:FormulaOfElementsCloseToIdentity}) with $\absolute{c_{\bu}}\leq\delta,\absolute{c_{\ba}}\leq\delta,\absolute{c_{\bou}}\leq\frac{\delta}{R}$, and $g^{\zfrak}g^{\tr}$ satisfying $d_{\mathbf{G}}(U_tg^{\zfrak}g^{\tr}U_{-t},\ebf_\Gbf)\leq\delta$ whenever $\absolute{t}\leq R.$
    
    In the case $R=\infty$, we denote
    \[
    \textstyle{\Kak(\infty,\delta,\Tilde{x})=\left\{U_{c_{\bu}}A_{c_{\ba}}g\Tilde{x}:\absolute{c_{\bu}},\absolute{c_{\ba}}\leq\delta;g\in\overline{B}^{\Gbf}_\delta\cap\mathcal{O},\, \log(g)\in\bigoplus_{j=1}^{n}\R\bx_{0,j}\right\}}.
    \]
    where $\mathcal{O}$ is a neighborhood of the identity such that $\log:\Ocal\to\gfrak$ is a diffeomorphism onto its image in $\gfrak$. 
    
    If $\xt=\ebf_\Gbf$, we write $\Kak(R,\delta)$ instead of $\Kak(R,\delta,\ebf_\Gbf)$. And for $x\in \mathbf{G}/\Gamma$, we denote
    \[
    \Kak(R,\delta,x)=\Kak(R,\delta)x=\{gx:g\in\Kak(R,\delta)\}.
    \]
    
\end{Definition}

\medskip

We have the following control over Kakutani-Bowen balls under conjugation by $A_{-p}:$
\begin{Lemma}\label{Lemma:Preliminaries:a_(-p)ga_pForgInKak(R,delta)}

    There exists $\denew\label{delta:Preliminaries:a_(-p)ga_pForgInKak(R,delta)}\in(0,10^{-3})$ such that the following holds. For $\delta\in(0,\deold{delta:Preliminaries:a_(-p)ga_pForgInKak(R,delta)}),R\geq1$, $c>0$, and $g\in\Kak(R,\delta)$, let $p=\frac{\log(cR)}{2}$ and $\overline{g}=A_{-p}gA_{p}$, write $\overline{g}=U_{\overline{c}_{\bu}}\Uh_{\overline{c}_{\bou}}A_{\overline{c}_{\ba}}\overline{g}^{\zfrak}\overline{g}^{\tr}$ as in (\ref{eq:Preliminaries:FormulaOfElementsCloseToIdentity}). Then
    \begin{equation*}
    \begin{gathered}
         \absolute{\overline{c}_{\bu}}\leq\delta/(c R),\qquad\absolute{\overline{c}_{\bou}}\leq\delta c,\qquad\absolute{\overline{c}_{\ba}}\leq\delta,\\
        d_{\mathbf{G}}\left(\overline{g}^{\zfrak},\ebf_\Gbf\right)<2\delta,\qquad d_{\mathbf{G}}\left(\overline{g}^{\tr},\ebf_\Gbf\right)\leq\Cold{C:Preliminaries:a_(-p)ga_pForgInKak(R,delta)}\delta\tilde{c}^{(\dim\Gbf)/2}R^{-0.5},
    \end{gathered}
    \end{equation*}
    where $\Cnew\label{C:Preliminaries:a_(-p)ga_pForgInKak(R,delta)}>1$ depends only on $\dim\Gbf$ and $\tilde{c}=\max\{c,c^{-1}\}.$
    
\end{Lemma}
This lemma is adjusted from \cite{lindenstrauss2025timechangerigidityunipotent}*{Lemma 3.12}, and the proof follows in a similar way, as it only relies on the above notion of Kakutani-Bowen balls, the relations between $\bu,\ba,\bou$ with each other and with other elements in the chain basis \eqref{eq:Preliminaries:ChainBasis}, and the fact that norms on the space of polynomials with degree at most $n$ (for some fixed $n\in\N$) are equivalent (for more details of this fact, see \cite{MR3994576}*{Lemma 2.7}).

\subsection{Auxiliary lemmas}\label{subSection:MeasureEstimateLemma}

Given $(\Gbf/\Gamma,\mu,U_t)$ an ergodic $\slfrak_2$-induced unipotent flow, with a corresponding $\slfrak_2(\R)$-subalgebra $\hfrak\subseteq\gfrak$, $\{\bu,\ba,\bou\}$ an $\slfrak_2$-triple, $\Hbf=\langle\exp(\hfrak)\rangle$, and a chain basis and metrics as in \textit{Subsection~\ref{subSection:Preliminaries:ChainBasisAndMetric}}, we provide some auxiliary lemmas from \cite{lindenstrauss2025timechangerigidityunipotent}*{Section~6,~7} in our general case, i.e. without the semisimple condition on $\Gbf$. These are only indirectly used in our paper, where certain results and arguments in \cite{lindenstrauss2025timechangerigidityunipotent} carry over to our case. There, the lemmas in this subsection serve as replacements of their corresponding semisimple versions in \cite{lindenstrauss2025timechangerigidityunipotent}.

\medskip

We begin with a Chevalley representation of the normalizer $\Lbf=N_\Gbf(\Hbf)$ of $\Hbf$ in $\Gbf$, i.e. a rational algebraic representation of $\mathbf{G}$, where $\Lbf$ is the set of elements fixing a vector:

\begin{Lemma}\label{Lemma:MeasureEstimateLemma:ChevalleyRepresentation}

    Let $\rho=\bigwedge^3\Ad:\mathbf{G}\to\GL\left(\bigwedge^3\gfrak\right)$ be the representation of $\mathbf{G}$ given for $g\in \mathbf{G}$ by
    $$\rho(g)(X\wedge Y\wedge Z)=(\Ad(g)X)\wedge(\Ad(g)Y)\wedge(\Ad(g)Z),\forall X,Y,Z\in\gfrak.$$
    Then $\mathbf{L}=\{g\in \mathbf{G}:\rho(g)v_\mathbf{L}=v_\mathbf{L}\}$, where $v_\mathbf{L}=\bu\wedge\ba\wedge\bou\in\bigwedge^3\gfrak.$
    
\end{Lemma}
\begin{proof}
Since $\mathbf{H}=\langle\exp(\hfrak)\rangle$, we get
\[
\Lbf=N_{\mathbf{G}}(\mathbf{H})=\{g\in \mathbf{G}:g\mathbf{H}g^{-1}=\mathbf{H}\}=\{g\in \mathbf{G}:\Ad(g)\hfrak=\hfrak\}.
\] 
This and the fact $\hfrak=\R\bu\oplus\R\ba\oplus\R\bou$ imply that
\[
\mathbf{L}=\{g\in \mathbf{G}:\rho(g)(\R v_\mathbf{L})=\R v_\mathbf{L}\},
\]
and so the inclusion $\{g\in \mathbf{G}:\rho(g)v_\mathbf{L}=v_\mathbf{L}\}\subseteq \mathbf{L}$ is obvious. Hence, we are left to show that if $g\in \mathbf{L}$ then $\rho(g)v_\mathbf{L}=v_\mathbf{L}.$

Fix $g\in \mathbf{L}.$ As $\Ad(g)\hfrak=\hfrak$, we get $\phi=\Ad(g)\vert_{\hfrak}\in\GL(\hfrak)$, and
\begin{equation}\label{eq:MeasureEstimateLemma:ChevalleyRepresentation:rho(g)v_N=det(phi)v_N}
    \begin{aligned}
        \rho(g)v_\mathbf{L}&=\Ad(g)\bu\wedge\Ad(g)\ba\wedge\Ad(g)\bou=\phi(\bu)\wedge\phi(\ba)\wedge\phi(\bou)\\
        &=\det(\phi)\bu\wedge\ba\wedge\bou=\det(\phi)v_\mathbf{L}.
    \end{aligned}
\end{equation}
Consider $\omega:\hfrak^3\to\R:\omega(X,Y,Z)=\kappa([X,Y],Z)$, where $\kappa$ is the Killing form on $\hfrak$, and $[\cdot,\cdot]$ is the Lie bracket. It is obvious that $\omega$ is $\R$-multilinear. Since $\kappa$ is ad-invariant and nondegenerate, we deduce that
\begin{equation}\label{eq:MeasureEstimateLemma:ChevalleyRepresentation:omegaIs3-form}
    \omega \text{ is a non-zero } 3\text{-form on } \hfrak.
\end{equation}
For $X,Y,Z\in\hfrak$, since $\phi$ is a Lie algebra automorphism on $\hfrak$, we have
\begin{align*}
    \det(\phi)\omega(X,Y,Z)&=\omega(\phi(X),\phi(Y),\phi(Z))=\kappa([\phi(X),\phi(Y)],\phi(Z))\\
    &=\kappa(\phi([X,Y]),\phi(Z))=\kappa([X,Y],Z)=\omega(X,Y,Z).
\end{align*}
This and (\ref{eq:MeasureEstimateLemma:ChevalleyRepresentation:omegaIs3-form}) imply that $\det(\phi)=1$, hence from (\ref{eq:MeasureEstimateLemma:ChevalleyRepresentation:rho(g)v_N=det(phi)v_N}), we deduce that $\rho(g)v_\mathbf{L}=v_\mathbf{L}.$ This finishes the proof of \textit{Lemma~\ref{Lemma:MeasureEstimateLemma:ChevalleyRepresentation}}.
\end{proof}

With the above Chevalley representation of $\mathbf{L}$ and $\Mbf=\Core_\Gbf(\Lbf)$ the normal core of $\Lbf$, we can state the following measure estimate lemma:
\begin{Lemma}\label{Lemma:MeasureEstimateLemma:MeasureEstimateLemma}

    Let $\delta>0,\rho=\bigwedge^3\Ad, v_\mathbf{L}=\bu\wedge\ba\wedge\bou$ be as in \textit{Lemma~\ref{Lemma:MeasureEstimateLemma:ChevalleyRepresentation}}, $\norm{\cdot}_{\bigwedge^3\gfrak}$ a norm on $\bigwedge^3\gfrak$. Then there exist $\denew\label{delta:MeasureEstimateLemma:ExponentOfBoundOnRepVectors},\denew\label{delta:MeasureEstimateLemma:BoundOnLatticeCounting}\in(0,1)$ and $\Cnew\label{C:MeasureEstimateLemma:ConstForLowerBoundOfDistOfgammaToNormalCore}\geq1$ such that the following holds. For a compact set $K\subset \mathbf{G}$, let $V_{\epsilon}$ be the set of $g\in K$ with $\norm{\rho(\gamma g^{-1})v_\Lbf-\rho(g^{-1})v_\Lbf}_{\bigwedge^3\gfrak}<\epsilon^{1/\deold{delta:MeasureEstimateLemma:ExponentOfBoundOnRepVectors}}$ for some $\gamma\in\Gamma$ that satisfies $\min_{x\in\Mbf}\norm{x-\gamma}\geq\Cold{C:MeasureEstimateLemma:ConstForLowerBoundOfDistOfgammaToNormalCore}\epsilon^{1/\delta}$ and $\norm{\gamma}\leq\epsilon^{-\deold{delta:MeasureEstimateLemma:ExponentOfBoundOnRepVectors}\deold{delta:MeasureEstimateLemma:BoundOnLatticeCounting}/2}$, then for sufficiently small $\epsilon>0$, $m_{G}(V_{\epsilon})\leq\epsilon^{\deold{delta:MeasureEstimateLemma:ExponentOfBoundOnRepVectors}/4}.$
    
\end{Lemma}
This is \cite{lindenstrauss2025timechangerigidityunipotent}*{Lemma 6.7}, and the proof follows similar steps as in \textit{Section 6} of the same paper. Note that the polynomial upper bound of $\absolute{\{\gamma\in\Gamma:\norm{\gamma}\leq T\}}$ with respect to $T>1$ in \cite{lindenstrauss2025timechangerigidityunipotent}*{(6.17)} can be obtained using a bound by $m_{\Gbf}(\{g\in\Gbf:\norm{g}\leq T+1\})$ (cf. \cite{MR2573139}*{Lemma 6.5}) and a polynomial bound of this set's $m_\Gbf$-measure (cf. \cite{MR2573139}*{Theorem 7.17}, see also \cite{MR1367214}*{p.~141, 142}).

\medskip

In the case $U_t$ is nonsplitting, to check the condition on $\min_{x\in\Mbf}\norm{x-\gamma}$ in \textit{Lemma~\ref{Lemma:MeasureEstimateLemma:MeasureEstimateLemma}}, we recall the following lemma in \cite{lindenstrauss2025timechangerigidityunipotent} and provide an argument that does not require $\Gbf$ to be semisimple:
\begin{Lemma}[\cite{lindenstrauss2025timechangerigidityunipotent}*{Lemma 7.3}]\label{Lemma:MeasureEstimateLemma:NormalcoreAndItsNeighborhood}

    Assume $U_t$ is nonsplitting. Then there exist $\enew\label{epsilon:MeasureEstimateLemma:NormalcoreAndItsNeighborhood}\in(0,1)$ and $\Cnew\label{C:MeasureEstimateLemma:NormalcoreAndItsNeighborhood}>1$ such that if $h\in\Hbf$ satisfies $hzp\in\Mbf$ for some $z\in\Zbf=C_\Gbf(\Hbf)$, $p\in B^\Gbf_\epsilon$, and $\epsilon\in(0,\eold{epsilon:MeasureEstimateLemma:NormalcoreAndItsNeighborhood})$, then the distance of $h$ to the center of $\Hbf$ is within $\Cold{C:MeasureEstimateLemma:NormalcoreAndItsNeighborhood}\epsilon$.
    
\end{Lemma}
\begin{proof}
We claim that 
\[
\Mbf\subseteq\Zbf.
\]
Indeed, let $J=\{1\leq j\leq n:q_j=0\}$, we see that the Lie algebras of $\Zbf$ and $\Lbf$ are $\zfrak=\bigoplus_{j\in J}\R\bx_{0,j}$ and $\lfrak=\hfrak\oplus\zfrak$, respectively. Let $\mfrak$ be the Lie algebra of $\Mbf$, it follows that $\mfrak\subseteq\lfrak$ and $\mfrak\lhd\gfrak$.

Since $\lfrak\neq\gfrak$ by our assumption that $U_t$ is nonsplitting, there exists $j_0$ such that $q_{j_0}\neq0$, and consequently, $\bx_{0,j_0}\notin\lfrak$.

Since $\mathfrak{m}\lhd\gfrak$, $\hfrak\lhd\lfrak=\hfrak\oplus\zfrak$, and $\mathfrak{m}\subseteq\lfrak$, it follows that $[\hfrak,\mathfrak{m}]$ is an ideal of $\hfrak\cong\slfrak_2(\R)$. If $[\hfrak,\mathfrak{m}]=\hfrak$, then
\[
q_{j_0}\bx_{0,j_0}=[\ba,\bx_{0,j_0}]\in\mathfrak{m}\subseteq\mathfrak{l},
\]
a contradiction. Hence, $[\hfrak,\mathfrak{m}]=0$, which implies $\mfrak\subseteq\zfrak$, and consequently, $\Mbf^\circ\subseteq\Zbf^\circ$.

For $m\in\Mbf$, consider $[\Hbf\ni h\mapsto f(h)=[h,m]=hmh^{-1}m^{-1}]$, we see that $f$ is continuous, $f(\ebf_\Gbf)=\ebf_\Gbf$, and $f(h)\in\Mbf\cap\Hbf$ for all $h\in\Hbf$. This implies
\[
f(h)\in\Mbf^\circ\cap\Hbf\subseteq\Zbf^\circ\cap\Hbf,\qquad\forall h\in\Hbf.
\]
And as the center of $\Hbf$ is of order at most $2$, we obtain $[h,m]=f(h)=\ebf_\Gbf$ for all $h\in\Hbf$, i.e. $m\in\Zbf$. Thus, the claim $\Mbf\subseteq\Zbf$ is proved.

Now, choose $\eold{epsilon:MeasureEstimateLemma:NormalcoreAndItsNeighborhood}\in(0,1)$ and let $h\in\Hbf$ be such that $hzp=m\in\Mbf$ for some $z\in\Zbf,p\in B^\Gbf_\epsilon$, $\epsilon\in(0,\eold{epsilon:MeasureEstimateLemma:NormalcoreAndItsNeighborhood})$. This and the facts $\Mbf\lhd\Gbf$, $\Mbf\subseteq\Zbf$, imply
\[
h=mp^{-1}z^{-1}=p^{-1}(pmp^{-1})z^{-1}\quad\text{and}\quad (pmp^{-1})z^{-1}\in\Zbf.
\]
Hence, it follows that
\[
\Ad(h)\vert_\hfrak=\Ad(p^{-1})\circ\Ad((pmp^{-1})z^{-1})\vert_\hfrak=\Ad(p^{-1})\vert_\hfrak.
\]
This, together with $d_{\Gbf}(p^{-1},\ebf_\Gbf)<\epsilon$ and $\Ad:\Hbf\to\Ad(\Hbf)$ being a covering whose kernel is the center of $\Hbf$, implies that by choosing a small enough $\eold{epsilon:MeasureEstimateLemma:NormalcoreAndItsNeighborhood}\in(0,1)$ and a large enough $\Cold{C:MeasureEstimateLemma:NormalcoreAndItsNeighborhood}>1$, depending only on $\Gbf$ and $\Hbf$, the distance of $h$ to the center of $\Hbf$ is within $\Cold{C:MeasureEstimateLemma:NormalcoreAndItsNeighborhood}\epsilon$.
\end{proof}

\subsection{Basic Lemma}\label{subSection:Preliminaries:BasicLemma}
In this subsection, we state the \textit{Basic Lemma}, which establishes a key property of an ergodic nonsplitting $\slfrak_2$-induced unipotent flow $(\Gbf/\Gamma,\mu,U_t)$. We split this into three parts, which in short are:
\begin{itemize}
    \item The first part states that if two ``good" points $x,y\in \mathbf{G}/\Gamma$ are $(\delta,\epsilon,R)$-two sides matchable, with $\delta,\epsilon>0$ sufficiently small and $R>0$ sufficiently large, then $x$ and $U_{\tilde{t}}y$ are related via a ``$(R,\delta)$-small" Kakutani-Bowen ball for some $\tilde{t}$ bounded by $R$.
    \item The second part gives a sharper conclusion compare to the first part, provided $x,y$ are $(\delta,\epsilon,R)$-two sides matchable with the same matching function for all $L_1\leq R\leq L_2<\infty$.
    \item The last part is the second part in the special case $L_2=\infty$.
\end{itemize}

This property is a combination of \textit{Lemma~8.4} and \textit{Corollary~8.5} in \cite{lindenstrauss2025timechangerigidityunipotent}. We fix a constant $\enew\label{epsilon:Preliminaries:CoveringLemma}\in(0,1)$ depending only on $(\Gbf/\Gamma,\mu,U_t)$ (in particular, it can be taken as $C_5^{-1}\epsilon_3$ in \cite{lindenstrauss2025timechangerigidityunipotent}*{Lemma~5.2}), recall the logarithm map $\log:\Ocal\to\gfrak$ in \textit{Definition~\ref{Def:Preliminaries:Kakutani-BowenBalls}}, and state the \textit{Basic Lemma} as follows:
\begin{Lemma}[Basic Lemma]\label{Lemma:BasicLemma:BasicLemma}

    For $\eta\in(0,10^{-3})$, there exist $\denew\label{delta:BasicLemma:BoundOndeltaInTwoSideMatching}>0$, $\Rnew\label{R:BasicLemma:LowerBoundOnRInTwoSideMatching}>1$, $\Cnew\label{C:BasicLemma:BoundOndeltaRatioInKak}>1$, and a compact set $\Knew\label{K:BasicLemma:GoodSetInBasicLemma}\subset \mathbf{G}/\Gamma$ with $\mu(\Kold{K:BasicLemma:GoodSetInBasicLemma})>1-\eta$ satisfying the following property. Let $\delta\in(0,\deold{delta:BasicLemma:BoundOndeltaInTwoSideMatching})$, $\epsilon\in(0,10^{-10}\eold{epsilon:Preliminaries:CoveringLemma})$, $x\in \mathbf{G}/\Gamma,y\in\Kold{K:BasicLemma:GoodSetInBasicLemma}$, and $\xt,\yt\in \mathbf{G}$ be lifts of $x,y$, respectively, with $d_{\mathbf{G}}(\xt,\yt)<\delta.$
    \begin{enumerate}[label=(\Alph*)]
        \item\label{item:BasicLemma:Statement:A} If $x,y$ are $(\delta,\epsilon,R)$-two sides matchable for some $R\geq\Rold{R:BasicLemma:LowerBoundOnRInTwoSideMatching}$ with a $C^1$-matching function $h$, then there exist $s\in[-R,-\frac{1}{2}R]$, $\sol\in[\frac{3}{2}R+s,R]$, $\Tilde{t}\in[-10^{-9}R,10^{-9}R]$, and $\gamma\in\Gamma$ such that
        \begin{equation*}
            \begin{gathered}
                 U_{h(s)}\xt\in\Kak(2R,\Cold{C:BasicLemma:BoundOndeltaRatioInKak}\delta)U_{s}\yt\gamma,\qquad U_{h(\sol)}\xt\in\Kak(2R,\Cold{C:BasicLemma:BoundOndeltaRatioInKak}\delta)U_{\sol}\yt\gamma,\\
                 \xt\in\Kak(2R,\Cold{C:BasicLemma:BoundOndeltaRatioInKak}\delta)U_{\Tilde{t}}\yt\gamma.
            \end{gathered}
        \end{equation*}
        \item\label{item:BasicLemma:Statement:B} If $x,y$ are $(\delta,\epsilon,R)$-two sides matchable for all $R\in[L_1,L_2]$ with the same $C^1$-matching function $h$, where $\Rold{R:BasicLemma:LowerBoundOnRInTwoSideMatching}\leq L_1\leq L_2<\infty$, then there exist $s\in[-L_2,-\textstyle{\frac{1}{2}}L_2]$, $\sol\in[\textstyle{\frac{3}{2}}L_2+s,L_2]$, $\Tilde{t}\in[-10^{-3}L_1,10^{-3}L_1]$, and $\gamma\in\Gamma$ such that
        \begin{equation*}
            \begin{gathered}
            U_{h(s)}\xt\in\Kak(2L_2,\Cold{C:BasicLemma:BoundOndeltaRatioInKak}\delta)U_s\yt\gamma,\qquad U_{h(\sol)}\xt\in\Kak(2L_2,\Cold{C:BasicLemma:BoundOndeltaRatioInKak}\delta)U_{\sol}\yt\gamma,\\
            \xt\in\Kak(1.99L_2,10\Cold{C:BasicLemma:BoundOndeltaRatioInKak}\delta)U_{\Tilde{t}}\yt\gamma.
        \end{gathered}
        \end{equation*}
        \item\label{item:BasicLemma:Statement:C} If $L_2=\infty$ in \ref{item:BasicLemma:Statement:B}, then $x=U_{c_{\bu}}A_{c_{\ba}}gy$ with $c_{\bu}\in[-10^{-2}L_1,10^{-2}L_1]$, $c_{\ba}\in[-10\Cold{C:BasicLemma:BoundOndeltaRatioInKak}\delta,10\Cold{C:BasicLemma:BoundOndeltaRatioInKak}\delta]$, and $g\in \overline{B}^{\mathbf{G}}_{10\Cold{C:BasicLemma:BoundOndeltaRatioInKak}\delta}$, $\log(g)\in\bigoplus_{j=1}^n\R\bx_{0,j}$. 
        Furthermore, $c_{\bu},c_{\ba},g$ are unique in the following sense: if 
        \[
        c'_{\bu}\in\R,\qquad c'_{\ba}\in[-10\Cold{C:BasicLemma:BoundOndeltaRatioInKak}\delta,10\Cold{C:BasicLemma:BoundOndeltaRatioInKak}\delta],\qquad g'\in \overline{B}^{\mathbf{G}}_{10\Cold{C:BasicLemma:BoundOndeltaRatioInKak}\delta}
        \]
        satisfying 
        \[
        \textstyle{\log(g')\in\bigoplus_{j=1}^n\R\bx_{0,j},\qquad x=U_{c'_{\bu}}A_{c'_{\ba}}g'y,}
        \]
        then 
        \[
        c'_{\bu}=c_{\bu},\qquad c'_{\ba}=c_{\ba},\qquad g'=g.
        \]
    \end{enumerate}
    
\end{Lemma}

As mentioned previously, this is exactly \cite{lindenstrauss2025timechangerigidityunipotent}*{Lemma~8.4 and~Corollary 8.5}. In particular, part \ref{item:BasicLemma:Statement:A} corresponds to \textit{Lemma 8.4}, while parts \ref{item:BasicLemma:Statement:B} and \ref{item:BasicLemma:Statement:C} correspond to \textit{Corollary 8.5}. For the detailed proof of \textit{Lemma~\ref{Lemma:BasicLemma:BasicLemma}}, we refer to \cite{lindenstrauss2025timechangerigidityunipotent}*{Subsections~8.6,~8.7} and the uniqueness argument in the proof of \cite{lindenstrauss2025timechangerigidityunipotent}*{Claim~10.6}.

\subsection{\texorpdfstring{$\epsilon$}{}-good time changes}\label{subSection:Preliminaries:epsilon-goodTimeChange}
We end this section by recalling the regularity condition $\epsilon$-good on time changes in \cite{lindenstrauss2025timechangerigidityunipotent}*{Section~4}:
\begin{Definition}[\cite{lindenstrauss2025timechangerigidityunipotent}*{Definition 4.2}]\label{Definition:epsilon-goodTimeChange}
    
    A time change $\tau\in L^1_+(X,\mu)$ for a flow $T_t$ on a probability space $(X,\mu)$ is said to be \textit{$\epsilon$-good} for some $\epsilon\in(0,10^{-3})$ if $\int_X\tau d\mu=1$, $\operatorname{esssup}\absolute{\tau-1}<\epsilon$, and there exists $Z\subseteq X$ of full measure such that for $x\in Z,\tau(T_tx)$ is a $C^{\infty}$ in $t$ and $\absolute{\tau(T_tx)-1}<\epsilon$ for all $t\in\R$.
    
\end{Definition}
Following \cite{lindenstrauss2025timechangerigidityunipotent}*{Lemmas 4.1, 4.3, Appendix B}, we deduce a reduction from considering arbitrary Kakutani factor maps to considering those that have $\epsilon$-good corresponding time changes:
\begin{Proposition}\label{Prop:epsilon-goodTimeChange:epsilon-goodTimeChange}
    Let $\psi:\Gbf_1/\Gamma_1\to\Gbf_2/\Gamma_2$ be a Kakutani factor map between ergodic $\slfrak_2$-induced unipotent flows $(\Gbf_i/\Gamma_i,\mu_i,U^{(i)}_t)$, $i=1,2$. Denote $A^{(1)}_t=\exp(t\ba^{(1)})$, where $\{\bu^{(1)},\ba^{(1)},\bou^{(1)}\}$ is an $\slfrak_2$-triple corresponding to $U^{(1)}_t$. Then for $\epsilon\in(0,10^{-3})$, there exist $p_0\in\R$, $\sigma:\Gbf_1/\Gamma_1\to\R$ measurable, and a Kakutani factor map $\Tilde{\psi}:\Gbf_1/\Gamma_1\to\Gbf_2/\Gamma_2$ between $(\Gbf_i/\Gamma_i,\mu_i,U^{(i)}_t)$, $i=1,2$, such that $\psi(\cdot)=\Tilde{\psi}(U^{(1)}_{\sigma(\cdot)}A^{(1)}_{p_0}\cdot)$ and the time change $\tau\in L^1_+(\Gbf_1/\Gamma_1,\mu_1)$ corresponding to $\Tilde{\psi}$ is $\epsilon$-good. Moreover, if $\psi$ is a Kakutani equivalence, then so is $\Tilde{\psi}$.
    
\end{Proposition}
The reason why we want the time change corresponding to a Kakutani factor map to be $\epsilon$-good is that then, the Kakutani factor map roughly preserves two sides matchability of two points. More explicitly:
\begin{Lemma}[\cite{lindenstrauss2025timechangerigidityunipotent}*{Lemma 4.4}]\label{Lemma:epsilon-goodTimeChange:PreserveTwoSideMatching}

    Let $T^{(i)}_t$ be an ergodic flow on a metric space $(X_i,d_i)$ with respect to a Borel probability measure $\mu_i$, $\psi:X_1\to X_2$ a Kakutani factor map between $(X_i,\mu_i,T^{(i)}_t),i=1,2$, where the corresponding time changes $\tau\in L^1_+(X_1,\mu_1)$ is $\epsilon$-good for some $\epsilon\in(0,10^{-3})$. Then for $\eta\in(0,\epsilon)$, there exist a compact set $\Knew\label{K:epsilon-goodTimeChange:CompactSet}\subseteq X_1$ with $\mu_1(\Kold{K:epsilon-goodTimeChange:CompactSet})>1-\eta$ and $\Rnew\label{R:epsilon-goodTimeChange:SmallestMatchingLength}>1$ such that the following holds. Given $\delta\in(0,\epsilon)$, there exists $\enew\label{epsilon:epsilon-goodTimeChange:psiUniformContinous}(\delta)\in(0,\delta)$ such that for $x,y\in\Kold{K:epsilon-goodTimeChange:CompactSet}$,
    \begin{itemize}
        \item if $d_1(x,y)<\eold{epsilon:epsilon-goodTimeChange:psiUniformContinous}(\delta)$ then $d_2(\psi(x),\psi(y))<\delta$,
        \item and if for some $R\geq\Rold{R:epsilon-goodTimeChange:SmallestMatchingLength},x$ and $y$ are $(\eold{epsilon:epsilon-goodTimeChange:psiUniformContinous}(\delta),\epsilon,R)$-two sides matchable with a $C^1$-matching function $h:[-R,R]\to[-R,R]$, then $\psi(x)$ and $\psi(y)$ are $(\delta,30\epsilon,R)$-two sides matchable with the $C^1$-matching function $\Tilde{h}$ given by
        \[
        \Tilde{h}(t)=\omega(x,h(\rho(y,t))),
        \]
        where $\omega$ and $\rho$ are determined by
        \[
        \omega(u,s)=\int_0^s\tau(T^{(1)}_ru)dr\quad\text{and}\quad s=\int_0^{\rho(u,s)}\tau(T^{(1)}_ru)dr.
        \]
    \end{itemize}
\end{Lemma}

\section{Structural properties of ergodic \texorpdfstring{$\slfrak_2$}{}-induced unipotent flows}\label{Section:StructureOfErgodicUnipotentFlow}
In this section, we provide some structural properties of the base linear algebraic group and the homogeneous space given an ergodic $\slfrak_2$-induced unipotent flow. The results and arguments in this section are based on the theory of linear algebraic groups and their lattices. For background, we refer to \cites{MR1102012,MR507234,MR776417,onishchik2000discrete}.

\medskip

The following setting are fixed for this section: 
\begin{itemize}
    \item $\Gbf$ is the identity component of a real linear algebraic group;
    \item Lie algebra $\gfrak=\Lie(\Gbf)$, and $\Gamma<\Gbf$ is a lattice in $\Gbf$;
    \item $\mu$ is the Haar probability measure on $\Gbf/\Gamma$;
    \item $U_t$ is an ergodic $\slfrak_2$-induced unipotent flow on $(\Gbf/\Gamma,\mu)$;
    \item $\hfrak$ is the $\slfrak_2(\R)$-subalgebra corresponding to $U_t$.
\end{itemize}

With these settings, the structural results in \cite{MR910358}*{Section 4C} can be stated as following:
\begin{Lemma}[\cite{MR910358}*{Proposition 4.19, 4.22, 4.23}]\label{Lemma:StructuralProperties:WitteResults}

    In this section's setting, 
    \begin{enumerate}[label=(\alph*)]
        \item\label{item:WitteResults:radG=nilG} the nilradical $\Rbf$ of $\Gbf$, i.e. the maximal connected normal nilpotent subgroup of $\Gbf$, coincides with the radical of $\Gbf$, i.e. the maximal connected normal solvable subgroup of $\Gbf$;
        \item\label{item:WitteResults:GammacapNilG/Z(G)islattice} $\Gamma\cap\Rbf$ and $\Gamma\cap Z(\Gbf)$ are lattices in $\Rbf$ and $Z(\Gbf)$, respectively;
        \item\label{item:WitteResults:GammaProjectDenselyOnCompactFactor}$\Gamma$ projects densely into the maximal compact semisimple factor of $\Gbf$;
        \item\label{item:WitteResults:[G,G]Z(G)=G} $[\Gbf,\Gbf]Z(\Gbf)=\Gbf$.
    \end{enumerate}
    
\end{Lemma}
Let $\Nbf$ be the unipotent radical of $\Gbf$, it follows from \textit{Lemma~\ref{Lemma:StructuralProperties:WitteResults}-\ref{item:WitteResults:radG=nilG}}, the \textit{Jordan decomposition} (cf. \cite{MR1102012}*{Theorem 4.4}), and the statement and proof of \cite{MR1102012}*{Theorem 10.6}, we have the following decomposition
\begin{equation}\label{eq:StructuralProperties:R=T_cxT_(nc)xN}
    \Rbf=\Tbf\Nbf=\Tbf_{\operatorname{c}}\Tbf_{\nc}\Nbf\cong\Tbf\times\Nbf\cong\Tbf_{\operatorname{c}}\times\Tbf_{\nc}\times\Nbf,
\end{equation}
where $\Tbf$ is the connected, central subgroup of $\Rbf$ consisting of semisimple elements, and 
\[
\Tbf_{\operatorname{c}}\cong\T^a,\qquad \Tbf_{\nc}\cong(\R_{>0})^b\label{page:T_c=a-torus}
\]
are the compact, noncompact factor of $\Tbf$, respectively ($a,b\geq0$). 

Apply \cite{MR1102012}*{Theorem 10.6} again, we obtain that $\Tbf$ is the unique maximal torus of $\Rbf$. Then by \cite{MR92928} and \cite{MR1102012}*{Proposition 11.23}, $\Tbf$ is the identity component of the center of a Levi subgroup of $\Gbf$. Since $\Tbf$ is central in the nilradical $\Rbf$ of $\Gbf$, we deduce that
\begin{equation}\label{eq:StructuralProperties:TisCentralInG}
    \Tbf\leq Z(\Gbf)^\circ,    
\end{equation}
where $Z(\Gbf)^\circ$ is the identity component of $Z(\Gbf)$.

Consider a Levi subgroup $\Sbf\Tbf$ that contains $\Hbf=\langle\exp(\hfrak)\rangle$, where $\Sbf$ is its semisimple factor. It follows that $\Sbf$ is the identity component of a semisimple real linear algebraic group. Moreover, $\Sbf$ is a maximal semisimple factor of $\Gbf$ so that $\Hbf\subseteq\Sbf$ and satisfying the following:
\begin{equation}\label{eq:StructuralProperties:G=(ST)ltimesN,ST=(ST_c)xT_(nc)}
    \begin{gathered}
        \Gbf=\Sbf\Rbf=(\Sbf\Tbf)\Nbf\cong(\Sbf\Tbf)\ltimes\Nbf,\qquad\Sbf\Tbf=\Sbf\Tbf_{\operatorname{c}}\Tbf_{\nc}\cong(\Sbf\Tbf_{\operatorname{c}})\times\Tbf_{\nc},\\
        \Sbf\cap\Rbf=\Sbf\cap\Tbf_{\operatorname{c}},\qquad \Sbf\cap\Rbf<Z(\Gbf)^\circ,\qquad 
        \absolute{\Sbf\cap\Rbf}<\infty.
    \end{gathered}
\end{equation}

Since $\Gamma\cap\Rbf$ is a lattice in $\Rbf$, where $\Rbf$ is the nilradical (and radical) of $\Gbf$ (cf. \textit{Lemma~\ref{Lemma:StructuralProperties:WitteResults}-\ref{item:P=[S,R]:1},\ref{item:P=[S,R]:2}}), by \cite{onishchik2000discrete}*{I.1 Theorem 4.3, 5.7}, the image of $\Gamma$ in $\Gbf/\Rbf\cong\Sbf/(\Sbf\cap\Rbf)$ under the natural map is a lattice in $\Sbf/(\Sbf\cap\Rbf)$, hence its preimage $\Gamma_\Sbf$ with respect to the natural map $\Sbf\to\Sbf/(\Sbf\cap\Rbf)$ is a lattice in $\Sbf$. We call $\Gamma_\Sbf$ the lattice in $\Sbf$ induced from $\Gamma$. And in the paper, we use this terminology to refer to the lattice construction above.

From now on, we also fix the notations $\Rbf, \Nbf, \Tbf,\Tbf_{\operatorname{c}},\Tbf_{\nc}, \Sbf,\Gamma_\Sbf$, along with the arguments above, in the beginning setting of this section.

\medskip

Now, let us consider the subgroup
\[
\Pbf=[\Sbf,\Rbf],
\]
i.e. the closed subgroup generated by $\{srs^{-1}r^{-1}:s\in\Sbf,r\in\Rbf\}$, of $\Rbf$. This subgroup has some useful properties, which we provide in the following:
\begin{Lemma}\label{Lemma:StructuralProperties:P=[S,R]}
    \begin{enumerate}[label=(\alph*)]
        \item\label{item:P=[S,R]:1} $\Pbf$ is a closed, simply connected, connected, nilpotent, normal subgroup of $\Nbf$ and of $\Gbf$ such that $[\Sbf,\Nbf]=\Pbf$;
        \item\label{item:P=[S,R]:2} $[\Sbf,\Pbf]=\Pbf$;
        \item\label{item:P=[S,R]:3}$[\Gbf,\Rbf]=\Pbf$ and $[\Gbf,\Gbf]=\Sbf\Pbf\cong\Sbf\ltimes\Pbf$;
        \item\label{item:P=[S,R]:4} $\Gamma\cap\Pbf$ is a lattice in $\Pbf$.
    \end{enumerate}
    
\end{Lemma}
\begin{proof}
\ref{item:P=[S,R]:1} By (\ref{eq:StructuralProperties:R=T_cxT_(nc)xN}) and (\ref{eq:StructuralProperties:TisCentralInG}), we deduce that
\[
\Pbf=[\Sbf,\Rbf]=\overline{\langle[\Sbf,\Tbf],[\Sbf,\Nbf]\rangle}=[\Sbf,\Nbf]
\]
is a closed, connected subgroup of $\Nbf$. With this and the fact that the unipotent radical $\Nbf$ of $\Gbf$ is simply connected, connected, nilpotent, we deduce that $\Pbf$ is simply connected and nilpotent. And the condition $\Pbf\lhd\Gbf$ can be directly checked.

\medskip

\noindent \ref{item:P=[S,R]:2} Let $\sfrak$, $\nfrak$, and $\pfrak$ be the Lie algebra of $\Sbf$, $\Nbf$, and $\Pbf$, respectively. Consider $\nfrak$ as an $\sfrak$-representation via the adjoint map. Since $\sfrak$ is semisimple, $\nfrak$ splits into irreducible $\sfrak$-representations. Hence, the subspace $[\sfrak,\nfrak]$ spanned by $\{[X,Y]:X\in\sfrak,Y\in\nfrak\}$ is the direct sum of nontrivial irreducible $\sfrak$-representations in $\nfrak$. This implies that
\begin{equation}\label{eq:StructuralProperties:P=[S,R]:[sfrak,[sfrak,nfrak]]}
    [\sfrak,[\sfrak,\nfrak]]=[\sfrak,\nfrak].
\end{equation}
Apply part \ref{item:P=[S,R]:1} to the subgroup $\Sbf\Nbf\cong\Sbf\ltimes\Nbf$, it follows that $\pfrak$ is the Lie algebra generated by $[\sfrak,\nfrak]$. Similarly for $\Pbf$ and $\Sbf\Pbf\cong\Sbf\ltimes\Pbf$, the Lie algebra $\Lie([\Sbf,\Pbf])$ of $[\Sbf,\Pbf]$ is generated by $[\sfrak,\pfrak]$. Thus, combining these with (\ref{eq:StructuralProperties:P=[S,R]:[sfrak,[sfrak,nfrak]]}) gives us $\Lie([\Sbf,\Pbf])=\pfrak$, which implies that $[\Sbf,\Pbf]=\Pbf$.

\medskip

\noindent\ref{item:P=[S,R]:3} We show $[\Rbf,\Rbf]\subseteq\Pbf$. With this, (\ref{eq:StructuralProperties:G=(ST)ltimesN,ST=(ST_c)xT_(nc)}) implies that
\begin{align*}
    &[\Gbf,\Rbf]=\overline{\langle[\Sbf,\Rbf],[\Rbf,\Rbf]\rangle}=\Pbf\\
    &[\Gbf,\Gbf]=\overline{\langle[\Sbf,\Sbf],[\Sbf,\Rbf],[\Rbf,\Rbf]\rangle}=\overline{\langle\Sbf,\Pbf\rangle}=\Sbf\Pbf\cong\Sbf\ltimes\Pbf.
\end{align*}
Since $\Pbf=[\Sbf,\Rbf]\lhd\Gbf$ by part \ref{item:P=[S,R]:1}, we can directly check that $\Sbf\Pbf\lhd\Gbf$. Projecting the relation $[\Gbf,\Gbf]Z(\Gbf)=\Gbf$ (cf. \textit{Lemma~\ref{Lemma:StructuralProperties:WitteResults}-\ref{item:WitteResults:[G,G]Z(G)=G}}) via the natural map $\Gbf\to\Gbf/(\Sbf\Pbf)=\Qbf'$, it follows that $[\Qbf',\Qbf']Z(\Qbf')=\Qbf'$. This implies $[\qfrak',\qfrak']+Z(\qfrak')=\qfrak'$, where $\qfrak'=\Lie(\Qbf')$. And as (\ref{eq:StructuralProperties:G=(ST)ltimesN,ST=(ST_c)xT_(nc)}) implies
\[
\Qbf'=\Gbf/(\Sbf\Pbf)\cong\Rbf/((\Sbf\cap\Tbf_{\operatorname{c}})\Pbf),
\]
we deduce $\Qbf'$ is nilpotent, and so by induction, $\qfrak'=Z(\qfrak')$. Hence, we get $\Qbf'\cong\Rbf/((\Sbf\cap\Tbf_{\operatorname{c}})\Pbf)$ is abelian, which shows $[\Rbf,\Rbf]\subseteq(\Sbf\cap\Tbf_{\operatorname{c}})\Pbf$. This, (\ref{eq:StructuralProperties:R=T_cxT_(nc)xN}), (\ref{eq:StructuralProperties:G=(ST)ltimesN,ST=(ST_c)xT_(nc)}), the fact that $\Pbf\subseteq\Nbf$ (cf. part \ref{item:P=[S,R]:1}), and the fact that $[\Rbf,\Rbf]$ is a connected subgroup of $\Rbf$ imply that $[\Rbf,\Rbf]\subseteq\Pbf$.

\medskip

\noindent\ref{item:P=[S,R]:4} Denote $\Gamma_\Rbf=\Gamma\cap\Rbf$ and $\Lambda=[\Gamma,\Gamma_\Rbf]$, it follows that $\Gamma_\Rbf$ is a lattice in $\Rbf$ (cf. \textit{Lemma~\ref{Lemma:StructuralProperties:WitteResults}-\ref{item:WitteResults:GammacapNilG/Z(G)islattice}}) and $\Lambda\leq\Gamma\cap\Pbf$ are discrete subgroups of $\Pbf$. Since $\Pbf$ is a connected subgroup of the unipotent radical $\Nbf$ of $\Gbf$, $\Pbf$ is Zariski closed. Hence, the Zariski closure $\overline{\Lambda}^{\operatorname{Zar}}$ of $\Lambda$ is contained in $\Pbf$.

Since $\Gamma_\Rbf$ is a lattice in (cf. \eqref{eq:StructuralProperties:R=T_cxT_(nc)xN})
\[
\Rbf=\Tbf\Nbf=\Tbf_{\operatorname{c}}\Tbf_{\nc}\Nbf\cong\T^a\times(\R_{>0})^b\times\Nbf,
\]
it follows from \cite{MR910358}*{Corollary 4.2} that $\overline{\Gamma_\Rbf}^{\operatorname{Zar}}$ contains $\Nbf$. And by \textit{Borel density theorem} (cf. \cite{MR776417}*{Theorem 3.2.5}) and \textit{Lemma~\ref{Lemma:StructuralProperties:WitteResults}-\ref{item:WitteResults:GammaProjectDenselyOnCompactFactor}}, we obtain that $\overline{\Gamma_\Sbf}^{\operatorname{Zar}}$ contains $\Sbf$.

Since $[\Gamma_\Rbf,\Gamma_\Rbf]\subseteq\Lambda$ and $\Nbf\subseteq\overline{\Gamma_\Rbf}^{\operatorname{Zar}}$, we have
\begin{equation}\label{eq:StructuralProperties:[N,N]subsetZarClosusreOfLambdasubsetP}
    [\Nbf,\Nbf]\subseteq\overline{\Lambda}^{\operatorname{Zar}}\subseteq \Pbf.
\end{equation}
Let $V=\Nbf/[\Nbf,\Nbf]$ and $q:\Rbf\to V$ the homomorphism defined by
\[
q(tn)=n[\Nbf,\Nbf],\qquad t\in\Tbf,\ n\in\Nbf.
\]
For $\gamma=sr\in\Gamma$, where $s\in\Gamma_\Sbf$ and $r\in\Rbf$, we have
\[
q([\gamma,\eta])=s\cdot q(\eta)-q(\eta)\in q(\overline{\Lambda}^{\operatorname{Zar}}),\qquad \eta\in\Gamma_\Rbf,
\]
where $s\cdot q(\eta)=q(s\eta s^{-1})$. As $q\left(\overline{\Lambda}^{\operatorname{Zar}}\right)$ is a vector subspace, $q\left(\overline{\Gamma_\Rbf}^{\operatorname{Zar}}\right)$ contains $q(\Nbf)=V$, and $\Sbf\subseteq\overline{\Gamma_\Sbf}^{\operatorname{Zar}}$, it follows that
\[
(s-\id)V\subseteq q(\overline{\Lambda}^{\operatorname{Zar}}),\qquad\forall s\in\Sbf.
\]
Hence, $q(\Pbf)=q([\Sbf,\Nbf])\subseteq q\left(\overline{\Lambda}^{\operatorname{Zar}}\right)$. This and \eqref{eq:StructuralProperties:[N,N]subsetZarClosusreOfLambdasubsetP} give us $\overline{\Lambda}^{\operatorname{Zar}}=\Pbf$. Then, by \cite{MR507234}*{Theorem~2.3}, $\Lambda$ is a cocompact lattice in $\Pbf$. And as $\Lambda\leq\Gamma\cap\Pbf$ and $\Gamma\cap\Pbf$ is discrete, $\Gamma\cap\Pbf$ is also a cocompact lattice in $\Pbf$.
\end{proof}

It follows from \textit{Lemma~\ref{Lemma:StructuralProperties:P=[S,R]}} that the subgroup $\Sbf\Pbf\cong\Sbf\ltimes\Pbf$ is perfect. In fact, $\Gbf$ is essentially the direct product of $\Sbf\Pbf$ and a connected abelian subgroup. We establish this in the following lemma by constructing this subgroup (or rather, its noncompact factor) to be ``compatible" with $\Gamma$:
\begin{Lemma}\label{Lemma:StructuralProperties:Q_(nc)=R^dCompatibleWithGamma}
    In this section's setting, let $\Pbf=[\Sbf,\Rbf]$. Then there exist closed subgroups $\Qbf_{\operatorname{c}}$ and $\Qbf_{\nc}$ of $Z(\Gbf)$ isomorphic to $\T^a$ ($a\geq0$ given on p.~\pageref{page:T_c=a-torus}) and $\R^d$ ($d\geq0$), such that the following decompositions hold:
    \[
    \Rbf=\Pbf\Qbf_{\operatorname{c}}\Qbf_{\nc}\cong\Pbf\times\Qbf_{\operatorname{c}}\times\Qbf_{\nc}\qquad \quad\Gbf=\Sbf\Pbf\Qbf_{\operatorname{c}}\Qbf_{\nc}\cong((\Sbf\ltimes\Pbf)\Qbf_{\operatorname{c}})\times\Qbf_{\nc}
    \]
    Moreover, $\Sbf\cap\Qbf_{\operatorname{c}}$ is a finite subgroup, and $\Gamma\cap\Qbf_{\nc}$ is a lattice in $\Qbf_{\nc}$.
\end{Lemma}
\begin{proof}
Consider $\Tilde{\Gbf}=(\Sbf\times\R^a\times(\R_{>0})^b)\ltimes\Nbf$ equipped with the group law
\[
(s_1,x_1,y_1,n_1)(s_2,x_2,y_2,n_2)=(s_1s_2,x_1+x_2,y_1y_2,n_1(s_1n_2s_1^{-1})),
\]
it follows that the map $\pi:\Tilde{\Gbf}\to\Gbf$ given by
\[
\pi(s,x,y,n)=nst_{\operatorname{c}}(x)t_{\nc}(y),\qquad \forall(s,x,y,n)\in\Tilde{\Gbf},
\]
where $t_{\nc}$ is the isomorphism $(\R_{>0})^b\to\Tbf_{\nc}$ and $t_{\operatorname{c}}:\R^a\to\Tbf_{\operatorname{c}}$ is the composition of the isomorphism $\T^a\xrightarrow[]{}\Tbf_{\operatorname{c}}$ and $[\R^a\to\T^a:x\mapsto x\Z^a]$, is a continuous surjective homomorphism, with
\begin{equation*}
    \ker\pi=\{(s,x,1,e_\Nbf):s\in\Sbf\cap\Tbf_{\operatorname{c}},x\in x_s+\Z^a\text{ for some }x_s\in t_{\operatorname{c}}^{-1}(s^{-1})\}.
\end{equation*}
Denote $\Tilde{\Gamma}=\pi^{-1}(\Gamma)$, it follows that $\Tilde{\Gamma}$ is a discrete subgroup of $\Tilde{\Gbf}$. Then, we can directly check that $\Tilde{\Gbf}/\Tilde{\Gamma}\to\Gbf/\Gamma:\Tilde{g}\Tilde{\Gamma}\mapsto\pi(\Tilde{g})\Gamma$ is a well-defined continuous bijective closed equivariant map, hence a homeomorphism between homogeneous spaces. This implies that $\Tilde{\Gamma}$ is a lattice in $\Tilde{\Gbf}$ and $U_t$ is an ergodic $\slfrak_2$-induced unipotent flow on $(\Tilde{\Gbf}/\Tilde{\Gamma},\Tilde{\mu})$, where $\Tilde{\mu}$ is the Haar probability measure on $\Tilde{\Gbf}/\Tilde{\Gamma}$. Thus, we can apply all the previous arguments in this section to $\Tilde{\Gbf}$ and obtain from \textit{Lemma~\ref{Lemma:StructuralProperties:WitteResults}}, \textit{Lemma~\ref{Lemma:StructuralProperties:P=[S,R]}}, and $\Tilde{\Gbf}$'s construction that
\begin{equation}\label{eq:StructuralProperties:Q_(nc)=R^d:Tilde(R),Z(Tilde(G))}
    \begin{aligned}
        &\Tilde{\Gamma}\cap\Tilde{\Rbf}\text{ is a lattice in }\Tilde{\Rbf}=\R^a\times(\R_{>0})^b\times\Nbf;\\
        &\Tilde{\Gamma}\cap Z(\Tilde{\Gbf})\text{ is a lattice in }Z(\Tilde{\Gbf});\\
        &\Pbf=[\Sbf,\Tilde{\Rbf}]\text{ satisfies }\Tilde{\Gamma}\cap\Pbf\text{ is a lattice in }\Pbf\\
        &\qquad\qquad\qquad\quad\quad\text{ and }\Tilde{\Gbf}=((\Sbf\times\{0\}\times\{1\})\ltimes\Pbf)Z(\Tilde{\Gbf}).
    \end{aligned}\tag{$\blacklozenge$}
\end{equation}

From the above construction and (\ref{eq:StructuralProperties:G=(ST)ltimesN,ST=(ST_c)xT_(nc)}), we see that once we construct a closed subgroup $\Tilde{\Qbf}$ of $Z(\Tilde{\Gbf})$ isomorphic to $\R^d$ ($d\geq 0$) such that $\Tilde{\Gamma}\cap\Tilde{\Qbf}$ is a lattice in $\Tilde{\Qbf}$ and the decompositions
\[
\Tilde{\Rbf}=(\R^a\times\Pbf)\Tilde{\Qbf}\cong\Pbf\times\R^a\times\Tilde{\Qbf}\quad\qquad\Tilde{\Gbf}=((\Sbf\times\R^a)\ltimes\Pbf)\Tilde{\Qbf}\cong(\Sbf\ltimes\Pbf)\times\R^a\times\Tilde{\Qbf}
\]
hold, the subgroups $\Qbf_{\operatorname{c}}=\Tbf_{\operatorname{c}}$ and $\Qbf_{\nc}=\pi(\Tilde{\Qbf})$ satisfy the requirements of \textit{Lemma~\ref{Lemma:StructuralProperties:Q_(nc)=R^dCompatibleWithGamma}}.

We construct the desired subgroup $\Tilde{\Qbf}$ of $Z(\Tilde{\Gbf})$ as follows. By construction, $\Tilde{\Rbf}$ is simply connected, connected, and nilpotent, so it is diffeomorphic to $\Lie(\Tilde{\Rbf})=\Lie(\Rbf)$ via the  diffeomorphisms $\exp$ and $\log$. By \eqref{eq:StructuralProperties:Q_(nc)=R^d:Tilde(R),Z(Tilde(G))}, it follows that $\pfrak=\Lie(\Pbf)$ and $Z(\Lie(\Tilde{\Gbf}))=\Lie(Z(\Tilde{\Gbf}))$ are $\log(\Tilde{\Gamma}\cap\Tilde{\Rbf})$-rational subspaces of $\Lie(\Tilde{\Rbf})$. And as the construction implies $\Z^a\subseteq\Tilde{\Gamma}\cap\R^a$, we also have $\Lie(\R^a)$ is a $\log(\Tilde{\Gamma}\cap\Tilde{\Rbf})$-rational subspace of $\Lie(\Tilde{\Rbf})$. Thus, there exists a $\log(\Tilde{\Gamma}\cap\Tilde{\Rbf})$-rational subspace $\qfrak$ of $\Lie(\Tilde{\Rbf})$ such that
\[
Z(\Lie(\Tilde{\Gbf}))=(\pfrak\cap Z(\Lie(\Tilde{\Gbf})))\oplus\Lie(\R^a)\oplus\qfrak.
\]
This and \eqref{eq:StructuralProperties:Q_(nc)=R^d:Tilde(R),Z(Tilde(G))} imply that with $\sfrak=\Lie(\Sbf)$,
\[
\Lie(\Tilde{\Rbf})=\Lie(\R^a)\oplus\pfrak\oplus\qfrak,\qquad\Lie(\Tilde{\Gbf})=(\sfrak\ltimes\pfrak)\oplus\Lie(\R^a)\oplus\qfrak.
\]
Now, letting $\Tilde{\Qbf}=\exp(\qfrak)$, we obtain our desired subgroup and finish the proof of \textit{Lemma~\ref{Lemma:StructuralProperties:Q_(nc)=R^dCompatibleWithGamma}}.
\end{proof}

It follows from \textit{Lemma~\ref{Lemma:StructuralProperties:Q_(nc)=R^dCompatibleWithGamma}} that $\Gbf/\Gamma$ is an extension of a homogeneous space over $\Sbf\Pbf$. In particular, \textit{Lemma~\ref{Lemma:StructuralProperties:Q_(nc)=R^dCompatibleWithGamma}} implies that 
\[
\Qbf=\Qbf_{\operatorname{c}}\Qbf_{\nc}\cong\Qbf_{\operatorname{c}}\times\Qbf_{\nc}\cong\T^a\times\R^d
\]
is a closed subgroup of $Z(\Gbf)$ such that
\begin{equation}\label{eq:StructuralProperties:G=SPQ}
    \begin{aligned}
        &\text{the decomposition }\Gbf=\Sbf\Pbf\Qbf\text{ holds, } \Sbf\Pbf\cap\Qbf=\Sbf\cap\Qbf_{\operatorname{c}}\text{ is a finite}\\
        &\text{subgroup, and } \Gamma_{\Qbf}=\Gamma\cap\Qbf\text{ is a lattice in }\Qbf.
    \end{aligned}\tag{$\bigstar$}
\end{equation}
This implies that the image of $\Gamma$ in $\Sbf\Pbf/(\Sbf\cap\Qbf_{\operatorname{c}})\cong\Gbf/\Qbf$ under the natural map(s) is a lattice in $\Sbf\Pbf/(\Sbf\cap\Qbf_{\operatorname{c}})$ (cf. \cite{onishchik2000discrete}*{I.1 Theorem 4.3, 5.7}). Let $\Gamma_{\Sbf\Pbf}$ be its preimage in $\Sbf\Pbf$ with respect to the natural map $\Sbf\Pbf\to\Sbf\Pbf/(\Sbf\cap\Qbf_{\operatorname{c}})$, it follows that $\Gamma_{\Sbf\Pbf}$ is a lattice in $\Sbf\Pbf$ containing $\Sbf\cap\Qbf_{\operatorname{c}}$. Then, choose a fundamental domain $F_{\Sbf\Pbf}\subset\Sbf\Pbf$ for $\Gamma_{\Sbf\Pbf}$, and for $(s,x)\in\Sbf\Pbf\times F_{\Sbf\Pbf}$, denote
\begin{itemize}
    \item $s\cdot x\in F_{\Sbf\Pbf}$ be such that $(s\cdot x)^{-1}sx\in\Gamma_{\Sbf\Pbf}$, i.e. induced from the $\Sbf\Pbf$-action on $\Sbf\Pbf/\Gamma_{\Sbf\Pbf}$;
    \item $\beta(s,x)=q\Gamma_{\Qbf}$, where $q\in\Qbf$ satisfies $(s\cdot x)^{-1}sxq\in\Gamma$;
\end{itemize}
we can identify $\Gbf/\Gamma$ with $F_{\Sbf\Pbf}\times(\Qbf/\Gamma_{\Qbf})$ using the $\Gbf$-equivariant map
\begin{equation}\label{eq:StructuralProperties:(x,yGamma_Q)mapstoxyGamma}
    F_{\Sbf\Pbf}\times(\Qbf/\Gamma_{\Qbf})\to\Gbf/\Gamma:(x,y\Gamma_{\Qbf})\mapsto xy\Gamma,
\end{equation}
where the $\Gbf$-action on $F_{\Sbf\Pbf}\times(\Qbf/\Gamma_{\Qbf})$ is given for $(x,y\Gamma_{\Qbf})\in F_{\Sbf\Pbf}\times(\Qbf/\Gamma_{\Qbf})$ and $g\in\Gbf$ by
\begin{equation}\label{eq:StructuralProperties:G-action}
    g\cdot(x,y\Gamma_{\Qbf})=(s\cdot x,(qy\Gamma_{\Qbf})\beta(s,x)^{-1}),\quad s\in\Sbf\Pbf,q\in\Qbf\text{ such that }g=sq.
\end{equation}
Moreover, rescale the Haar measure $m_{\Sbf\Pbf}$ on $\Sbf\Pbf$ so that $m_{\Sbf\Pbf}(F_{\Sbf\Pbf})=1$ and let $\mu_\Qbf$ be the Haar probability measure on $\Qbf/\Gamma_{\Qbf}$, the above $\Gbf$-equivariant map is also a measurable isomorphism between $(F_{\Sbf\Pbf}\times(\Qbf/\Gamma_{\Qbf}),m_{\Sbf\Pbf}\times\mu_\Qbf)$ and $(\Gbf/\Gamma,\mu)$.
\begin{Remark}\label{Remark:StructuralProperties:Settings}

    For the remainder of this paper, for convenience, whenever we refer to the setting of this section, we mean all the following:
    \begin{itemize}
        \item the setting mentioned in the beginning of this section;
        \item the notations and arguments following \textit{Lemma~\ref{Lemma:StructuralProperties:WitteResults}};
        \item the notations and results of \textit{Lemma~\ref{Lemma:StructuralProperties:WitteResults}}, \textit{Lemma~\ref{Lemma:StructuralProperties:P=[S,R]}}, \textit{Lemma~\ref{Lemma:StructuralProperties:Q_(nc)=R^dCompatibleWithGamma}};
        \item and the notations and arguments following \textit{Lemma~\ref{Lemma:StructuralProperties:Q_(nc)=R^dCompatibleWithGamma}}.
    \end{itemize}
    
\end{Remark}
\begin{Remark}\label{Remark:StructuralProperties:SplittingImpliesLooselyKronecker}
    On the alternative proof mentioned in \textit{Subsection \ref{subSection:Splitting<=>LooselyKronecker:=>}}, one can show that in the case $U_t$ is splitting, $\Pbf=[\Sbf,\Rbf]$ is in fact trivial. The idea is to show $\Sbf$ acts trivially on $\Pbf$ using the ergodicity of the $\Hbf$-action on $\Sbf\Pbf/\Gamma_{\Sbf\Pbf}$ and the facts that $\Sbf$ can be considered as $\Hbf\times\Sbf'$ and $\Hbf$ acts trivially on $\Pbf$ (by the splitting assumption). Then, the decomposition provided by (\ref{eq:StructuralProperties:(x,yGamma_Q)mapstoxyGamma}) and (\ref{eq:StructuralProperties:G-action}) implies that the ergodic flow $(\Gbf/\Gamma,\mu,U_t)$ is a compact extension of $(\Sbf/\Gamma_\Sbf,\mu_\Sbf,U_t)$ (as $\Pbf$ is trivial). Since $\Sbf$ is semisimple, we can apply \cite{MR4255058}*{Corollary 1.2} to deduce that $(\Sbf/\Gamma_\Sbf,\mu_\Sbf,U_t)$ is loosely Kronecker. This and \cite{MR653094}*{Section 7, Theorem 7.3} then shows that $(\Gbf/\Gamma,\mu,U_t)$ is loosely Kronecker.
\end{Remark}

\section{Spectral gap condition and an \texorpdfstring{$\SL_2(\R)$}{}-ergodic theorem}\label{Section:SpectralGapCondition}
In this section, continuing with the setting in \textit{Section~\ref{Section:StructureOfErgodicUnipotentFlow}}, we rigorously define the spectral gap condition mentioned in \textit{Section~\ref{Section:Introduction}} and establish an \textit{$\SL_2(\R)$-ergodic theorem} for $\Sbf_{\operatorname{c}}\Qbf$-invariant functions, which extends \cite{lindenstrauss2025timechangerigidityunipotent}*{Corollary 3.15}. This is one new idea we need in this paper to conduct the renormalization step in \textit{Section~\ref{Section:Renormalization}}, as mentioned in \textit{Section~\ref{Section:Introduction}}. Most of the arguments in this section are centered on the theory of unitary representations. For background and notations, we refer to \cites{MR776417,MR4995524}.

\medskip

In addition to the setting in \textit{Section~\ref{Section:StructureOfErgodicUnipotentFlow}}, we denote $\Sbf_{\nc}$ the noncompact factor of $\Sbf$, i.e. the almost direct product of noncompact simple factors, $\Sbf_{\operatorname{c}}$ the compact factor of $\Sbf$, i.e. the almost direct product of compact simple factors, and $A_t=\exp(t\ba)$, where $\{\bu,\ba,\bou\}$ is an $\slfrak_2$-triple corresponding to $U_t=\exp(t\bu)$ (cf. \textit{Subsection \ref{subSection:Preliminaries:HorocyclicUnipotentFlow}}). 

Now, we state the \textit{$\SL_2(\R)$-ergodic theorem} for $\Sbf_{\operatorname{c}}\Qbf$-invariant functions as follows:
\begin{Proposition}[$\SL_2(\R)$-ergodic theorem]\label{Prop:SpectralGapCondition:SL_2-ergodicTheorem}

    For any $\epsilon\in(0,1)$ and any $\Sbf_{\operatorname{c}}\Qbf$-invariant $f\in L^2(\Gbf/\Gamma,\mu)$, $\mu$-a.e. $x\in\Gbf/\Gamma$ satisfies
    \[
    \lim_{n\to\infty}\frac{1}{m_\Hbf(B^{\Hbf,\norm{\cdot}}_\epsilon)}\int_{B^{\Hbf,\norm{\cdot}}_\epsilon}f(A_nhx)dm_\Hbf(h)=\int_{\Gbf/\Gamma}fd\mu,
    \]
    where $B^{\Hbf,\norm{\cdot}}_\epsilon=\{h\in\Hbf:\norm{h-e_\Gbf}<\epsilon\}$, $m_\Hbf$ is a Haar measure on $\Hbf$.
    
\end{Proposition}
Recall that from \textit{Section~\ref{Section:StructureOfErgodicUnipotentFlow}}, we have the decomposition $\Gbf=\Sbf\Pbf\Qbf$ in \eqref{eq:StructuralProperties:G=SPQ} on p.~\pageref{eq:StructuralProperties:G=SPQ}, the identification between $(\Gbf/\Gamma,\mu)$ and $(F_{\Sbf\Pbf}\times(\Qbf/\Gamma_{\Qbf}),m_{\Sbf\Pbf}\times\mu_\Qbf)$ via (\ref{eq:StructuralProperties:(x,yGamma_Q)mapstoxyGamma}), and the $\Gbf$-action on $F_{\Sbf\Pbf}\times(\Qbf/\Gamma_\Qbf)$ given in (\ref{eq:StructuralProperties:G-action}). These allow us to identify the space of $\Sbf_{\operatorname{c}}\Qbf$-invariant functions in $L^2(\Gbf/\Gamma,\mu)$ with the space of $\Sbf_{\operatorname{c}}$-invariant functions in $L^2(F_{\Sbf\Pbf},m_{\Sbf\Pbf})$, which can be further identified with the space of $\Sbf_{\operatorname{c}}$-invariant functions in $L^2(\Sbf\Pbf/\Gamma_{\Sbf\Pbf},\mu_{\Sbf\Pbf})$. Here, $\mu_{\Sbf\Pbf}$ is the Haar probability measure on $\Sbf\Pbf/\Gamma_{\Sbf\Pbf}$. Thus, for a given $\Sbf_{\operatorname{c}}\Qbf$-invariant function $f\in L^2(\Gbf/\Gamma,\mu)$, there exists a $\Sbf_{\operatorname{c}}$-invariant function $f_{\Sbf\Pbf}\in L^2(\Sbf\Pbf/\Gamma_{\Sbf\Pbf},\mu_{\Sbf\Pbf})$ such that
\begin{equation}\label{eq:SpectralGapCondition:int_(G/Gamma)f=int_(SP/Gamma_(SP))f_(SP)}
    \int_{\Gbf/\Gamma}fd\mu=\int_{\Sbf\Pbf/\Gamma_{\Sbf\Pbf}}f_{\Sbf\Pbf}d\mu_{\Sbf\Pbf},
\end{equation}
and since $\Hbf\subseteq\Sbf_{\nc}\subseteq\Sbf$, $\mu$-a.e. $g\Gamma\in\Gbf/\Gamma$ satisfy
\begin{equation}\label{eq:SpectralGapCondition:int_(B^H_epsilon)f=int_(B^H_epsilon)f_(SP)}
    \int_{B^{\Hbf,\norm{\cdot}}_\epsilon}f(A_nhg\Gamma)dm_{\Hbf}(h)=\int_{B^{\Hbf,\norm{\cdot}}_\epsilon}f_{\Sbf\Pbf}(A_nhs\Gamma_{\Sbf\Pbf})dm_\Hbf(h),    
\end{equation}
where $s\in\Sbf\Pbf$ satisfies $g=sq$ for some $q\in\Qbf$. 

Let $\Hcal_1$ be the space of $\Sbf_{\operatorname{c}}$-invariant functions in $L^2_0(\Sbf\Pbf/\Gamma_{\Sbf\Pbf},\mu_{\Sbf\Pbf})$. As $\Sbf_{\operatorname{c}}$ and $\Sbf_{\nc}$ commute, we can consider $\Hcal_1$ as a unitary $\Sbf_{\nc}$-representation by
\[
(s\cdot f)(x)=f(s^{-1}x),\qquad\forall x\in\Sbf\Pbf/\Gamma_{\Sbf\Pbf},\quad f\in\Hcal_1,\quad s\in\Sbf_{\nc}.
\]

We claim that to prove \textit{Proposition~\ref{Prop:SpectralGapCondition:SL_2-ergodicTheorem}}, it is enough to prove the following:
\begin{Lemma}\label{Lemma:SpectralGapCondition:SpectralGapOnSP}

    $\mathbbm{1}_\Hbf\nprec\Hcal_1\vert_{\Hbf}.$
    
\end{Lemma}
Indeed, since $(\Sbf\Pbf/\Gamma_{\Sbf\Pbf},\mu_{\Sbf\Pbf},U_t)$ is a factor of the ergodic $\slfrak_2$-induced flow $(\Gbf/\Gamma,\mu,U_t)$, it is also an ergodic $\slfrak_2$-induced unipotent flow. This implies that the $\Hbf$-action on $\Hcal_1$ is ergodic. Hence, once we prove this lemma, by the arguments in the proof of \cite{lindenstrauss2025timechangerigidityunipotent}*{Corollary 3.15}, it follows that for $\mu_{\Sbf\Pbf}$-a.e. $s\Gamma_{\Sbf\Pbf}\in\Sbf\Pbf/\Gamma_{\Sbf\Pbf}$,
\[
\lim_{n\to\infty}\frac{1}{m_\Hbf(B^{\Hbf,\norm{\cdot}}_\epsilon)}\int_{B^{\Hbf,\norm{\cdot}}_\epsilon}f_{\Sbf\Pbf}(A_nhs\Gamma_{\Sbf\Pbf})dm_\Hbf(h)=\int_{\Sbf\Pbf/\Gamma_{\Sbf\Pbf}}f_{\Sbf\Pbf}d\mu_{\Sbf\Pbf}.
\]
This, (\ref{eq:SpectralGapCondition:int_(G/Gamma)f=int_(SP/Gamma_(SP))f_(SP)}), and (\ref{eq:SpectralGapCondition:int_(B^H_epsilon)f=int_(B^H_epsilon)f_(SP)}) imply that $\mu$-a.e. $g\Gamma\in\Gbf/\Gamma$ satisfy
\[
\lim_{n\to\infty}\frac{1}{m_\Hbf(B^{\Hbf,\norm{\cdot}}_\epsilon)}\int_{B^{\Hbf,\norm{\cdot}}_\epsilon}f(A_nhx)dm_\Hbf(h)=\int_{\Gbf/\Gamma}fd\mu,
\]
proving \textit{Proposition~\ref{Prop:SpectralGapCondition:SL_2-ergodicTheorem}}.

\medskip

For the remainder of this section, we prove \textit{Lemma~\ref{Lemma:SpectralGapCondition:SpectralGapOnSP}}.
\begin{proof}[Proof of \textit{Lemma~\ref{Lemma:SpectralGapCondition:SpectralGapOnSP}}]
By construction, $\Pbf$ is the nilradical of $\Sbf\Pbf$, so it follows that $\Gamma_{\Sbf\Pbf}\cap\Pbf$ is a lattice in $\Pbf$ (cf. \textit{Lemma~\ref{Lemma:StructuralProperties:WitteResults}-\ref{item:WitteResults:GammacapNilG/Z(G)islattice}}). And since $\Sbf_{\operatorname{c}}$ is compact, $\Gamma_{\Sbf\Pbf}\cap\Sbf_{\operatorname{c}}\Pbf$ is a lattice in $\Sbf_{\operatorname{c}}\Pbf$. This let us consider $\Gamma_{\Sbf_{\nc}'}$ the lattice in $\Sbf_{\nc}'=\Sbf_{\nc}/(\Sbf_{\nc}\cap\Sbf_{\operatorname{c}})$ given as the image of $\Gamma_{\Sbf\Pbf}$ under the natural map
\[
\Sbf\Pbf\to\Sbf_{\nc}':s_{\nc}s_{\operatorname{c}}p\mapsto s_{\nc}(\Sbf_{\nc}\cap\Sbf_{\operatorname{c}}),\qquad\forall s_{\nc}\in\Sbf_{\nc},\quad s_{\operatorname{c}}\in\Sbf_{\operatorname{c}},\quad p\in\Pbf.
\]
And by an argument similar to that for identifying the subspace of $\Qbf$-invariant functions in $L^2(\Gbf/\Gamma,\mu)$ with $L^2(\Sbf\Pbf/\Gamma_{\Sbf\Pbf},\mu_{\Sbf\Pbf})$ above, we can identify the subspace $\Hcal_2$ of $\Pbf$-invariant functions in $\Hcal_1$, i.e. the space of $\Sbf_{\operatorname{c}}\Pbf$-invariant functions in $L^2_0(\Sbf\Pbf/\Gamma_{\Sbf\Pbf},\mu_{\Sbf\Pbf})$, with $L^2_0(\Sbf_{\nc}'/\Gamma_{\Sbf_{\nc}'},\mu_{\Sbf_{\nc}'})$. Here, $\mu_{\Sbf_{\nc}'}$ is the Haar probability measure on $\Sbf_{\nc}'/\Gamma_{\Sbf_{\nc}'}$.

Since $\Sbf_{\nc}'$ is semisimple without compact factors and $\Hbf\subseteq\Sbf_{\nc}$ is locally isomorphic to $\SL_2(\R)$, it follows from \cite{MR1719827}*{Lemma 3.1}, the arguments in \cite{MR2486935}*{p.~284,~285}, \cite{MR1779896}*{Theorem C}, and \cite{MR4995524}*{Theorem~9.80} that
\[
\mathbbm{1}_\Hbf\nprec L^2_0(\Sbf_{\nc}'/\Gamma_{\Sbf_{\nc}'},\mu_{\Sbf_{\nc}'})\vert_\Hbf,\qquad\text{i.e.}\quad\mathbbm{1}_\Hbf\nprec\Hcal_2\vert_\Hbf.
\]
This gives us \textit{Lemma~\ref{Lemma:SpectralGapCondition:SpectralGapOnSP}} in the case $\Pbf=\{e_\Pbf\}$, since then $\Hcal_2=\Hcal_1$. And in the case $\Pbf\neq\{e_\Pbf\}$, the proof of \textit{Lemma~\ref{Lemma:SpectralGapCondition:SpectralGapOnSP}} is reduced to showing
\begin{equation}\label{eq:SpectralGapCondition:1_HnprecHcal_2^perp}
    \mathbbm{1}_\Hbf\nprec\Hcal_2^\perp\vert_\Hbf,
\end{equation}
where $\Hcal_2^\perp$ is the orthogonal complement of $\Hcal_2$ in $\Hcal_1$.

Let $\Hcal$ be the space of $\Pbf$-invariant functions in $L^2_0(\Sbf\Pbf/\Gamma_{\Sbf\Pbf},\mu_{\Sbf\Pbf})$ and $\Hcal^\perp$ its orthogonal complement in $L^2_0(\Sbf\Pbf/\Gamma_{\Sbf\Pbf},\mu_{\Sbf\Pbf})$. It is clear that $\Hcal_2\subseteq\Hcal$, and since $\Pbf\lhd\Sbf\Pbf$, $\Hcal$ and $\Hcal^\perp$ are $\Sbf_{\operatorname{c}}$-invariant. With these, we deduce that $\Hcal_2^\perp$ is a subspace of $\Hcal^\perp$. Hence, to obtain \eqref{eq:SpectralGapCondition:1_HnprecHcal_2^perp}, it suffices to show that
\begin{equation}\label{eq:SpectralGapCondition:1_HnprecHcal^perp}
    \mathbbm{1}_\Hbf\nprec\Hcal^\perp\vert_\Hbf,
\end{equation}

Our idea of showing \eqref{eq:SpectralGapCondition:1_HnprecHcal^perp} consists of two steps. The first step is to apply \cite{MR3054613}*{Theorem 1.2} to $(\Dbf,\Ebf)$, where $\Dbf$ is the normal closure of $\Hbf$ in $\Sbf$ and $\Ebf=[\Dbf,\Pbf]$, to obtain $\mathbbm{1}_\Dbf\nprec\Hcal^\perp\vert_\Dbf$. Then, in the second step, we argue using \cite{MR1779896}*{Theorem C} to deduce (\ref{eq:SpectralGapCondition:1_HnprecHcal^perp}).

\medskip

We begin with the first step. As mentioned, we want to apply \cite{MR3054613}*{Theorem 1.2} to $(\Dbf,\Ebf)$. In particular, this result shows $(\Dbf,\Ebf)$ has strong relative property (T), i.e. for any unitary representation $\rho$ of $\Dbf\ltimes\Ebf$, if $\mathbbm{1}_\Dbf\prec\rho\vert_\Dbf$, then there exist nontrivial $\Ebf$-invariant vectors (cf. \cite{MR3054613}*{p.~10}), provided that $\Dbf$ is a connected semisimple Lie group without compact factors and $[\Dbf,\Ebf]=\Ebf$ (cf. \cite{MR3054613}*{Theorem 1.2}). Since the first item follows from $\Dbf$'s construction and the fact that $\Hbf$ is connected, locally isomorphic to $\SL_2(\R)$, and the second item follows from an argument similar to that in the proof of \textit{Lemma~\ref{Lemma:StructuralProperties:P=[S,R]}-\ref{item:P=[S,R]:2}}, to show $\mathbbm{1}_\Dbf\nprec\Hcal^\perp\vert_\Dbf$, we only need to show that $\Hcal^\perp$ has no nontrivial $\Ebf$-invariant vectors.

\medskip

We begin the proof of the fact $\Hcal^\perp$ has no nontrivial $\Ebf$-invariant vectors with the following claim:
\begin{Claim}\label{Claim:SpectralGapCondition:PsubsetEGamma_(SP)}

    $\Pbf\subseteq\overline{\Ebf\Gamma_{\Sbf\Pbf}}$.
    
\end{Claim}
\begin{proof}
By construction, we can directly check $\Ebf\lhd\Sbf\Pbf$. Hence, $\overline{\Ebf\Gamma_{\Sbf\Pbf}}$ is a closed subgroup of $\Sbf\Pbf$. Let $F=\overline{\Ebf\Gamma_{\Sbf\Pbf}}\cap\Pbf$, it follows that $F$ is a closed subgroup of $\Pbf$ and $\Ebf$, $\Gamma_\Pbf=\Gamma_{\Sbf\Pbf}\cap\Pbf$ are closed subgroups of $F$. Let $F^\circ$ be the identity component of $F$, we have $\Gamma_\Pbf\cap F^\circ$ is a lattice in $F^\circ$ because $F/\Gamma_\Pbf$ is compact.

Let $V=\Pbf/[\Pbf,\Pbf]$ and $q:\Pbf\to V$ the natural map. We have $\Gamma_\Pbf\cap[\Pbf,\Pbf]$ is a lattice in $[\Pbf,\Pbf]$, $\Gamma_V=q(\Gamma_\Pbf)$ is a lattice in $V$, and since $\Gamma_\Pbf\cap F^\circ$ is a lattice in $F^\circ$ (shown above), it follows that with $W=q(F^\circ)$, $\Gamma_V\cap W$ is a lattice in $W$. This then implies via \cite{onishchik2000discrete}*{I.~1~Theorem~4.3,~5.7} that $\Gamma_{V/W}=\Gamma_VW/W$ is a lattice in $V/W$.

The $\Sbf$-action on $\Pbf$ by conjugation naturally induces an $\Sbf$-action on $V$. We claim that $W$ is $\Sbf$-invariant. Indeed, for $d\in\Dbf$ and $f\in F^\circ$, we have
\begin{equation}\label{eq:SpectralGapCondition:d.f[P,P]}
    d\cdot f[\Pbf,\Pbf]=dfd^{-1}[\Pbf,\Pbf]=(f[\Pbf,\Pbf])([d,f^{-1}]^{-1}[\Pbf,\Pbf])\in q(F^\circ)=W,
\end{equation}
since $[\Dbf,\Pbf]=\Ebf\subseteq F^\circ$. And for $\gamma\in\Gamma_\Sbf$ and $f\in F^\circ$, let $\sigma\in\Pbf$ be such that $\gamma^{-1}\sigma\in\Gamma_{\Sbf\Pbf}$, we have
\begin{equation}\label{eq:SpectralGapCondition:gamma.f[P,P]}
    \gamma\cdot\ f[\Pbf,\Pbf]=\gamma f\gamma^{-1}[\Pbf,\Pbf]=(\gamma^{-1}\sigma)^{-1}f(\gamma^{-1}\sigma)[\Pbf,\Pbf]\in q(F^\circ)=W,
\end{equation}
since $V$ is abelian and conjugation by $(\gamma^{-1}\sigma)^{-1}$ maps $F^\circ$ to itself. Hence, $W$ is $\Dbf\Gamma_\Sbf$-invariant. And since $\overline{\Dbf\Gamma_\Sbf}=\Sbf$, it follows that $W$ is $\Sbf$-invariant.

We can now consider the induced $\Sbf$-action on $V/W$. We claim that $\Gamma_{V/W}$ is $\Sbf$-invariant. Indeed, by computations as in \eqref{eq:SpectralGapCondition:d.f[P,P]} and \eqref{eq:SpectralGapCondition:gamma.f[P,P]}, it follows that $\Dbf$ acts trivially on $V/W$ and $\Gamma_V$ is $\Gamma_\Sbf$-invariant. Hence, $\Gamma_{V/W}$ is $\Dbf\Gamma_\Sbf$-invariant. And as $\overline{\Dbf\Gamma_\Sbf}=\Sbf$, we deduce that $\Gamma_{V/W}$ is $\Sbf$-invariant.

Now, we have the $\Sbf$-action on $V/W$ preserves $\Gamma_{V/W}$, which is a lattice in $V/W$. This and the fact that $\Sbf$ is connected imply that $\Sbf$ acts trivially on $V/W$. Hence, $[\Sbf,V]\subseteq W$. And since $[\Sbf,\Pbf]=\Pbf$ (cf. \textit{Lemma~\ref{Lemma:StructuralProperties:P=[S,R]}-\ref{item:P=[S,R]:2}}), it follows that $V=[\Sbf,V]\subseteq W$. Thus, $V=W=q(F^\circ)$.

Let $\ffrak$ be the Lie algebra of $F^\circ$, from $q(F^\circ)=V=\Pbf/[\Pbf,\Pbf]$, it follows that $\pfrak=\ffrak+[\pfrak,\pfrak]$. With this and the fact that $\pfrak$ is nilpotent, an induction argument gives us $\pfrak=\ffrak$, which implies that $\Pbf=F^\circ\subseteq\overline{\Ebf\Gamma_{\Sbf\Pbf}}\cap\Pbf$. This finishes the proof of \textit{Claim \ref{Claim:SpectralGapCondition:PsubsetEGamma_(SP)}}.
\end{proof}
Now, we show that $\Hcal^\perp$ has no nontrivial $\Ebf$-invariant vectors as follows.

Since $\Hcal$ is the subspace of $\Pbf$-invariant functions in $L^2_0(\Sbf\Pbf/\Gamma_{\Sbf\Pbf},\mu_{\Sbf\Pbf})$, it is enough to show that if $f\in\Hcal^{\perp}$ is $\Ebf$-invariant, then $f$ is $\Pbf$-invariant.

Consider an $\Ebf$-invariant function $f\in L^2_0(\Sbf\Pbf/\Gamma_{\Sbf\Pbf},\mu_{\Sbf\Pbf})$ and lift $f$ to $\Tilde{f}=f\circ\pi\in L^2_{\operatorname{loc}}(\Sbf\Pbf,m_{\Sbf\Pbf})$, where $\pi:\Sbf\Pbf\to\Sbf\Pbf/\Gamma_{\Sbf\Pbf}$ is the natural map, it follows that
\begin{equation}\label{eq:SpectralGapCondition:Tilde(f)gamma=Tilde(f)}
    \text{for }\gamma\in\Gamma_{\Sbf\Pbf},\qquad\Tilde{f}(g\gamma)=\Tilde{f}(g), \quad m_{\Sbf\Pbf}\text{-a.e. }g\in\Sbf\Pbf.
\end{equation} 
For $x\in\Sbf\Pbf$, denote $L_x$ the left multiplication by $x$ on $\Sbf\Pbf$ and $\Sbf\Pbf/\Gamma_{\Sbf\Pbf}$. Since $f$ is $\Ebf$-invariant, $\Tilde{f}\circ L_x=\Tilde{f}$ for all $x\in\Ebf$. This and the fact $\Ebf\lhd\Sbf\Pbf$ imply that
\[
\text{for }m_\Ebf\text{-a.e. } x\in\Ebf,\qquad \Tilde{f}\circ R_x=\Tilde{f}\, (m_{\Sbf\Pbf}\text{-a.e.}),
\] 
where $R_x$ is the right multiplication by $x$ on $\Sbf\Pbf$. Since the map
\[
\Sbf\Pbf\to L^1_{\operatorname{loc}}(\Sbf\Pbf,m_{\Sbf\Pbf}):x\mapsto \Tilde{f}\circ R_x
\]
is continuous, we deduce that $\Tilde{f}\circ R_x=\Tilde{f}$ for all $x\in\Ebf$. This and (\ref{eq:SpectralGapCondition:Tilde(f)gamma=Tilde(f)}) imply that $\Tilde{f}\circ R_x=\Tilde{f}$ for all $x\in\Ebf\Gamma_{\Sbf\Pbf}$, and thus for all $x\in\overline{\Ebf\Gamma_{\Sbf\Pbf}}$. Hence, \textit{Claim \ref{Claim:SpectralGapCondition:PsubsetEGamma_(SP)}} implies that 
\[
\Tilde{f}\circ R_x=\Tilde{f},\qquad \forall x\in\Pbf.
\]
And as $\Pbf\lhd\Sbf\Pbf$, we get $\Tilde{f}\circ L_x=\Tilde{f}$ for all $x\in\Pbf$, giving us $f$ is $\Pbf$-invariant.

Therefore, we conclude the first step and obtain $\mathbbm{1}_\Dbf\nprec\Hcal^\perp\vert_\Dbf$.

\medskip

We continue with the second step. Denote $\rho$ the $\Dbf$-representation $\Hcal^\perp\vert_\Dbf$, we argue using \cite{MR1779896}*{Theorem C} to deduce (\ref{eq:SpectralGapCondition:1_HnprecHcal^perp}) from the fact $\mathbbm{1}_\Dbf\nprec\rho$ shown above as follows.

By construction, $\Dbf=\Sbf_1\cdots\Sbf_k$ as an almost direct product of simple factors of $\Sbf$, where each simple factor $\Sbf_j$ satisfies $\hfrak=\Lie(\Hbf)$ has nontrivial projection on $\sfrak_j=\Lie(\Sbf_j)$. Denote $\Tilde{\Dbf}=\Sbf_1\times\cdots\times\Sbf_k$ and
\[
\varpi:\Tilde{\Dbf}\to\Dbf:(s_1,\dots,s_k)\mapsto s_1\cdots s_k,
\]
it follows that $\varpi$ is a finite covering homomorphism. Then, $\Tilde{\Hbf}=\varpi^{-1}(\Hbf)^\circ$ is a connected subgroup of $\Tilde{\Dbf}$ locally isomorphic to $\SL_2(\R)$, whose projection to each simple factor $\Sbf_j$ is also locally isomorphic to $\SL_2(\R)$, and the unitary representation $\Tilde{\rho}=\rho\circ\varpi$ of $\Tilde{\Dbf}$ satisfies $\mathbbm{1}_{\Tilde{\Dbf}}\nprec\Tilde{\rho}$ (since $\mathbbm{1}_\Dbf\nprec\rho$). This allows us to apply \cite{MR1779896}*{Theorem C}, deducing that $\mathbbm{1}_{\Tilde{\Hbf}}\nprec\Tilde{\rho}\vert_{\Tilde{\Hbf}}$. As $\varpi\vert_{\Tilde{\Hbf}}:\Tilde{\Hbf}\to\Hbf$ is a finite covering, we obtain (\ref{eq:SpectralGapCondition:1_HnprecHcal^perp}), finishing the second step.

\medskip

Thus, by the previous arguments, from (\ref{eq:SpectralGapCondition:1_HnprecHcal^perp}) we deduce (\ref{eq:SpectralGapCondition:1_HnprecHcal_2^perp}), with which we finish the proof of \textit{Lemma~\ref{Lemma:SpectralGapCondition:SpectralGapOnSP}}.
\end{proof}

\section{Behavior of \texorpdfstring{$\psi$}{} on close points and along \texorpdfstring{$A^{(1)}_t$}{}-orbits}\label{Section:BehaviorOfpsiOnClosePointsAndAlongA_t}
In this section, we recall the behaviors of a Kakutani factor map between two ergodic $\slfrak_2$-induced unipotent flows acting on two close points and on two points on the same orbit under the ``diagonal flow", which are established in \cite{lindenstrauss2025timechangerigidityunipotent}.

\medskip

We recall and fix the following notations and assumptions, which will be present for the rest of the paper, where $i=1,2$:
\begin{enumerate}[label=(\ref{Section:BehaviorOfpsiOnClosePointsAndAlongA_t}:\roman*)]
    \item\label{item:BehaviorOfpsis:1} $\Gbf_i$ is the identity component of a real linear algebraic group, $\Gamma_i<\Gbf_i$ is a lattice, and $\mu_i$ is the Haar probability measure on $\Gbf_i/\Gamma_i$;
    \item\label{item:BehaviorOfpsis:2} $(\Gbf_i/\Gamma_i,\mu_i,U^{(i)}_t)$ is an ergodic $\slfrak_2$-induced unipotent flow;
    \item\label{item:BehaviorOfpsis:3} the Lie subalgebra $\hfrak_i\subseteq\gfrak_i=\Lie(\Gbf_i)$ is the $\slfrak_2(\R)$-subalgebra corresponding to $U^{(i)}_t$, and $\Hbf_i=\langle\exp(\hfrak_i)\rangle$;
    \item\label{item:BehaviorOfpsis:4}  $\bu^{(i)},\ba^{(i)},\bou^{(i)}\in\hfrak_i$ form a $\slfrak_2$-triple with
    \[
    U^{(i)}_t=\exp(t\bu^{(i)}),\quad A^{(i)}_t=\exp(t\ba^{(i)}),\quad \Uh^{(i)}_t=\exp(t\bou^{(i)});
    \]
    \item\label{item:BehaviorOfpsis:5} assume that the flow $U^{(2)}_t$ on $\Gbf_2/\Gamma_2$ is nonsplitting;
    \item\label{item:BehaviorOfpsis:6} with respect to $\{\bu^{(i)},\ba^{(i)},\bou^{(i)}\}$, fix a chain basis for $\gfrak_{i}$
    \begin{equation}\label{eq:BehaviorOfpsi:ChainBasis}
        \left\{\bu^{(i)},\ba^{(i)},\bou^{(i)},\bx^{(i)}_{0,1},\dots,\bx^{(i)}_{q_1^{(i)},1},\dots,\bx^{(i)}_{0,n^{(i)}},\dots,\bx^{(i)}_{q_{n^{(i)}}^{(i)},n^{(i)}}\right\}
    \end{equation}
    (cf. (\ref{eq:Preliminaries:ChainBasis})), and denote (cf. \textit{Subsection~\ref{subSection:Preliminaries:KakBall}})
    \[
    \Zbf_i=C_{\Gbf_i}(\Hbf_i),\quad\textstyle{\zfrak_i=\oplus_{j\in J_i}\R\bx^{(i)}_{0,j}},\quad J_i=\{1\leq j\leq n^{(i)}:q^{(i)}_j=0\},\quad i=1,2;
    \] 
    \item\label{item:BehaviorOfpsis:7} with respect to the chain basis (\ref{eq:BehaviorOfpsi:ChainBasis}) of $\gfrak_i$, define $d_{\mathbf{G}_{i}},d_{\mathbf{G}_{i}/\Gamma_i}$ as in \textit{Subsection \ref{subSection:Preliminaries:ChainBasisAndMetric}}, and the Kakutani-Bowen balls in $\mathbf{G}_{i}$ as in \textit{Definition \ref{Def:Preliminaries:Kakutani-BowenBalls}}, $i=1,2$;
    \item\label{item:BehaviorOfpsis:8}$\psi:\mathbf{G}_{1}/\Gamma_1\to\mathbf{G}_{2}/\Gamma_2$ is a Kakutani factor map between the flows $(\mathbf{G}_{i}/\Gamma_i,\mu_i,U^{(i)}_t)$, $i=1,2$, whose time change $\tau\in L^1_+(\Gbf_1/\Gamma_1,\mu_1)$ is $10^{-20}\eold{epsilon:Preliminaries:CoveringLemma}$-good (cf. \textit{Subsection~\ref{subSection:Preliminaries:epsilon-goodTimeChange}}), and $\eold{epsilon:Preliminaries:CoveringLemma}$ given on p.~\pageref{epsilon:Preliminaries:CoveringLemma} for $(\Gbf_2/\Gamma_2,\mu_2,U^{(2)}_t)$;
    \item\label{item:BehaviorOfpsis:9}$\Znew\label{Z:BehaviorOfpsi:GoodSetWhereTheConjugationRelationHold}$ is a full measure $U^{(1)}_t$-invariant subset of $\mathbf{G}_{1}/\Gamma_1$ such that $\omega,\rho$ given by
    \begin{equation}\label{eq:BehaviorOfpsi:rho(x,s)andomega(x,s)}
        \omega(x,s)=\int_0^s\tau(U^{(1)}_rx)dr,\quad s=\int_0^{\rho(x,s)}\tau(U^{(1)}_rx)dr,
    \end{equation}
    are well-defined for all $x\in\Zold{Z:BehaviorOfpsi:GoodSetWhereTheConjugationRelationHold}$, $s\in\R$, and
    \begin{equation}\label{eq:BehaviorOfpsi:|partial_svarpi_i(x,s)-1|IsSmall}
        \begin{gathered}
            [\R\to\R:s\mapsto\tau(U^{(1)}_sx)]\in C^{\infty}(\R),\\
            \absolute{\tau(U^{(1)}_sx)-1}<10^{-20}\eold{epsilon:Preliminaries:CoveringLemma},\quad\forall x\in\Zold{Z:BehaviorOfpsi:GoodSetWhereTheConjugationRelationHold},\; s\in\R,
        \end{gathered}
    \end{equation}
    \begin{equation}\label{eq:BehaviorOfpsi:ConjugationRelation}
        \begin{aligned}
        &\psi(U^{(1)}_{\rho(x,s)}x)=U^{(2)}_s\psi(x),\\
        &\psi(U^{(1)}_sx)=U^{(2)}_{\omega(x,s)}\psi(x),\quad\forall x\in\Zold{Z:BehaviorOfpsi:GoodSetWhereTheConjugationRelationHold},\;s\in\R.
        \end{aligned}
    \end{equation}
    We also assume every subset of $\mathbf{G}_{1}/\Gamma_1$ to be contained in $\Zold{Z:BehaviorOfpsi:GoodSetWhereTheConjugationRelationHold}$, if nothing is stated otherwise.
\end{enumerate}

\subsection{Behavior of \texorpdfstring{$\psi$}{} acting on close points}\label{subSection:BehaviorOfpsi:OnClosePoints}
We first recall in this subsection the behavior of $\psi$ action on two close points. Roughly speaking, for two points $x,y$ in a ``good" subset of large measure in $\Gbf_1/\Gamma_1$ such that $x$ lies in a small Kakutani-Bowen ball around $y$, $\psi(x)$ also lies in a small Kakutani-Bowen ball around $\psi(y)$. This is \cite{lindenstrauss2025timechangerigidityunipotent}*{Lemma 8.1}, which we state as follows:
\begin{Proposition}\label{Proposition:BehaviorOfpsi:psi(Kak)=Kak}

    For $\eta\in(0,10^{-3})$, there exist $\denew\label{delta:BehaviorOfpsi:BoundOndeltaInTwoSideMatching}>0,\Rnew\label{R:BehaviorOfpsi:BoundOnRInTwoSideMatching}>1$, and a compact set $\Knew\label{K:BehaviorOfpsi:GoodSetInLemmapsi(Kak)=Kak}\subset \mathbf{G}_{1}/\Gamma_1$ with $\mu_1(\Kold{K:BehaviorOfpsi:GoodSetInLemmapsi(Kak)=Kak})>1-5\eta$ satisfying the following property. For $\delta\in(0,\deold{delta:BehaviorOfpsi:BoundOndeltaInTwoSideMatching}]$, there exists $\enew\label{epsilon:BehaviorOfpsi:distanceBoundOnG_(d1)/Gamma_1}(\delta)\in(0,\delta)$ such that if $x,y\in\Kold{K:BehaviorOfpsi:GoodSetInLemmapsi(Kak)=Kak}$ and $x\in\Kak(R,\eold{epsilon:BehaviorOfpsi:distanceBoundOnG_(d1)/Gamma_1}(\delta),y)$ for $\Rold{R:BehaviorOfpsi:BoundOnRInTwoSideMatching}\leq R\leq\infty$, then
    \[
    \psi(x)\in\Kak(1.9R,\delta,\psi(y)).
    \]
    Moreover, when $R=\infty$, the same conclusion holds, and 
    \[
    c_{\bu}\in[-\delta,\delta],\quad c_{\ba}\in[-\delta,\delta],\quad g\in \overline{B}^{\mathbf{G}_{2}}_{\delta}\cap\log^{-1}\left(\textstyle{\bigoplus_{j=1}^{n^{(2)}}\R\bx^{(2)}_{0,j}}\right)
    \]
    that satisfy $\psi(x)=U^{(2)}_{c_{\bu}}A^{(2)}_{c_{\ba}}g\psi(y)$ are unique in the sense as in \textit{Lemma~\ref{Lemma:BasicLemma:BasicLemma}-\ref{item:BasicLemma:Statement:C}}.
    
\end{Proposition}

The proof of this proposition is essentially a combination of \textit{Lemma~\ref{Lemma:BasicLemma:BasicLemma}} and \textit{Lemma~\ref{Lemma:epsilon-goodTimeChange:PreserveTwoSideMatching}}. The main idea is as follows. For $x$ lies in a small Kakutani-Bowen ball around $y$, we show that $x$ and $y$ are two sides matchable. Then, \textit{Lemma~\ref{Lemma:epsilon-goodTimeChange:PreserveTwoSideMatching}} implies $\psi(x)$ and $\psi(y)$ are two sides matchable. From here, \textit{Lemma~\ref{Lemma:BasicLemma:BasicLemma}} gives us the conclusion of \textit{Proposition~\ref{Proposition:BehaviorOfpsi:psi(Kak)=Kak}}. For more details, we refer to \cite{lindenstrauss2025timechangerigidityunipotent}*{Subsection 8.3}.

\subsection{Behavior of \texorpdfstring{$\psi$}{} along $A^{(1)}_t$-orbits}\label{subSection:BehaviorOfpsi:AlongA_t} In this subsection, we recall the behavior of $\psi$ acting on two points on the same $A^{(1)}_t$-orbit, deducing the existence of an element $g_{\psi}\in\exp\left(\R\bu^{(2)}\oplus\bigoplus_{j\notin J_2}\R\bx^{(2)}_{0,j}\right)$ such that $(g_{\psi}\psi)( A^{(1)}_px)$, $A^{(2)}_p(g_{\psi}\psi)(x)$ are on the same $U^{(2)}_t$-orbit. To better suit our purposes in the latter part of the paper, we state the result in an adjusted form, similar to that of \cite{lindenstrauss2025timechangerigidityunipotent}*{Proposition 10.1}, as follows:
\begin{Proposition}\label{Proposition:BehaviorOfpsiA_t}

    There exists $g_{\psi}\in\exp\left(\R\bu^{(2)}\oplus\bigoplus_{j\notin J_2}\R\bx^{(2)}_{0,j}\right)$ such that $\tilde{\psi}=g_{\psi}\psi:\Gbf_1/\Gamma_1\to\Gbf_2/\Gamma_2$ is a Kakutani factor map between $(\Gbf_i/\Gamma_i,\mu_i,U^{(i)}_t)$, $i=1,2$, whose corresponding time change $\tau$ is $10^{-20}\eold{epsilon:Preliminaries:CoveringLemma}$-good, satisfying
    \[
    \tilde{\psi}(A^{(1)}_1x)=A^{(2)}_1U^{(2)}_{f(x)}\tilde{\psi}(x),\qquad \forall x\in\Zold{Z:BehaviorOfpsiA_t:GoodSetForpsi(A_px)=U_()exp(p(a+X_0))psi(x)},
    \]
    where $\Znew\label{Z:BehaviorOfpsiA_t:GoodSetForpsi(A_px)=U_()exp(p(a+X_0))psi(x)}$ is a full measure $U^{(1)}_t$-invariant set in $\mathbf{G}_{1}/\Gamma_1$, and $f:\Zold{Z:BehaviorOfpsiA_t:GoodSetForpsi(A_px)=U_()exp(p(a+X_0))psi(x)}\to\R$ is a measurable function.
    
\end{Proposition}
In the following, we provide the main ideas for the proof of this result, based on \cite{lindenstrauss2025timechangerigidityunipotent}*{Section 10}. We refer to the same reference for the detailed arguments.

Before we begin, let us define, in addition to \ref{item:BehaviorOfpsis:1}-\ref{item:BehaviorOfpsis:9}, the set
\begin{enumerate}[label=(\ref{subSection:BehaviorOfpsi:AlongA_t}:\roman*)]
    \item\label{item:BehaviorOfpsiA_t:1} $\Znew\label{Z:BehaviorOfpsiUnderA_t:GenericPointsOfGoodSetInLemmapsi(Kak)=Kak}$ is the set of generic points $x$ of $\mathbf{1}_{\Kold{K:BehaviorOfpsi:GoodSetInLemmapsi(Kak)=Kak}}$ w.r.t. $U^{(1)}_t$ in $\Zold{Z:BehaviorOfpsi:GoodSetWhereTheConjugationRelationHold}\subseteq\mathbf{G}_{1}/\Gamma_1$, where $\Kold{K:BehaviorOfpsi:GoodSetInLemmapsi(Kak)=Kak}$ is given in \textit{Proposition \ref{Proposition:BehaviorOfpsi:psi(Kak)=Kak}} for some $\eta\in(0,10^{-3})$, such that $\psi(x)$ is not $U^{(2)}_t$-periodic (i.e. $U^{(2)}_t\psi(x)=\psi(x)$ iff $t=0$).
\end{enumerate}
As the flow $U^{(i)}_t$ on $(\Gbf_i/\Gamma_i,\mu_i)$ is ergodic, $i=1,2$ (cf. \ref{item:BehaviorOfpsis:2}), it follows from the \textit{Pointwise ergodic theorem} and \ref{item:BehaviorOfpsis:8} that $\Zold{Z:BehaviorOfpsiUnderA_t:GenericPointsOfGoodSetInLemmapsi(Kak)=Kak}$ is $U^{(1)}_t$-invariant and of full measure in $\Gbf_1/\Gamma_1$. This set helps us in the following and in \textit{Section \ref{Section:BehaviorOfpsiW_t}}.

\medskip

The proof of \textit{Proposition~\ref{Proposition:BehaviorOfpsiA_t}} is split into four steps.

The first step is to obtain the local relation between $\psi A^{(1)}_p$ and $\psi$. In particular, we apply \textit{Proposition~\ref{Proposition:BehaviorOfpsi:psi(Kak)=Kak}} and \eqref{eq:BehaviorOfpsi:ConjugationRelation} to show that for small $p$,
\begin{equation}\label{eq:BehaviorOfpsi:AlongA_t:LocalRelation}
    \psi(A^{(1)}_px)=U^{(2)}_{c_{\bu}(x,p)}A^{(2)}_{c_{\ba}(p)}g(p)\psi(x),\qquad\mu_1\text{-a.e. }x\in\Gbf_1/\Gamma_1.
\end{equation}
For more details, we refer to \cite{lindenstrauss2025timechangerigidityunipotent}*{Proof of Lemma 10.4}.

In the second step, we show that in \eqref{eq:BehaviorOfpsi:AlongA_t:LocalRelation},
\begin{equation}\label{eq:BehaviorOfpsi:AlongA_t:c_a(p)=p}
    c_{\ba}(p)=p.
\end{equation}
To do so, we apply \eqref{eq:BehaviorOfpsi:AlongA_t:LocalRelation}, \eqref{eq:BehaviorOfpsi:ConjugationRelation}, and \ref{item:BehaviorOfpsiA_t:1} to obtain an equation in $s$ that connects $\omega$, $\rho$ from \eqref{eq:BehaviorOfpsi:rho(x,s)andomega(x,s)} with $c_u(\cdot,\cdot)$, $c_\ba(\cdot)$ from \eqref{eq:BehaviorOfpsi:AlongA_t:LocalRelation}. Divide both sides of this equation and let $s\to\infty$ then give us \eqref{eq:BehaviorOfpsi:AlongA_t:c_a(p)=p}. For more details, see \cite{lindenstrauss2025timechangerigidityunipotent}*{Proof of Proposition 10.3}.

For the third step, with $J_2$ in \ref{item:BehaviorOfpsis:6}, we show that in \eqref{eq:BehaviorOfpsi:AlongA_t:LocalRelation},
\begin{equation}\label{eq:BehaviorOfpsi:AlongA_t:g(p)HasNoZ_2-part}
    g(p)\in\exp\left(\textstyle{\bigoplus_{j\notin J_2}\R\bx^{(2)}_{0,j}}\right).
\end{equation}
To do so, we start from a relation between the flows $U^{(1)}_t$, $A^{(1)}_t$, $\Uh^{(1)}_t$ given by \textit{Lemma~\ref{Lemma:Preliminaries:SL_2Computation}}, then argue using \textit{Proposition~\ref{Proposition:BehaviorOfpsi:psi(Kak)=Kak}}, \eqref{eq:BehaviorOfpsi:AlongA_t:LocalRelation}, \eqref{eq:BehaviorOfpsi:AlongA_t:c_a(p)=p}, and \textit{Lemma~\ref{Lemma:Preliminaries:a_(-p)ga_pForgInKak(R,delta)}}. This shows that with $p(z_n)=\frac{1}{2}\log\left(\frac{\absolute{e^p-1}}{z_n\delta}\right)$ for some small $\delta>0$ and well-chosen small $z_n>0$, $A^{(2)}_{-p(z_n)}g(p)A^{(2)}_{p(z_n)}\xrightarrow[n\to\infty]{}e_{\Gbf_2}$, giving us \eqref{eq:BehaviorOfpsi:AlongA_t:g(p)HasNoZ_2-part}. For more details, see \cite{lindenstrauss2025timechangerigidityunipotent}*{Subsection 10.2}.

Lastly, in the fourth step, we combine the previous 3 steps and the \textit{Baker-Campbell-Hausdorff formula} to activate the \textit{Banach fixed-point theorem}, proving the existence of $g_\psi$ in \textit{Proposition \ref{Proposition:BehaviorOfpsiA_t}}. For more details, see \cite{lindenstrauss2025timechangerigidityunipotent}*{Subsection 10.3}.
\begin{Remark}\label{Remark:BehaviorOfpsi:Assume(g_(psi))=I_(d_2+1)}
    We see that in \textit{Proposition~\ref{Proposition:BehaviorOfpsiA_t}}, if $\psi$ is a Kakutani equivalence, then so is $\Tilde{\psi}=g_\psi\psi$.    
\end{Remark}

\section{Behavior of \texorpdfstring{$\psi$}{} with respect to the compact semisimple factor and the abelian factor}\label{Section:BehaviorOfpsiW_t}
In this section, continuing with the settings in \textit{Section~\ref{Section:BehaviorOfpsiOnClosePointsAndAlongA_t}}, and consider in addition the structural settings in \textit{Section~\ref{Section:StructureOfErgodicUnipotentFlow}}, we study the behavior of $\psi$ with respect to $\Wbf$-the product of the compact semisimple factor and the abelian factor. Our goal is to show that up to $U^{(2)}_t$-translations, $\psi$ is $\Wbf$-equivariant via a continuous homomorphism $\Tilde{\varpi}:\widetilde{\Wbf}\to\Zbf_2=C_{\Gbf_2}(\Hbf_2)$, where $\widetilde{\Wbf}$ is a universal cover of $\Wbf$. This is one new idea of this paper to help with the renormalization procedure in \textit{Section~\ref{Section:Renormalization}}, as mentioned in \textit{Section~\ref{Section:Introduction}}.

\medskip

Before stating the main result of this section, we recall and fix the notations and assumptions \ref{item:BehaviorOfpsis:1}-\ref{item:BehaviorOfpsis:9} in \textit{Section~\ref{Section:BehaviorOfpsiOnClosePointsAndAlongA_t}}, $\Kold{K:BehaviorOfpsi:GoodSetInLemmapsi(Kak)=Kak}$, $\deold{delta:BehaviorOfpsi:BoundOndeltaInTwoSideMatching}$, $\Rold{R:BehaviorOfpsi:BoundOnRInTwoSideMatching}$, $\eold{epsilon:BehaviorOfpsi:distanceBoundOnG_(d1)/Gamma_1}(\cdot)$ in \textit{Proposition~\ref{Proposition:BehaviorOfpsi:psi(Kak)=Kak}} for a fixed $\eta\in(0,10^{-3})$, and \ref{item:BehaviorOfpsiA_t:1} in p.~\pageref{item:BehaviorOfpsiA_t:1}. And based on the settings in \textit{Section~\ref{Section:StructureOfErgodicUnipotentFlow}} and \textit{Proposition~\ref{Proposition:BehaviorOfpsiA_t}}, we add:
\begin{enumerate}[label=(\ref{Section:BehaviorOfpsiW_t}:\roman*)]
    \item\label{item:BehaviorOfpsiW_t:1} There exist a full measure $U^{(1)}_t$-invariant set $\Zold{Z:BehaviorOfpsiA_t:GoodSetForpsi(A_px)=U_()exp(p(a+X_0))psi(x)}$ in $\Gbf_1/\Gamma_1$ and a measurable function $f:\Zold{Z:BehaviorOfpsiA_t:GoodSetForpsi(A_px)=U_()exp(p(a+X_0))psi(x)}\to\R$ such that
    \[
    \psi(A^{(1)}_1x)=A^{(2)}_1U^{(2)}_{f(x)}\psi(x),\qquad\forall x\in\Zold{Z:BehaviorOfpsiA_t:GoodSetForpsi(A_px)=U_()exp(p(a+X_0))psi(x)}.
    \]    
    \item\label{item:BehaviorOfpsiW_t:2} $\Gbf_1=\Sbf_1\Pbf_1\Qbf_{1,\operatorname{c}}\Qbf_{1,\nc}\cong((\Sbf_1\ltimes\Pbf_1)\Qbf_{1,\operatorname{c}})\times\Qbf_{1,\nc}$;
    \item\label{item:BehaviorOfpsiW_t:3} $\Rbf_1=\Pbf_1\Qbf_{1,\operatorname{c}}\Qbf_{1,\nc}\cong\Pbf_1\times\Qbf_{1,\operatorname{c}}\times\Qbf_{1,\nc}$ is the nilradical (and the radical) of $\Gbf_1$;
    \item\label{item:BehaviorOfpsiW_t:4} $\Sbf_1$ is the identity component of a semisimple real linear algebraic group, and a maximal connected semisimple factor of $\Gbf_1$ that contains $\Hbf_1=\langle\exp(\hfrak_1)\rangle$, $\Sbf_{1,\nc}$ its noncompact factor, i.e. the almost direct product of noncompact simple factors, and $\Sbf_{1,\operatorname{c}}$ its compact factor, i.e. the almost direct product of compact simple factors;
    \item\label{item:BehaviorOfpsiW_t:5} $\Qbf_{1,\operatorname{c}}\cong\T^{a_1}$, $\Qbf_{1,\nc}\cong\R^{d_1}$, and 
    \[
    \Qbf_1=\Qbf_{1,\operatorname{c}}\Qbf_{1,\nc}\cong\Qbf_{1,\operatorname{c}}\times\Qbf_{1,\nc}\cong\T^{a_1}\times\R^{d_1}
    \]
    are connected subgroups of $Z(\Gbf_1)$, where $a_1,d_1\geq0$;
    \item\label{item:BehaviorOfpsiW_t:6} $\Sbf_1\Pbf_1\cap\Qbf_1=\Sbf_1\cap\Qbf_{1,\operatorname{c}}$ is finite and central in $\Gbf_1$;
    \item\label{item:BehaviorOfpsiW_t:7} $\Gamma_{\Qbf_1}=\Gamma_1\cap\Qbf_1$ and $\Gamma_{1,\nc}=\Gamma_1\cap\Qbf_{1,\nc}$ are lattices in $\Qbf_1$ and $\Qbf_{1,\nc}$, respectively;
    \item\label{item:BehaviorOfpsiW_t:8} $\Wbf=\Sbf_{1,\operatorname{c}}\Qbf_1$, $F_\Wbf\subseteq\Wbf$ is a precompact Borel fundamental domain for $\Gamma_{1,\nc}$ (which is a lattice in $\Wbf$), and $m_\Wbf$ is a Haar measure on $\Wbf$ such that $m_{\Wbf}(F_\Wbf)=1$.
    \item\label{item:BehaviorOfpsiW_t:9} $\widetilde{\Wbf}$ is a universal cover of $\Wbf$, $\pi_\Wbf:\widetilde{\Wbf}\to\Wbf$ is the natural map, $\Gamma_{\widetilde{\Wbf}}=\pi_\Wbf^{-1}(\Gamma_{1,\nc})$, $F_{\widetilde{\Wbf}}\subseteq\widetilde{\Wbf}$ is a precompact Borel fundamental domain for $\Gamma_{\widetilde{\Wbf}}$ (which is a lattice in $\widetilde{\Wbf}$), and $m_{\widetilde{\Wbf}}$ is a Haar measure on $\widetilde{\Wbf}$ such that $m_{\widetilde{\Wbf}}(F_{\widetilde{\Wbf}})=1$.
\end{enumerate}

Now, we state the main result of this section as follows:
\begin{Proposition}\label{Prop:BehaviorOfpsiW_t}

There exist a full measure $U^{(1)}_t$ and $A_1^{(1)}$-invariant set $\Znew\label{Z:BehaviorOfpsiW_t:ResultOfInduction}\subseteq\Gbf_1/\Gamma_1$ and a Kakutani factor map $\Tilde{\psi}:\Gbf_1/\Gamma_1\to\Gbf_2/\Gamma_2$ between two systems $(\Gbf_i/\Gamma_i,\mu_i,U^{(i)}_t)$, $i=1,2$, whose time change $\Tilde{\tau}:\Zold{Z:BehaviorOfpsiW_t:ResultOfInduction}\to\R$ is $10^{-20}\eold{epsilon:Preliminaries:CoveringLemma}$-good, such that for $x\in\Zold{Z:BehaviorOfpsiW_t:ResultOfInduction}$,
    \begin{equation*}
        \begin{gathered}
            \textstyle{\Tilde{\psi}\left(gx\right)=\Tilde{\varpi}(\Tilde{g})\Tilde{\psi}(x),\qquad \forall g\in\Wbf,\quad\Tilde{g}\in\pi_{\Wbf}^{-1}(g);}\\
            \Tilde{\psi}(x)=U^{(2)}_{\Tilde{b}(x)}\psi(x);\qquad \Tilde{\psi}(A^{(1)}_1x)=A^{(2)}_1U^{(2)}_{\Tilde{f}(x)}\Tilde{\psi}(x);
        \end{gathered}
    \end{equation*}
    where $\Tilde{\varpi}:\widetilde{\Wbf}\to\Zbf_2=C_{\Gbf_2}(\Hbf_2)$ is a continuous homomorphism such that
    \[
    \Tilde{\varpi}(\gamma)\in\Gamma_2\cap Z(\Gbf_2),\qquad\forall\gamma\in\Gamma_{\widetilde{\Wbf}};
    \]
    and $\Tilde{b}, \Tilde{f}:\Zold{Z:BehaviorOfpsiW_t:ResultOfInduction}\to\R$ are some measurable functions.
\end{Proposition}
Notice that when $\Wbf=\Sbf_{1,\operatorname{c}}\Qbf_1$ is trivial, all the requirements in \textit{Proposition~\ref{Prop:BehaviorOfpsiW_t}} are automatically satisfied with $\Tilde{\psi}=\psi$, $\tilde{b}\equiv0$, $\tilde{f}=f$ and $\Zold{Z:BehaviorOfpsiW_t:ResultOfInduction}=\bigcap_{n\in\Z}A_n^{(1)}(\Zold{Z:BehaviorOfpsi:GoodSetWhereTheConjugationRelationHold}\cap \Zold{Z:BehaviorOfpsiA_t:GoodSetForpsi(A_px)=U_()exp(p(a+X_0))psi(x)})$. Thus, for the remainder of this section, let us assume $\Wbf$ is nontrivial.

\medskip

\textit{Idea of proof.} Our proof consists of two steps. First, we study the behavior of $\psi$ acting on two points that differ by some element in $\Wbf$, using a strategy similar to that in the sketch proof of \textit{Proposition \ref{Proposition:BehaviorOfpsiA_t}}. For the second step, we use behavior of $\psi$ in the first step to obtain a cohomological equation, which we solve and deduce that \textit{Proposition~\ref{Prop:BehaviorOfpsiW_t}} holds when the $\epsilon$-good condition of time change is ignored. We then use convolution to recover this property.

\medskip

The latter part of this section is dedicated to the details arguments following the \textit{Idea of proof}.

\subsection{Behavior of \texorpdfstring{$\psi$}{} with respect to \texorpdfstring{$\Wbf$}{}}\label{subSection:BehaviorOfpsiW_t:psi(W_tx)=...psi(x)}
In this subsection, we employ a strategy similar to that in the sketch proof of \textit{Proposition \ref{Proposition:BehaviorOfpsiA_t}} to deduce the relation between $\psi(g(\cdot))$ and $\psi$ for $g\in\Wbf$. The detailed statement is as follows:
\begin{Proposition}\label{Prop:BehaviorOfpsiW_t:psi(W_tx)=...psi(x)}

    There exist full measure $U^{(1)}_t$-invariant sets $\Znew\label{Z:BehaviorOfpsiW_t:psi(W_tx)=U_(f_j(x,s))g(s)psi(x)}(g)$ in $\Gbf_1/\Gamma_1$ for $g\in\Wbf$ and $\Tilde{\varpi}:\widetilde{\Wbf}\to \Zbf_2=C_{\Gbf_2}(\Hbf_2)$ a continuous homomorphism such that 
    \[
    \Tilde{\varpi}(\gamma)\in\Gamma_2\cap Z(\Gbf_2),\qquad\forall\gamma\in\Gamma_{\widetilde{W}};
    \]
    $\{(x,g):x\in\Zold{Z:BehaviorOfpsiW_t:psi(W_tx)=U_(f_j(x,s))g(s)psi(x)}(g),g\in\Wbf\}$ is of full measure in $(\Gbf_1/\Gamma_1)\times\Wbf$; 
    and for every $g\in\Wbf$,
    \begin{equation*}
        \begin{gathered}
        \Zold{Z:BehaviorOfpsiW_t:psi(W_tx)=U_(f_j(x,s))g(s)psi(x)}(g)\subseteq\Zold{Z:BehaviorOfpsiUnderA_t:GenericPointsOfGoodSetInLemmapsi(Kak)=Kak}\cap g^{-1}\Zold{Z:BehaviorOfpsiUnderA_t:GenericPointsOfGoodSetInLemmapsi(Kak)=Kak};\qquad\Zold{Z:BehaviorOfpsiW_t:psi(W_tx)=U_(f_j(x,s))g(s)psi(x)}(g)=\Zold{Z:BehaviorOfpsiW_t:psi(W_tx)=U_(f_j(x,s))g(s)psi(x)}(g\gamma)\subseteq\Zold{Z:BehaviorOfpsiW_t:psi(W_tx)=U_(f_j(x,s))g(s)psi(x)}(e_{\Gbf_1}),\quad\forall\gamma\in\Gamma_{1,\nc};\\
        \psi(gx)=U^{(2)}_{f(x,g)}\Tilde{\varpi}(\Tilde{g})\psi(x),\qquad \forall x\in\Zold{Z:BehaviorOfpsiW_t:psi(W_tx)=U_(f_j(x,s))g(s)psi(x)}(g),\quad\Tilde{g}\in\pi_{\Wbf}^{-1}(g);
        \end{gathered}
    \end{equation*}
    where $f:\{(x,g):x\in\Zold{Z:BehaviorOfpsiW_t:psi(W_tx)=U_(f_j(x,s))g(s)psi(x)}(g),g\in\Wbf\}\to\R$ is a measurable function satisfying that for $r\in\R$, $g\in\Wbf$, $\gamma\in\Gamma_{1,\nc}$,
       \begin{equation*}
        \begin{gathered}
        f(U^{(1)}_rx,g)-f(x,g)=\omega(gx,r)-\omega(x,r),\quad \forall x\in\Zold{Z:BehaviorOfpsiW_t:psi(W_tx)=U_(f_j(x,s))g(s)psi(x)}(g);\\
        f(x,\gamma)=0,\quad \forall x\in\Zold{Z:BehaviorOfpsiW_t:psi(W_tx)=U_(f_j(x,s))g(s)psi(x)}(\gamma);\qquad f(x,g\gamma)=f(x,g),\quad \forall x\in\Zold{Z:BehaviorOfpsiW_t:psi(W_tx)=U_(f_j(x,s))g(s)psi(x)}(g).
    \end{gathered}
    \end{equation*}
    
\end{Proposition}
Let us begin the proof of \textit{Proposition \ref{Prop:BehaviorOfpsiW_t:psi(W_tx)=...psi(x)}}.

First, we apply \textit{Proposition \ref{Proposition:BehaviorOfpsi:psi(Kak)=Kak}} to obtain the following local relation:
\begin{Lemma}\label{Lemma:BehaviorOfpsiW_t:psi(W_sx)=...ForSmalls}

    For $\delta\in(0,\deold{delta:BehaviorOfpsi:BoundOndeltaInTwoSideMatching}]$ and $g\in\overline{B}^{\Wbf}_{\eold{epsilon:BehaviorOfpsi:distanceBoundOnG_(d1)/Gamma_1}(\delta)}$, define
    \[
    \Znew\label{Z:BehaviorOfpsiW_t:pointsSatisfyingpsi(W_tx)=...psi(x)}(g)=\Zold{Z:BehaviorOfpsiUnderA_t:GenericPointsOfGoodSetInLemmapsi(Kak)=Kak}\cap g^{-1}\Zold{Z:BehaviorOfpsiUnderA_t:GenericPointsOfGoodSetInLemmapsi(Kak)=Kak}.
    \]
    Then for $x\in\Zold{Z:BehaviorOfpsiW_t:pointsSatisfyingpsi(W_tx)=...psi(x)}(g)$,
    \[
    \psi(gx)=U^{(2)}_{c_{\bu}(x,g)}A^{(2)}_{c_{\ba}(x,g)}\varpi(x,g)\psi(x)
    \]
    where
    \[
    c_{\bu}(x,g)\in\R,\quad\absolute{c_{\ba}(x,g)}\leq\delta,\quad \varpi(x,g)\in \overline{B}^{\Gbf_2}_\delta\cap\log^{-1}\left(\textstyle{\bigoplus_{k=1}^{n^{(2)}}\R\bx^{(2)}_{0,k}}\right),
    \]
    and they are unique in the sense as in \textit{Lemma~\ref{Lemma:BasicLemma:BasicLemma}-\ref{item:BasicLemma:Statement:C}}. Moreover, the set
    \[
    \left\{(x,g):x\in\Zold{Z:BehaviorOfpsiW_t:pointsSatisfyingpsi(W_tx)=...psi(x)}(g),g\in\overline{B}^{\Wbf}_{\eold{epsilon:BehaviorOfpsi:distanceBoundOnG_(d1)/Gamma_1}(\deold{delta:BehaviorOfpsi:BoundOndeltaInTwoSideMatching})}\right\}
    \]
    is of full measure in $(\Gbf_1/\Gamma_1)\times\overline{B}^{\Wbf}_{\eold{epsilon:BehaviorOfpsi:distanceBoundOnG_(d1)/Gamma_1}(\deold{delta:BehaviorOfpsi:BoundOndeltaInTwoSideMatching})}$, and on it, $c_{\bu}(\cdot,\cdot)$, $c_{\ba}(\cdot,\cdot)$, $\varpi(\cdot,\cdot)$ are measurable.
    
\end{Lemma}
\begin{proof}
Fix $\delta\in(0,\deold{delta:BehaviorOfpsi:BoundOndeltaInTwoSideMatching}]$ and $g\in\overline{B}^{\Wbf}_{\eold{epsilon:BehaviorOfpsi:distanceBoundOnG_(d1)/Gamma_1}(\delta)}$. For $x\in\Zold{Z:BehaviorOfpsiW_t:pointsSatisfyingpsi(W_tx)=...psi(x)}(g)$, \ref{item:BehaviorOfpsiA_t:1} implies that $x$ is a generic point of both $\mathbf{1}_{\Kold{K:BehaviorOfpsi:GoodSetInLemmapsi(Kak)=Kak}}$ and $\mathbf{1}_{g^{-1}\Kold{K:BehaviorOfpsi:GoodSetInLemmapsi(Kak)=Kak}}$ w.r.t. the ergodic flow $U^{(1)}_t$. Hence, there exists $t(x)>0$ such that
\[
U^{(1)}_{t(x)}x\in\Kold{K:BehaviorOfpsi:GoodSetInLemmapsi(Kak)=Kak},\qquad gU^{(1)}_{t(x)}x\in\Kold{K:BehaviorOfpsi:GoodSetInLemmapsi(Kak)=Kak}.
\]
And since $gU^{(1)}_{t(x)}x\in\Kak(\infty,\eold{epsilon:BehaviorOfpsi:distanceBoundOnG_(d1)/Gamma_1}(\delta),U^{(1)}_{t(x)}x)$, it follows from \textit{Proposition \ref{Proposition:BehaviorOfpsi:psi(Kak)=Kak}} that for some $c_\bu\in[-\delta,\delta]$, $c_\ba\in[-\delta,\delta]$, $g'\in\overline{B}^{\Gbf_2}_\delta\cap\log^{-1}\left(\oplus_{k=1}^{n^{(2)}}\R\bx^{(2)}_{0,k}\right)$, unique in the sense of \textit{Lemma \ref{Lemma:BasicLemma:BasicLemma}-\ref{item:BasicLemma:Statement:C}},
\[
\psi(gU^{(1)}_{t(x)}x)=U^{(2)}_{c_\bu}A^{(2)}_{c_\ba}g'\psi(U^{(1)}_{t(x)}x).
\]
This and \eqref{eq:BehaviorOfpsi:ConjugationRelation} give us the first part of \textit{Lemma \ref{Lemma:BehaviorOfpsiW_t:psi(W_sx)=...ForSmalls}}.

\medskip

The second part of \textit{Lemma \ref{Lemma:BehaviorOfpsiW_t:psi(W_sx)=...ForSmalls}} follows once we show 
\[
\left\{(x,g):x\in\Zold{Z:BehaviorOfpsiW_t:pointsSatisfyingpsi(W_tx)=...psi(x)}(g),g\in\overline{B}^{\Wbf}_{\eold{epsilon:BehaviorOfpsi:distanceBoundOnG_(d1)/Gamma_1}(\deold{delta:BehaviorOfpsi:BoundOndeltaInTwoSideMatching})}\right\}
\]
is measurable in $(\Gbf_1/\Gamma_1)\times\Wbf$. And this holds because the above set is 
\[
\Phi^{-1}(\Zold{Z:BehaviorOfpsiUnderA_t:GenericPointsOfGoodSetInLemmapsi(Kak)=Kak})\cap\left(\Zold{Z:BehaviorOfpsiUnderA_t:GenericPointsOfGoodSetInLemmapsi(Kak)=Kak}\times\overline{B}^{\Wbf}_{\eold{epsilon:BehaviorOfpsi:distanceBoundOnG_(d1)/Gamma_1}(\deold{delta:BehaviorOfpsi:BoundOndeltaInTwoSideMatching})}\right),
\]
where $\Phi:(\Gbf_1/\Gamma_1)\times\Wbf\to\Gbf_1/\Gamma_1$ is given by $(x,g)\mapsto gx$.

\medskip

Thus, we finish the proof of \textit{Lemma \ref{Lemma:BehaviorOfpsiW_t:psi(W_sx)=...ForSmalls}}.
\end{proof}

In the above lemma, we see that the second part also holds when $\eold{epsilon:BehaviorOfpsi:distanceBoundOnG_(d1)/Gamma_1}(\deold{delta:BehaviorOfpsi:BoundOndeltaInTwoSideMatching})$ is replaced by some $\epsilon\in[0,\eold{epsilon:BehaviorOfpsi:distanceBoundOnG_(d1)/Gamma_1}(\deold{delta:BehaviorOfpsi:BoundOndeltaInTwoSideMatching})]$, and for $g\in\overline{B}^{\Wbf}_{\eold{epsilon:BehaviorOfpsi:distanceBoundOnG_(d1)/Gamma_1}(\deold{delta:BehaviorOfpsi:BoundOndeltaInTwoSideMatching})}$,
\[
c_{\bu}(\cdot,g),\qquad c_{\ba}(\cdot,g),\qquad \varpi(\cdot,g)
\]
are measurable on $\Gbf_1/\Gamma_1$. 

\medskip

Next, we show that, in fact, $c_{\ba}(\cdot,g)$ and $\varpi(\cdot,g)$ depend only on $g$:
\begin{Lemma}\label{Lemma:BehaviorOfpsiW_t:c_a(.,s),g(.,s)areU_r-invariant}

    Let $c_{0,k}:\left\{(x,s):x\in\Zold{Z:BehaviorOfpsiW_t:pointsSatisfyingpsi(W_tx)=...psi(x)}(g),g\in\overline{B}^{\Wbf}_{\eold{epsilon:BehaviorOfpsi:distanceBoundOnG_(d1)/Gamma_1}(\deold{delta:BehaviorOfpsi:BoundOndeltaInTwoSideMatching})}\right\}\to\R$ be a measurable function for all $1\leq k\leq n^{(2)}$ such that
    \[
    \log(\varpi(x,g))=\textstyle{\sum_{k=1}^{n^{(2)}}c_{0,k}(x,g)\bx^{(2)}_{0,k}}.
    \]
    For $g\in\overline{B}^{\Wbf}_{\eold{epsilon:BehaviorOfpsi:distanceBoundOnG_(d1)/Gamma_1}(\deold{delta:BehaviorOfpsi:BoundOndeltaInTwoSideMatching})}$, define
    \begin{equation}\label{eq:BehaviorOfpsiW_t:DefOfc_a(s),g(s)}
        \begin{aligned}
            c_{\ba}(g)&=\int_{\Gbf_1/\Gamma_1}c_{\ba}(x,g)d\mu_1(x);\\
            c_{0,k}(g)&=\int_{\Gbf_1/\Gamma_1}c_{0,k}(x,g)d\mu_1(x),\qquad \forall 1\leq k\leq n^{(2)};\\
            \varpi(g)&=\exp\left(\textstyle{\sum_{k=1 }^{n^{(2)}}c_{0,k}(g)\bx^{(2)}_{0,k}}\right).
        \end{aligned}
    \end{equation}
    Then $c_{\ba}(\cdot),c_{0,k}(\cdot)$ for $1\leq k\leq n^{(2)}$, and $\varpi(\cdot)$ are measurable maps defined on $\overline{B}^{\Wbf}_{\eold{epsilon:BehaviorOfpsi:distanceBoundOnG_(d1)/Gamma_1}(\deold{delta:BehaviorOfpsi:BoundOndeltaInTwoSideMatching})}$ such that for all $\delta\in(0,\deold{delta:BehaviorOfpsi:BoundOndeltaInTwoSideMatching}]$ and $g\in\overline{B}^{\Wbf}_{\eold{epsilon:BehaviorOfpsi:distanceBoundOnG_(d1)/Gamma_1}(\delta)}$, 
    \[
    c_{\ba}(g)\in[-\delta,\delta],\qquad \varpi(g)\in \overline{B}^{\Gbf_2}_\delta\cap\log^{-1}\left(\textstyle{\bigoplus_{k=1}^{n^{(2)}}\R\bx^{(2)}_{0,k}}\right).
    \]
    Moreover, for $g\in\overline{B}^{\Wbf}_{\eold{epsilon:BehaviorOfpsi:distanceBoundOnG_(d1)/Gamma_1}(\deold{delta:BehaviorOfpsi:BoundOndeltaInTwoSideMatching})}$,
    \[
    \Znew\label{Z:BehaviorOfpsiW_t:SetOnWhichc_a(.,p),g(.,p)areConstants}(g)=\left\{x\in\Zold{Z:BehaviorOfpsiW_t:pointsSatisfyingpsi(W_tx)=...psi(x)}(g):c_{\ba}(x,g)=c_{\ba}(g),c_{0,k}(x,g)=c_{0,k}(g),1\leq k\leq n^{(2)}\right\}
    \]
    is a full measure $U^{(1)}_t$-invariant set in $\Zold{Z:BehaviorOfpsiW_t:pointsSatisfyingpsi(W_tx)=...psi(x)}(g)$; the set
    \[
    \left\{(x,g):x\in\Zold{Z:BehaviorOfpsiW_t:SetOnWhichc_a(.,p),g(.,p)areConstants}(g),g\in\overline{B}^{\Wbf}_{\eold{epsilon:BehaviorOfpsi:distanceBoundOnG_(d1)/Gamma_1}(\deold{delta:BehaviorOfpsi:BoundOndeltaInTwoSideMatching})}\right\}
    \]
    is of full measure in $(\Gbf_1/\Gamma_1)\times\overline{B}^{\Wbf}_{\eold{epsilon:BehaviorOfpsi:distanceBoundOnG_(d1)/Gamma_1}(\deold{delta:BehaviorOfpsi:BoundOndeltaInTwoSideMatching})}$; and also for every $ x\in\Zold{Z:BehaviorOfpsiW_t:SetOnWhichc_a(.,p),g(.,p)areConstants}(g)$ and $g\in\overline{B}^{\Wbf}_{\eold{epsilon:BehaviorOfpsi:distanceBoundOnG_(d1)/Gamma_1}(\deold{delta:BehaviorOfpsi:BoundOndeltaInTwoSideMatching})}$, we have
    \[
    c_{\ba}(x,g)=c_{\ba}(g),\qquad \varpi(x,g)=\varpi(g).
    \]
    
\end{Lemma}
\begin{proof}
By \textit{Lemma~\ref{Lemma:BehaviorOfpsiW_t:psi(W_sx)=...ForSmalls}}, the comments following its proof, and \textit{Fubini's theorem}, we obtain the first part of \textit{Lemma~\ref{Lemma:BehaviorOfpsiW_t:c_a(.,s),g(.,s)areU_r-invariant}}.

\medskip

For the second part, we recall that $\Zold{Z:BehaviorOfpsiUnderA_t:GenericPointsOfGoodSetInLemmapsi(Kak)=Kak}$ is $U^{(1)}_t$-invariant (cf. \ref{item:BehaviorOfpsiA_t:1}), hence so is $\Zold{Z:BehaviorOfpsiW_t:pointsSatisfyingpsi(W_tx)=...psi(x)}(g)=\Zold{Z:BehaviorOfpsiUnderA_t:GenericPointsOfGoodSetInLemmapsi(Kak)=Kak}\cap g^{-1}\Zold{Z:BehaviorOfpsiUnderA_t:GenericPointsOfGoodSetInLemmapsi(Kak)=Kak}$ for all $g\in\overline{B}^{\Wbf}_{\eold{epsilon:BehaviorOfpsi:distanceBoundOnG_(d1)/Gamma_1}(\deold{delta:BehaviorOfpsi:BoundOndeltaInTwoSideMatching})}$. This, \textit{Lemma~\ref{Lemma:BehaviorOfpsiW_t:psi(W_sx)=...ForSmalls}}, and \eqref{eq:BehaviorOfpsi:ConjugationRelation} imply that for $g\in\overline{B}^{\Wbf}_{\eold{epsilon:BehaviorOfpsi:distanceBoundOnG_(d1)/Gamma_1}(\deold{delta:BehaviorOfpsi:BoundOndeltaInTwoSideMatching})}$, $x\in\Zold{Z:BehaviorOfpsiW_t:pointsSatisfyingpsi(W_tx)=...psi(x)}(g)$, and $t\in\R$,
\begin{align*}
    &U^{(2)}_{c'}A^{(2)}_{c_{\ba}(U^{(1)}_tx,g)}\varpi(U^{(1)}_tx,g)\psi(x)\\
    =\;&\psi(gU^{(1)}_tx)=\psi(U^{(1)}_tgx)=U^{(2)}_{c''}A^{(2)}_{c_{\ba}(x,s)}\varpi(x,g)\psi(x),
\end{align*}
where $c',c''$ are some real numbers. By the uniqueness part in \textit{Lemma~\ref{Lemma:BehaviorOfpsiW_t:psi(W_sx)=...ForSmalls}}, we deduce that $c_{\ba}(\cdot,g)$, $\varpi(\cdot,g)$, and consequently, $c_{0,k}(\cdot,g)$, $1\leq k\leq n^{(2)}$, are $U^{(1)}_t$-invariant. Since the flow $U^{(1)}_t$ on $(\Gbf_1/\Gamma_1,\mu_1)$ is ergodic and these functions are measurable, this gives us the second part of \textit{Lemma~\ref{Lemma:BehaviorOfpsiW_t:c_a(.,s),g(.,s)areU_r-invariant}}.

\medskip

Thus, we finish the proof of \textit{Lemma~\ref{Lemma:BehaviorOfpsiW_t:c_a(.,s),g(.,s)areU_r-invariant}}.
\end{proof}
Similar as for \textit{Lemma~\ref{Lemma:BehaviorOfpsiW_t:psi(W_sx)=...ForSmalls}}, the second part of \textit{Lemma~\ref{Lemma:BehaviorOfpsiW_t:c_a(.,s),g(.,s)areU_r-invariant}} also holds when $\eold{epsilon:BehaviorOfpsi:distanceBoundOnG_(d1)/Gamma_1}(\deold{delta:BehaviorOfpsi:BoundOndeltaInTwoSideMatching})$ is replaced by some $\epsilon\in(0,\eold{epsilon:BehaviorOfpsi:distanceBoundOnG_(d1)/Gamma_1}(\deold{delta:BehaviorOfpsi:BoundOndeltaInTwoSideMatching})]$.

\medskip

We continue with showing that in \textit{Lemma~\ref{Lemma:BehaviorOfpsiW_t:c_a(.,s),g(.,s)areU_r-invariant}}, $\varpi(g)$ commutes with $\Hbf_2$:
\begin{Lemma}\label{Lemma:BehaviorOfpsiW_t:g(s)CommuteH_2}

    For $g\in\overline{B}^{\Wbf}_{\eold{epsilon:BehaviorOfpsi:distanceBoundOnG_(d1)/Gamma_1}(\deold{delta:BehaviorOfpsi:BoundOndeltaInTwoSideMatching})}$, let $c_{\bu}(\cdot,g)$ be given in \textit{Lemma~\ref{Lemma:BehaviorOfpsiW_t:psi(W_sx)=...ForSmalls}} and $c_{\ba}(g)$, $\varpi(g)$ given in \textit{Lemma~\ref{Lemma:BehaviorOfpsiW_t:c_a(.,s),g(.,s)areU_r-invariant}}. Then
    \[
    \textstyle{\varpi(g)\in\exp\left(\bigoplus_{k\in J_2}\R\bx^{(2)}_{0,k}\right).}
    \]
\end{Lemma}
\begin{proof}
With $\varpi(g)=\exp\left(\sum_{k=1}^{n^{(2)}}c_{0,k}(g)\bx^{(2)}_{0,k}\right)$ (cf.~\eqref{eq:BehaviorOfpsiW_t:DefOfc_a(s),g(s)}), since $g$ commutes with $\Hbf_1$, \textit{Lemma~\ref{Lemma:BehaviorOfpsiW_t:c_a(.,s),g(.,s)areU_r-invariant}} and \ref{item:BehaviorOfpsiW_t:1} imply that $\mu_1$-a.e. $x\in\Gbf_1/\Gamma_1$ satisfy
\begin{align*}
    &U^{(2)}_{c_{\bu}(x,g)}A^{(2)}_{c_{\ba}(g)}\exp\left(\textstyle{\sum_{k=1}^{n^{(2)}}c_{0,k}(g)\bx^{(2)}_{0,k}}\right)\psi(x)=\psi(gx)\\
    =\;&\psi\left(A^{(1)}_{-p}gA^{(1)}_px\right)=U^{(2)}_{\tilde{c}}A^{(2)}_{c_{\ba}(g)}\exp\left(\textstyle{\sum_{k=1}^{n^{(2)}}e^{-q^{(2)}_kp}c_{0,k}(g)\bx^{(2)}_{0,k}}\right)\psi(x),
\end{align*}
where $\tilde{c}$ is a real constant. By the uniqueness part in \textit{Lemma~\ref{Lemma:BehaviorOfpsiW_t:psi(W_sx)=...ForSmalls}}, 
\[
\exp\left(\textstyle{\sum_{k=1}^{n^{(2)}}c_{0,k}(g)\bx^{(2)}_{0,k}}\right)=\exp\left(\textstyle{\sum_{k=1}^{n^{(2)}}e^{-q^{(2)}_kp}c_{0,k}(g)\bx^{(2)}_{0,k}}\right),
\]
which implies that $c_{0,k}(g)=0$ whenever $q^{(2)}_k\neq0$, i.e. $k\notin J_2$. Thus, 
\[
\varpi(g)=\exp\left(\textstyle{\sum_{k\in J_2}c_{0,k}(g)\bx^{(2)}_{0,k}}\right),
\]
and we finish the proof of \textit{Lemma~\ref{Lemma:BehaviorOfpsiW_t:g(s)CommuteH_2}}.
\end{proof}
With the above lemmas, we show that $c_\ba$ and $\varpi$ are local homomorphisms:
\begin{Lemma}\label{Lemma:BehaviorOfpsiW_t:A_(c_a(s))g(s)IsA1-parasubgroup}

    Let $c_{\ba}(\cdot)$ and $\varpi(\cdot)$ be given in \textit{Lemma~\ref{Lemma:BehaviorOfpsiW_t:c_a(.,s),g(.,s)areU_r-invariant}}. Then for every $g_1,g_2\in\overline{B}^{\Wbf}_{0.5\eold{epsilon:BehaviorOfpsi:distanceBoundOnG_(d1)/Gamma_1}(\deold{delta:BehaviorOfpsi:BoundOndeltaInTwoSideMatching}/10)}$,
    \[
    c_{\ba}(g_1g_2)=c_{\ba}(g_1)+c_{\ba}(g_2),\qquad \varpi(g_1g_2)=\varpi(g_1)\varpi(g_2).
    \]
    
\end{Lemma}
\begin{proof}
For $x\in g_2^{-1}\Zold{Z:BehaviorOfpsiW_t:SetOnWhichc_a(.,p),g(.,p)areConstants}(g_1)\cap\Zold{Z:BehaviorOfpsiW_t:SetOnWhichc_a(.,p),g(.,p)areConstants}(g_2)\cap\Zold{Z:BehaviorOfpsiW_t:SetOnWhichc_a(.,p),g(.,p)areConstants}(g_1g_2)$, \textit{Lemmas~\ref{Lemma:BehaviorOfpsiW_t:psi(W_sx)=...ForSmalls},~\ref{Lemma:BehaviorOfpsiW_t:c_a(.,s),g(.,s)areU_r-invariant},~\ref{Lemma:BehaviorOfpsiW_t:g(s)CommuteH_2}}, and the \textit{Baker-Campbell-Hausdorff formula} imply that
\begin{align*}
    &U^{(2)}_{c_{\bu}(x,g_1g_2)}A^{(2)}_{c_{\ba}(g_1g_2)}\varpi(g_1g_2)\psi(x)=\psi((g_1g_2)x)\\
    =\;&\psi(g_1(g_2x))=U^{(2)}_{c_{\bu}(g_2x,g_1)}A^{(2)}_{c_{\ba}(g_1)}\varpi(g_1)\psi(g_2x)\\
    =\;&U^{(2)}_{c_{\bu}(g_2x,g_1)}A^{(2)}_{c_{\ba}(g_1)}\varpi(g_1)U^{(2)}_{c_{\bu}(x,g_2)}A^{(2)}_{c_{\ba}(g_2)}\varpi(g_2)\psi(x)\\
    =\;&U^{(2)}_{c_{\bu}(g_2x,g_1)+c_{\bu}(x,g_2)e^{2c_{\ba}(g_1)}}A^{(2)}_{c_{\ba}(g_1)+c_{\ba}(g_2)}\left(\varpi(g_1)\varpi(g_2)\right)\psi(x),
\end{align*}
where $\varpi(g_1)\varpi(g_2)\in \overline{B}^{\Gbf_2}_{\deold{delta:BehaviorOfpsi:BoundOndeltaInTwoSideMatching}/5}\cap\exp\left(\bigoplus_{k\in J_{2}}\R\bx^{(2)}_{0,k}\right)$. Hence, it follows from the uniqueness part in \textit{Lemma~\ref{Lemma:BehaviorOfpsiW_t:psi(W_sx)=...ForSmalls}} that for $g_1,g_2\in\overline{B}^{\Wbf}_{0.5\eold{epsilon:BehaviorOfpsi:distanceBoundOnG_(d1)/Gamma_1}(\deold{delta:BehaviorOfpsi:BoundOndeltaInTwoSideMatching}/10)}$, we have
\[
c_{\ba}(g_1g_2)=c_{\ba}(g_1)+c_{\ba}(g_2),\qquad\varpi(g_1g_2)=\varpi(g_1)\varpi(g_2).
\]
This finishes the proof of \textit{Lemma~\ref{Lemma:BehaviorOfpsiW_t:A_(c_a(s))g(s)IsA1-parasubgroup}}.
\end{proof}
Now, let $\denew\label{delta:BehaviorOfpsiW_t:B^(tilde(W))_delta==B^W_delta}\in(0,\frac{1}{2}\eold{epsilon:BehaviorOfpsi:distanceBoundOnG_(d1)/Gamma_1}(\frac{\deold{delta:BehaviorOfpsi:BoundOndeltaInTwoSideMatching}}{10}))$ be small enough such that $\pi_\Wbf: B^{\widetilde{\Wbf}}_{\deold{delta:BehaviorOfpsiW_t:B^(tilde(W))_delta==B^W_delta}}\to B^\Wbf_{\deold{delta:BehaviorOfpsiW_t:B^(tilde(W))_delta==B^W_delta}}$ is a bijective isometry. We define for $\Tilde{g}\in  B^{\widetilde{\Wbf}}_{\deold{delta:BehaviorOfpsiW_t:B^(tilde(W))_delta==B^W_delta}}$,
\[
\Znew\label{Z:BehaviorOfpsiW_t:psi(pi_W(tilde(g)x))=...psi(x)}(\Tilde{g})=\Zold{Z:BehaviorOfpsiW_t:SetOnWhichc_a(.,p),g(.,p)areConstants}(\pi_\Wbf(\Tilde{g})),\quad\Tilde{c}_\ba(\Tilde{g})=c_\ba(\pi_\Wbf(\Tilde{g})),\quad\Tilde{\varpi}(\Tilde{g})=\varpi(\pi_\Wbf(\Tilde{g})),
\]
\[
\Tilde{c}_\bu(x,\Tilde{g})=c_\bu(x,\pi_\Wbf(\Tilde{g})),\qquad\forall x\in\Zold{Z:BehaviorOfpsiW_t:psi(pi_W(tilde(g)x))=...psi(x)}(\Tilde{g}).
\]
By \textit{Lemmas~\ref{Lemma:BehaviorOfpsiW_t:psi(W_sx)=...ForSmalls},~\ref{Lemma:BehaviorOfpsiW_t:c_a(.,s),g(.,s)areU_r-invariant}}, and the comments following them, the above construction implies that for all $\Tilde{g}\in  B^{\widetilde{\Wbf}}_{\deold{delta:BehaviorOfpsiW_t:B^(tilde(W))_delta==B^W_delta}}$, $\Zold{Z:BehaviorOfpsiW_t:psi(pi_W(tilde(g)x))=...psi(x)}(\Tilde{g})\subseteq\Zold{Z:BehaviorOfpsiUnderA_t:GenericPointsOfGoodSetInLemmapsi(Kak)=Kak}\cap\pi_\Wbf(\Tilde{g})^{-1}\Zold{Z:BehaviorOfpsiUnderA_t:GenericPointsOfGoodSetInLemmapsi(Kak)=Kak}$ is a full measure $U^{(1)}_t$-invariant set and
\[
\psi(\pi_{\Wbf}(\Tilde{g})x)=U^{(2)}_{\Tilde{c}_\bu(x,\Tilde{g})}A^{(2)}_{\Tilde{c}_\ba(\Tilde{g})}\Tilde{\varpi}(\Tilde{g})\psi(x),\qquad\forall x\in\Zold{Z:BehaviorOfpsiW_t:psi(pi_W(tilde(g)x))=...psi(x)}(\Tilde{g}).
\]
Moreover, the set
\[
\left\{(x,\Tilde{g}):x\in\Zold{Z:BehaviorOfpsiW_t:psi(pi_W(tilde(g)x))=...psi(x)}(\Tilde{g}),\Tilde{g}\in  B^{\widetilde{\Wbf}}_{\deold{delta:BehaviorOfpsiW_t:B^(tilde(W))_delta==B^W_delta}}\right\}
\]
is of full measure in $(\Gbf_1/\Gamma_1)\times B^{\widetilde{\Wbf}}_{\deold{delta:BehaviorOfpsiW_t:B^(tilde(W))_delta==B^W_delta}}$, $\Tilde{c}_\bu$ is measurable on this set, and $\Tilde{c}_\bu(\cdot,\Tilde{g})$ is measurable on $\Gbf_1/\Gamma_1$ for all $\Tilde{g}\in B^{\widetilde{\Wbf}}_{\deold{delta:BehaviorOfpsiW_t:B^(tilde(W))_delta==B^W_delta}}$.

Furthermore, since $\widetilde{\Wbf}$ is a universal cover of $\Wbf$ (cf. \ref{item:BehaviorOfpsiW_t:9}), by the above construction, \textit{Lemma~\ref{Lemma:BehaviorOfpsiW_t:A_(c_a(s))g(s)IsA1-parasubgroup}}, and the continuity of $c_\ba$, $\varpi$ at $e_{\Gbf_1}$ (cf. \textit{Lemmas~\ref{Lemma:BehaviorOfpsiW_t:psi(W_sx)=...ForSmalls},~\ref{Lemma:BehaviorOfpsiW_t:c_a(.,s),g(.,s)areU_r-invariant}}), we can extend
\[
\Tilde{c}_\ba: B^{\widetilde{\Wbf}}_{\deold{delta:BehaviorOfpsiW_t:B^(tilde(W))_delta==B^W_delta}}\to\R,\qquad\Tilde{\varpi}: B^{\widetilde{\Wbf}}_{\deold{delta:BehaviorOfpsiW_t:B^(tilde(W))_delta==B^W_delta}}\to\Zbf_2
\]
to continuous homomorphisms
\begin{equation}\label{eq:tildeVarpiDef}
\Tilde{c}_\ba:\widetilde{\Wbf}\to\R,\qquad\Tilde{\varpi}:\widetilde{\Wbf}\to\Zbf_2. 
\end{equation}

\medskip

With these, using induction, we obtain the following:
\begin{Lemma}\label{Lemma:BehaviorOfpsiW_t:GoodSet(tilde(g))Andc_u(.,tilde(g))}

    There exist $\Zold{Z:BehaviorOfpsiW_t:psi(pi_W(tilde(g)x))=...psi(x)}(\Tilde{g})\subseteq\Zold{Z:BehaviorOfpsiUnderA_t:GenericPointsOfGoodSetInLemmapsi(Kak)=Kak}\cap\pi_\Wbf(\Tilde{g})^{-1}\Zold{Z:BehaviorOfpsiUnderA_t:GenericPointsOfGoodSetInLemmapsi(Kak)=Kak}$ full measure $U^{(1)}_t$-invariant sets for all $\Tilde{g}\in\widetilde{\Wbf}$ such that the set
    \[
    \left\{(x,\Tilde{g}):x\in\Zold{Z:BehaviorOfpsiW_t:psi(pi_W(tilde(g)x))=...psi(x)}(\Tilde{g}),\Tilde{g}\in\widetilde{\Wbf}\right\}
    \]
    is of full measure in $(\Gbf_1/\Gamma_1)\times\widetilde{\Wbf}$, and for $\Tilde{g}\in\widetilde{\Wbf}$,
    \[
    \psi(\pi_\Wbf(\Tilde{g})x)=U^{(2)}_{\Tilde{c}_\bu(x,\Tilde{g})}A^{(2)}_{\Tilde{c}_\ba(\Tilde{g})}\Tilde{\varpi}(\Tilde{g})\psi(x),\qquad\forall x\in\Zold{Z:BehaviorOfpsiW_t:psi(pi_W(tilde(g)x))=...psi(x)}(\Tilde{g}).
    \]
    Here, $\Tilde{c}_\bu:\left\{(x,\Tilde{g}):x\in\Zold{Z:BehaviorOfpsiW_t:psi(pi_W(tilde(g)x))=...psi(x)}(\Tilde{g}),\Tilde{g}\in\widetilde{\Wbf}\right\}\to\R$ is measurable, and for $\Tilde{g}\in\widetilde{\Wbf}$, $\Tilde{c}_\bu(\cdot,\Tilde{g})$ is measurable on $\Gbf_1/\Gamma_1$.
    
\end{Lemma}
\begin{proof}
Based on the fact that $\widetilde{\Wbf}=\bigcup_{n\in\N}\left( B^{\widetilde{\Wbf}}_{\deold{delta:BehaviorOfpsiW_t:B^(tilde(W))_delta==B^W_delta}}\right)^n$, we prove this lemma using induction. In particular, we show that this lemma holds when $\widetilde{\Wbf}$ is replaced by $\left( B^{\widetilde{\Wbf}}_{\deold{delta:BehaviorOfpsiW_t:B^(tilde(W))_delta==B^W_delta}}\right)^n$, for all $n\in\N$.

The base case $n=1$ follows from the above arguments above. Now, assuming \textit{Lemma~\ref{Lemma:BehaviorOfpsiW_t:GoodSet(tilde(g))Andc_u(.,tilde(g))}} holds when $\widetilde{\Wbf}$ is replaced by $\left( B^{\widetilde{\Wbf}}_{\deold{delta:BehaviorOfpsiW_t:B^(tilde(W))_delta==B^W_delta}}\right)^n$, $n\in\N$, we show that it also holds when $\widetilde{\Wbf}$ is replaced by $\left( B^{\widetilde{\Wbf}}_{\deold{delta:BehaviorOfpsiW_t:B^(tilde(W))_delta==B^W_delta}}\right)^{n+1}$.

Fix a countable dense subset $K$ in $B^{\widetilde{\Wbf}}_{\deold{delta:BehaviorOfpsiW_t:B^(tilde(W))_delta==B^W_delta}}$. For $\Tilde{g}\in\left( B^{\widetilde{\Wbf}}_{\deold{delta:BehaviorOfpsiW_t:B^(tilde(W))_delta==B^W_delta}}\right)^{n+1}\setminus\left( B^{\widetilde{\Wbf}}_{\deold{delta:BehaviorOfpsiW_t:B^(tilde(W))_delta==B^W_delta}}\right)^n$, let
\[
K(\Tilde{g})=\left\{\Tilde{h}\in K:\Tilde{g}\Tilde{h}^{-1}\in\left(B^{\widetilde{\Wbf}}_{\deold{delta:BehaviorOfpsiW_t:B^(tilde(W))_delta==B^W_delta}}\right)^n\right\}
\]
and note that $K(\Tilde{g})\neq\emptyset$ since $B^{\widetilde{\Wbf}}_{\deold{delta:BehaviorOfpsiW_t:B^(tilde(W))_delta==B^W_delta}}\cap\left(B^{\widetilde{\Wbf}}_{\deold{delta:BehaviorOfpsiW_t:B^(tilde(W))_delta==B^W_delta}}\right)^n\Tilde{g}\neq\emptyset$ is open, we define
\[
\textstyle{\Zold{Z:BehaviorOfpsiW_t:psi(pi_W(tilde(g)x))=...psi(x)}(\Tilde{g})=\bigcap_{\Tilde{h}\in K(\Tilde{g})}\left(\pi_{\Wbf}(\Tilde{h})^{-1}\Zold{Z:BehaviorOfpsiW_t:psi(pi_W(tilde(g)x))=...psi(x)}(\Tilde{g}\Tilde{h}^{-1})\cap\Zold{Z:BehaviorOfpsiW_t:psi(pi_W(tilde(g)x))=...psi(x)}(\Tilde{h})\right),}
\]
it follows from the induction hypothesis that $\Zold{Z:BehaviorOfpsiW_t:psi(pi_W(tilde(g)x))=...psi(x)}(\Tilde{g})\subseteq\Zold{Z:BehaviorOfpsiUnderA_t:GenericPointsOfGoodSetInLemmapsi(Kak)=Kak}\cap\pi_\Wbf(\Tilde{g})^{-1}\Zold{Z:BehaviorOfpsiUnderA_t:GenericPointsOfGoodSetInLemmapsi(Kak)=Kak}$ is a full measure $U^{(1)}_t$-invariant set for all $\Tilde{g}\in\left( B^{\widetilde{\Wbf}}_{\deold{delta:BehaviorOfpsiW_t:B^(tilde(W))_delta==B^W_delta}}\right)^{n+1}$, and
\[
\left\{(x,\Tilde{g}):x\in\Zold{Z:BehaviorOfpsiW_t:psi(pi_W(tilde(g)x))=...psi(x)}(\Tilde{g}),\Tilde{g}\in\left( B^{\widetilde{\Wbf}}_{\deold{delta:BehaviorOfpsiW_t:B^(tilde(W))_delta==B^W_delta}}\right)^{n+1}\right\}
\]
is of full measure in $(\Gbf_1/\Gamma_1)\times\left( B^{\widetilde{\Wbf}}_{\deold{delta:BehaviorOfpsiW_t:B^(tilde(W))_delta==B^W_delta}}\right)^{n+1}$.

Also by the induction hypothesis, for all $\Tilde{g}\in\left( B^{\widetilde{\Wbf}}_{\deold{delta:BehaviorOfpsiW_t:B^(tilde(W))_delta==B^W_delta}}\right)^{n+1}\setminus \left( B^{\widetilde{\Wbf}}_{\deold{delta:BehaviorOfpsiW_t:B^(tilde(W))_delta==B^W_delta}}\right)^n$ and $x\in \Zold{Z:BehaviorOfpsiW_t:psi(pi_W(tilde(g)x))=...psi(x)}(\Tilde{g})$, we have
\begin{align*}
    \psi(\pi_\Wbf(\Tilde{g})x)&=\psi(\pi_{\Wbf}(\Tilde{g}_1)\pi_{\Wbf}(\Tilde{g}_2)x)=U^{(2)}_{\Tilde{c}_\bu(\pi_{\Wbf}(\Tilde{g}_2)x,\Tilde{g}_1)}A^{(2)}_{\Tilde{c}_\ba(\Tilde{g}_1)}\Tilde{\varpi}(\Tilde{g}_1)\psi(\pi_{\Wbf}(\Tilde{g}_2)x)\\
    &=U^{(2)}_{\Tilde{c}_\bu(\pi_{\Wbf}(\Tilde{g}_2)x,\Tilde{g}_1)}A^{(2)}_{\Tilde{c}_\ba(\Tilde{g}_1)}\Tilde{\varpi}(\Tilde{g}_1)U^{(2)}_{\Tilde{c}_\bu(x,\Tilde{g}_2)}A^{(2)}_{\Tilde{c}_\ba(\Tilde{g}_2)}\Tilde{\varpi}(\Tilde{g}_2)\psi(x)\\
    &=U^{(2)}_{\Tilde{c}_\bu(\pi_{\Wbf}(\Tilde{g}_2)x,\Tilde{g}_1)+\Tilde{c}_\bu(x,\Tilde{g}_2)e^{2\Tilde{c}_\ba(\Tilde{g}_1)}}A^{(2)}_{\Tilde{c}_\ba(\Tilde{g})}\Tilde{\varpi}(\Tilde{g})\psi(x),
\end{align*}
for all $\Tilde{g}_2\in K(\Tilde{g})$, $\Tilde{g}_1=\Tilde{g}\Tilde{g}_2^{-1}$. With this and \ref{item:BehaviorOfpsiA_t:1}, it follows that
\[
\Tilde{c}_\bu(x,\Tilde{g})=\Tilde{c}_\bu(\pi_{\Wbf}(\Tilde{g}_2)x,\Tilde{g}_1)+\Tilde{c}_\bu(x,\Tilde{g}_2)e^{2\Tilde{c}_\ba(\Tilde{g}_1)},\quad\Tilde{g}_2\in K(\Tilde{g}),\Tilde{g}_1=\Tilde{g}\Tilde{g}_2^{-1},
\]
is well-defined for all $x\in\Zold{Z:BehaviorOfpsiW_t:psi(pi_W(tilde(g)x))=...psi(x)}(\Tilde{g})$, $\Tilde{g}\in\left( B^{\widetilde{\Wbf}}_{\deold{delta:BehaviorOfpsiW_t:B^(tilde(W))_delta==B^W_delta}}\right)^{n+1}\setminus \left( B^{\widetilde{\Wbf}}_{\deold{delta:BehaviorOfpsiW_t:B^(tilde(W))_delta==B^W_delta}}\right)^n$. Then, by the induction hypothesis, we deduce that the remaining requirements of \textit{Lemma~\ref{Lemma:BehaviorOfpsiW_t:GoodSet(tilde(g))Andc_u(.,tilde(g))}} hold when $\widetilde{\Wbf}$ is replaced by $\left( B^{\widetilde{\Wbf}}_{\deold{delta:BehaviorOfpsiW_t:B^(tilde(W))_delta==B^W_delta}}\right)^{n+1}$.

Thus, by induction, we finish the proof of \textit{Lemma~\ref{Lemma:BehaviorOfpsiW_t:GoodSet(tilde(g))Andc_u(.,tilde(g))}}.
\end{proof}
From \textit{Lemma~\ref{Lemma:BehaviorOfpsiW_t:GoodSet(tilde(g))Andc_u(.,tilde(g))}} and the simple fact (cf. \ref{item:BehaviorOfpsiW_t:7})
\begin{equation}\label{eq:BehaviorOfpsiW_t:W_1x=x}
    \gamma x=x, \qquad\forall\gamma\in\Gamma_{1,\nc},\quad x\in\Gbf_1/\Gamma_1,
\end{equation}
we show the following:
\begin{Lemma}\label{Lemma:BehaviorOfpsiW_t:tilde(c)_a,tilde(varpi)(gamma)}

    For $\Tilde{\gamma}\in\Gamma_{\widetilde{\Wbf}}=\pi_\Wbf^{-1}(\Gamma_{1,\nc})$,
    \begin{gather*}
        \Tilde{c}_\ba(\Tilde{\gamma})=0,\qquad\Tilde{\varpi}(\Tilde{\gamma})\in\Gamma_2\cap Z(\Gbf_2),\\
        \Tilde{c}_\bu(x,\Tilde{\gamma})=0,\qquad\forall x\in\Zold{Z:BehaviorOfpsiW_t:psi(pi_W(tilde(g)x))=...psi(x)}(\Tilde{\gamma}).
    \end{gather*}
    Consequently, $\Tilde{c}_\ba\equiv0$. Moreover, we have, for $\Tilde{g}\in\widetilde{\Wbf}$, $x\in\Zold{Z:BehaviorOfpsiW_t:psi(pi_W(tilde(g)x))=...psi(x)}(\Tilde{g})$,
    \[
    \psi(\pi_{\Wbf}(\Tilde{g})x)=U^{(2)}_{\Tilde{c}_\bu(x,\Tilde{g})}\Tilde{\varpi}(\Tilde{g})\psi(x);
    \]
    and for $\Tilde{g}_1$, $\Tilde{g}_2\in\widetilde{\Wbf}$, $x\in\pi_{\Wbf}(\Tilde{g}_2)^{-1}\Zold{Z:BehaviorOfpsiW_t:psi(pi_W(tilde(g)x))=...psi(x)}(\Tilde{g}_1)\cap\Zold{Z:BehaviorOfpsiW_t:psi(pi_W(tilde(g)x))=...psi(x)}(\Tilde{g}_2)\cap\Zold{Z:BehaviorOfpsiW_t:psi(pi_W(tilde(g)x))=...psi(x)}(\Tilde{g}_1\Tilde{g}_2)$,
    \[
    \Tilde{c}_\bu(x,\Tilde{g}_1\Tilde{g}_2)=\Tilde{c}_\bu(\pi_{\Wbf}(\Tilde{g}_2)x,\Tilde{g}_1)+\Tilde{c}_\bu(x,\Tilde{g}_2).
    \]
    
\end{Lemma}
\begin{proof}
Fix $\Tilde{\gamma}\in\Gamma_{\widetilde{\Wbf}}=\pi_\Wbf^{-1}(\Gamma_{1,\nc})$. By \textit{Lemma~\ref{Lemma:BehaviorOfpsiW_t:GoodSet(tilde(g))Andc_u(.,tilde(g))}} and \eqref{eq:BehaviorOfpsiW_t:W_1x=x}, it follows that
\begin{equation}\label{eq:BehaviorOfpsiW_t:psi(x)=psi(pi_W(gamma)x)=...psi(x)}
    \psi(x)=\psi(\pi_\Wbf(\Tilde{\gamma})x)=U^{(2)}_{\Tilde{c}_\bu(x,\Tilde{\gamma})}A^{(2)}_{\Tilde{c}_\ba(\Tilde{\gamma})}\Tilde{\varpi}(\Tilde{\gamma})\psi(x),\quad\forall x\in\Zold{Z:BehaviorOfpsiW_t:psi(pi_W(tilde(g)x))=...psi(x)}(\Tilde{\gamma}).
\end{equation}
This and \ref{item:BehaviorOfpsis:8} imply that 
\begin{equation}\label{eq:BehaviorOfpsiW_t:A_(tilde(c)_a(gamma))tilde(varpi)(gamma)=U_tgGammag^(-1)}
    \textstyle{A^{(2)}_{\Tilde{c}_\ba(\Tilde{\gamma})}\Tilde{\varpi}(\Tilde{\gamma})\in\bigcap_{m_{\Gbf_2}\text{-a.e. }g\in\Gbf_2}\left(U^{(2)}_\R\cdot g\Gamma_2g^{-1}\right).}
\end{equation}
The intersection on the right-hand side is determined as follows:
\begin{Claim}\label{Claim:BehaviorOfpsiW_t:IntersectionOfU_R.gGammag^(-1)}

    $\bigcap_{m_{\Gbf_2}\text{-a.e. }g\in\Gbf_2}\left(U^{(2)}_\R\cdot g\Gamma_2g^{-1}\right)=U^{(2)}_\R\cdot\left(\Gamma_2\cap Z(\Gbf_2)\right)$.
    
\end{Claim}
\begin{proof}
Fix $g_0\in\bigcap_{m_{\Gbf_2}\text{-a.e. }g\in\Gbf_2}\left(U^{(2)}_\R\cdot g\Gamma_2g^{-1}\right)$. We will show that 
\[
g_0\in U^{(2)}_\R\cdot\left(\Gamma_2\cap Z(\Gbf_2)\right).
\]
For $\sigma\in\Gamma_2$, denote 
\[
E_\sigma=\left\{g\in\Gbf_2:g\sigma g^{-1}\in U^{(2)}_\R\cdot g_0\right\},
\]
it follows that $\cup_{\sigma\in\Gamma_2}E_\sigma$ is of full measure in $\Gbf_2$. Hence, there exists $\sigma_0\in\Gamma_2$ such that $m_{\Gbf_2}(E_{\sigma_0})>0$. Since $E_{\sigma_0}$ is the set of vanishing points of some polynomials on $\Gbf_2$, it follows that these polynomials vanish on a subset of $\Gbf_2$ having positive Haar measure. By an argument using \cite{MR4070868}, we deduce that these polynomials vanish on $\Gbf_2$, i.e. $E_{\sigma_0}=\Gbf_2$. This implies $g_0=U^{(2)}_{t_0}\sigma_0$ for some $t_0\in\R$, and so
\[
g\sigma_0 g^{-1}\sigma_0^{-1}\in U^{(2)}_\R\subseteq\Hbf_2=\langle\exp(\hfrak_2)\rangle,\qquad\forall g\in\Gbf_2.
\]
Thus, $\langle[\Gbf_2,\sigma_0]\rangle$ is a connected normal subgroup of $\Hbf_2$ contained in $U^{(2)}_\R$, which implies $[\Gbf_2,\sigma_0]=\{\ebf_{\Gbf_2}\}$, i.e. $\sigma_0\in Z(\Gbf_2)$. This shows that 
\[
g_0=U^{(2)}_{t_0}\sigma_0\in U^{(2)}_\R\cdot\left(\Gamma_2\cap Z(\Gbf_2)\right).
\]
Hence, $\bigcap_{m_{\Gbf_2}\text{-a.e. }g\in\Gbf_2}\left(U^{(2)}_\R\cdot g\Gamma_2g^{-1}\right)\subseteq U^{(2)}_\R\cdot\left(\Gamma_2\cap Z(\Gbf_2)\right)$. And since the reversed containment is obvious, we finish the proof of \textit{Claim \ref{Claim:BehaviorOfpsiW_t:IntersectionOfU_R.gGammag^(-1)}}.
\end{proof}
By \eqref{eq:BehaviorOfpsiW_t:A_(tilde(c)_a(gamma))tilde(varpi)(gamma)=U_tgGammag^(-1)} and \textit{Claim~\ref{Claim:BehaviorOfpsiW_t:IntersectionOfU_R.gGammag^(-1)}}, we have, for some $r_{\Tilde{\gamma}}\in\R$ and $\sigma_{\Tilde{\gamma}}\in\Gamma_2\cap Z(\Gbf_2)$,
\[
U^{(2)}_{r_{\Tilde{\gamma}}}A^{(2)}_{\Tilde{c}_\ba(\Tilde{\gamma})}\Tilde{\varpi}(\Tilde{\gamma})=\sigma_{\Tilde{\gamma}}.
\]
To both sides of this equation, conjugating by $A^{(2)}_p$ gives us $r_{\Tilde{\gamma}}=0$, then conjugating by $U^{(2)}_t$ gives us
\[
\Tilde{c}_\ba(\Tilde{\gamma})=0,\qquad\Tilde{\varpi}(\Tilde{\gamma})=\sigma_{\Tilde{\gamma}}\in\Gamma_2\cap Z(\Gbf_2).
\]
This, \eqref{eq:BehaviorOfpsiW_t:psi(x)=psi(pi_W(gamma)x)=...psi(x)}, and \ref{item:BehaviorOfpsiA_t:1} give us the first part of \textit{Lemma~\ref{Lemma:BehaviorOfpsiW_t:tilde(c)_a,tilde(varpi)(gamma)}}.

\medskip

From the facts that $\Tilde{c}_\ba:\widetilde{\Wbf}\to\R$ is a continuous homomorphism, that $\Tilde{c}_\ba(\Tilde{\gamma})=0$ for all $\Tilde{\gamma}\in\Gamma_{\widetilde{\Wbf}}$ shown above,  and that $\Gamma_{\widetilde{\Wbf}}=\pi_{\Wbf}^{-1}(\Gamma_{1,\nc})$, we obtain $\Tilde{c}_\ba\equiv0$. This, \textit{Lemma~\ref{Lemma:BehaviorOfpsiW_t:GoodSet(tilde(g))Andc_u(.,tilde(g))}}, and \ref{item:BehaviorOfpsiA_t:1} give us the last part of \textit{Lemma~\ref{Lemma:BehaviorOfpsiW_t:tilde(c)_a,tilde(varpi)(gamma)}}.

\medskip

Thus, we finish the proof of \textit{Lemma~\ref{Lemma:BehaviorOfpsiW_t:tilde(c)_a,tilde(varpi)(gamma)}}.
\end{proof}
With \textit{Lemma~\ref{Lemma:BehaviorOfpsiW_t:tilde(c)_a,tilde(varpi)(gamma)}}, we conclude the proof of \textit{Proposition~\ref{Prop:BehaviorOfpsiW_t:psi(W_tx)=...psi(x)}} as follows:
\begin{proof}[Proof of \textit{Proposition~\ref{Prop:BehaviorOfpsiW_t:psi(W_tx)=...psi(x)}}]
Let $\Tilde{\varpi}:\widetilde{\Wbf}\to\Zbf_2=C_{\Gbf_2}(\Hbf_2)$ be given in the arguments preceding \textit{Lemma~\ref{Lemma:BehaviorOfpsiW_t:GoodSet(tilde(g))Andc_u(.,tilde(g))}}. By \textit{Lemma~\ref{Lemma:BehaviorOfpsiW_t:tilde(c)_a,tilde(varpi)(gamma)}}, we have
\[
\Tilde{\varpi}(\Tilde{\gamma})\in\Gamma_2\cap Z(\Gbf_2),\qquad\forall\Tilde{\gamma}\in\Gamma_{\widetilde{\Wbf}}.
\]

For $g\in\Wbf$, define
\[
\textstyle{\Zold{Z:BehaviorOfpsiW_t:psi(W_tx)=U_(f_j(x,s))g(s)psi(x)}(g)=\bigcap_{\Tilde{g}\in\pi_{\Wbf}^{-1}(g\Gamma_{1,\nc})}}\Zold{Z:BehaviorOfpsiW_t:psi(pi_W(tilde(g)x))=...psi(x)}(\Tilde{g})\cap\bigcap_{\Tilde{\gamma}\in\Gamma_{\widetilde{\Wbf}}}\Zold{Z:BehaviorOfpsiW_t:psi(pi_W(tilde(g)x))=...psi(x)}(\Tilde{\gamma}).
\]
Then, by \textit{Lemmas~\ref{Lemma:BehaviorOfpsiW_t:GoodSet(tilde(g))Andc_u(.,tilde(g))}}, it follows that $\Zold{Z:BehaviorOfpsiW_t:psi(W_tx)=U_(f_j(x,s))g(s)psi(x)}(g)\subseteq\Gbf_1/\Gamma_1$ is a full measure $U^{(1)}_t$-invariant set for all $g\in\Wbf$, $\{(x,g):x\in\Zold{Z:BehaviorOfpsiW_t:psi(W_tx)=U_(f_j(x,s))g(s)psi(x)}(g),g\in\Wbf\}$ is of full measure in $(\Gbf_1/\Gamma_1)\times\Wbf$, and for all $g\in\Wbf$,
\[
\Zold{Z:BehaviorOfpsiW_t:psi(W_tx)=U_(f_j(x,s))g(s)psi(x)}(g)\subseteq\Zold{Z:BehaviorOfpsiUnderA_t:GenericPointsOfGoodSetInLemmapsi(Kak)=Kak}\cap g^{-1}\Zold{Z:BehaviorOfpsiUnderA_t:GenericPointsOfGoodSetInLemmapsi(Kak)=Kak},\qquad\Zold{Z:BehaviorOfpsiW_t:psi(W_tx)=U_(f_j(x,s))g(s)psi(x)}(g)=\Zold{Z:BehaviorOfpsiW_t:psi(W_tx)=U_(f_j(x,s))g(s)psi(x)}(g\gamma)\subseteq\Zold{Z:BehaviorOfpsiW_t:psi(W_tx)=U_(f_j(x,s))g(s)psi(x)}(e_{\Gbf_1}),\quad\forall\gamma\in\Gamma_{1,\nc}.
\]

For $g\in\Wbf$ and $x\in\Zold{Z:BehaviorOfpsiW_t:psi(W_tx)=U_(f_j(x,s))g(s)psi(x)}(g)$, we define
\[
f(x,g)=\Tilde{c}_\bu(x,\Tilde{g}),\qquad\Tilde{g}\in\pi_\Wbf^{-1}(g).
\]
By \textit{Lemma~\ref{Lemma:BehaviorOfpsiW_t:tilde(c)_a,tilde(varpi)(gamma)}}, we deduce that $f:\{(x,g):x\in\Zold{Z:BehaviorOfpsiW_t:psi(W_tx)=U_(f_j(x,s))g(s)psi(x)}(g),g\in\Wbf\}\to\R$ is a well-defined measurable function such that for $g\in\Wbf$, $\gamma\in\Gamma_{1,\nc}$,
\begin{gather*}
    f(x,\gamma)=0,\quad\forall x\in\Zold{Z:BehaviorOfpsiW_t:psi(W_tx)=U_(f_j(x,s))g(s)psi(x)}(\gamma);\qquad f(x,g\gamma)=f(x,g),\quad\forall x\in\Zold{Z:BehaviorOfpsiW_t:psi(W_tx)=U_(f_j(x,s))g(s)psi(x)}(g);\\
    \psi(gx)=U^{(2)}_{f(x,g)}\Tilde{\varpi}(\Tilde{g})\psi(x),\qquad\forall x\in\Zold{Z:BehaviorOfpsiW_t:psi(W_tx)=U_(f_j(x,s))g(s)psi(x)}(g),\quad\Tilde{g}\in\pi_\Wbf^{-1}(g).
\end{gather*}
With this last equation, \eqref{eq:BehaviorOfpsi:ConjugationRelation}, and \ref{item:BehaviorOfpsiA_t:1}, we obtain
\[
f(U^{(1)}_rx,g)-f(x,g)=\omega(gx,r)-\omega(x,r),\quad\forall r\in\R,\quad g\in\Wbf,\quad x\in\Zold{Z:BehaviorOfpsiW_t:psi(W_tx)=U_(f_j(x,s))g(s)psi(x)}(g).
\]

Thus, we finish the proof of \textit{Proposition~\ref{Prop:BehaviorOfpsiW_t:psi(W_tx)=...psi(x)}}.
\end{proof}

\subsection{Reduction via solving cohomological equations}\label{subSection:BehaviorOfpsiW_t:SolvingCohomologicalEquation}
In this subsection, we find a measurable $U^{(2)}_t$-translation to ``cancel" the $U^{(2)}_t$-translation in $\psi(g(\cdot))$, i.e. ``cancel" $U^{(2)}_{f_j(x,s)}$ in the relation between $\psi(g(\cdot))$ and $\psi$ in \textit{Proposition~\ref{Prop:BehaviorOfpsiW_t:psi(W_tx)=...psi(x)}}, and thus proving \textit{Proposition~\ref{Prop:BehaviorOfpsiW_t}}.

\medskip

Before going into the proof of \textit{Proposition~\ref{Prop:BehaviorOfpsiW_t}}, let us (heuristically) describe the main idea. By \textit{Proposition~\ref{Prop:BehaviorOfpsiW_t:psi(W_tx)=...psi(x)}}, we see that for a measurable function $\Tilde{b}:\Gbf_1/\Gamma_1\to\R$,
\[
U^{(2)}_{\Tilde{b}(gx)}\psi(gx)=U^{(2)}_{\Tilde{b}(gx)+f(x,g)-\Tilde{b}(x)}\Tilde{\varpi}(\Tilde{g})\left(U^{(2)}_{\Tilde{b}(x)}\psi(x)\right),\quad\Tilde{g}\in\pi_\Wbf^{-1}(g).
\]
Hence, to ``cancel" the $U^{(2)}_t$-translations, we look to solve the cohomological equation
\begin{equation}\label{eq:BehaviorOfpsiW_t:CohomologicalEquation}
    \Tilde{b}(x)-\Tilde{b}(gx)=f(x,g).
\end{equation}
Notice that by \textit{Proposition~\ref{Prop:BehaviorOfpsiW_t:psi(W_tx)=...psi(x)}} and \eqref{eq:BehaviorOfpsiW_t:W_1x=x}, $f$ can be considered as a cocycle $(\Gbf_1/\Gamma_1)\times(\Wbf/\Gamma_{1,\nc})\to\R$. This is necessary for \eqref{eq:BehaviorOfpsiW_t:CohomologicalEquation} to have a solution. Moreover, it also implies that $\{f(x,\cdot)_*m_{\Wbf/\Gamma_{1,\nc}}:x\in\Gbf_1/\Gamma_1\}$, where $m_{\Wbf/\Gamma_{1,\nc}}$ is the probability Haar measure on $\Wbf/\Gamma_{1,\nc}$, is a family of Borel probability measures on $\R$ such that
\[
(f(gx,\cdot)_*m_{\Wbf/\Gamma_{1,\nc}})(\cdot)=(f(x,\cdot)_*m_{\Wbf/\Gamma_{1,\nc}})(\cdot+f(x,g)).
\]
With this property, following a construction in \cite{MR1784202}*{p.~365}, we get
\begin{equation}\label{eq:BehaviorOfpsiW_t:IdeaForConstructingb_j}
    \textstyle{\Tilde{b}(x)=\inf\{t\in\R:(f(x,\cdot)_*m_{\Wbf/\Gamma_{1,\nc}})((-\infty,t])\geq\frac{1}{2}\}}
\end{equation}
is a solution to (\ref{eq:BehaviorOfpsiW_t:CohomologicalEquation}). On the other hand, we also want $\Tilde{\psi}(x)=U^{(2)}_{\Tilde{b}(x)}\psi(x)$ to be a Kakutani factor map between $(\Gbf_i/\Gamma_i,\mu_i,U^{(i)}_t),i=1,2$, with a $10^{-20}\eold{epsilon:Preliminaries:CoveringLemma}$-good time change $\Tilde{\tau}$. Hence, from
\[
\Tilde{\psi}(U^{(1)}_rx)=U^{(2)}_{\Tilde{b}(U^{(1)}_rx)-\Tilde{b}(x)+\omega(x,r)}\Tilde{\psi}(x),
\]
we want $\Tilde{b}$ to satisfy that $\Tilde{b}(U^{(1)}_rx)$ is $C^\infty$ in $r\in\R$. To obtain this desired regularity, we use the classical convolution techniques. More specifically, let $\phi$ be a positive smooth function on $\R$ such that $\int_\R\phi\,dx=1$, and consider the new family $\{(f(x,\cdot)_*m_{\Wbf/\Gamma_{1,\nc}})*(\phi dm_\R):x\in\Gbf_1/\Gamma_1\}$ of probability measures. With this family, the construction in (\ref{eq:BehaviorOfpsiW_t:IdeaForConstructingb_j}) becomes
\[
\Tilde{b}(x)\in\R\text{ such that }\int_{\Wbf/\Gamma_{1,\nc}}\varphi(\Tilde{b}(x)-f(x,g))dm_{\Wbf/\Gamma_{1,\nc}}(g)=\frac{1}{2},
\]
where $\varphi(t)=\int_{-\infty}^t\phi(s)ds$. By the \textit{Implicit function theorem}, we see that $\displaystyle{\left.\frac{d}{dr}\right\vert_{r=r_0}\Tilde{b}(U^{(1)}_rx)}$ exists and is continuous in $r$. Then, we perform another convolution to obtain $C^\infty$ regularity.

\medskip

Now, we proceed to the rigorous proof of \textit{Proposition~\ref{Prop:BehaviorOfpsiW_t}}, where we first construct $\Tilde{b}$, then use it to define $\Tilde{\psi}$.

\medskip

Recalling $\Zold{Z:BehaviorOfpsiW_t:psi(W_tx)=U_(f_j(x,s))g(s)psi(x)}(g)$ given in \textit{Proposition~\ref{Prop:BehaviorOfpsiW_t:psi(W_tx)=...psi(x)}}, we define two auxiliary sets
\[
\Znew\label{Z:BehaviorOfpsiW_t:FubiniConsequence}=\{x\in\Gbf_1/\Gamma_1:\Wbf_x=\{g\in\Wbf:x\in\Zold{Z:BehaviorOfpsiW_t:psi(W_tx)=U_(f_j(x,s))g(s)psi(x)}(g)\}\text{ is of full measure in }\Wbf\},
\]
\[    \Znew\label{Z:BehaviorOfpsiW_t:DomainOfSolForCohomologicalEq}=\{gx:x\in\Zold{Z:BehaviorOfpsiW_t:FubiniConsequence},g\in\Wbf_x\}.
\]
It is clear that $\Zold{Z:BehaviorOfpsiW_t:FubiniConsequence}$ is a full measure $U^{(1)}_t$-invariant set in $\Zold{Z:BehaviorOfpsiW_t:psi(W_tx)=U_(f_j(x,s))g(s)psi(x)}(e_{\Gbf_1})$ and
\[
\{(x,g):x\in\Zold{Z:BehaviorOfpsiW_t:FubiniConsequence},g\in \Wbf_x\}
\]
is a full measure set in $(\Gbf_1/\Gamma_1)\times\Wbf$, contained in 
\[
\{(x,g):x\in\Zold{Z:BehaviorOfpsiW_t:psi(W_tx)=U_(f_j(x,s))g(s)psi(x)}(g),g\in\Wbf\}.
\]
With $\Zold{Z:BehaviorOfpsiW_t:DomainOfSolForCohomologicalEq}$ as the domain, we can rigorously construct a solution $\hat{b}$ for the cohomological equation (\ref{eq:BehaviorOfpsiW_t:CohomologicalEquation}) that has $\displaystyle{\left.\frac{d}{dr}\right\vert_{r=r_0}\hat{b}(U^{(1)}_rx)}$ satisfying certain conditions. This serves as the ``base" for our convolution construction of $\Tilde{b}$ in \textit{Proposition~\ref{Prop:BehaviorOfpsiW_t}}.

Recall from \ref{item:BehaviorOfpsiW_t:8} that $F_\Wbf\subseteq\Wbf$ is a fundamental domain for $\Gamma_{1,\nc}$ and $m_{\Wbf}$ is a Haar measure on $\Wbf$ such that $m_{\Wbf}(F_\Wbf)=1$. For notation simplicity, we only write $dg$ in place of $dm_{\Wbf}(g)$. And for $g\in\Wbf$, we write $\Tilde{g}$ to denote an element in $\pi_\Wbf^{-1}(g)$, in the remainder of this section, starting with the following lemma:
\begin{Lemma}\label{Lemma:BehaviorOfpsiW_t:SolveCohomologicalEquation}

    Let $\phi=\frac{1}{\sqrt{2\pi}}e^{-t^2/2}$ and $\varphi(t)=\int_{-\infty}^t\phi(s)\,ds$. For $
    x\in\Zold{Z:BehaviorOfpsiW_t:DomainOfSolForCohomologicalEq}$, we define $\hat{b}(x)$ by
    \begin{equation}\label{eq:BehaviorOfpsiW_t:DefOfSolForCohomologicalEq}
        \int_{F_\Wbf}\varphi\left(\hat{b}(x)-(f(y,gh)-f(y,h))\right)dg=\frac{1}{2},
    \end{equation}
    where $y\in\Zold{Z:BehaviorOfpsiW_t:FubiniConsequence},h\in \Wbf_y$ satisfy $x=hy$.
    Then $\hat{b}:\Zold{Z:BehaviorOfpsiW_t:DomainOfSolForCohomologicalEq}\to\R$ is a well-defined measurable function such that for $x\in\Zold{Z:BehaviorOfpsiW_t:FubiniConsequence}$,
    \begin{equation*}
        \begin{gathered}
            \hat{b}(x)-\hat{b}(hx)=f(x,h),\quad \forall h\in \Wbf_x;\\
            \left.\frac{d}{dr}\right\vert_{r=r_0}\hat{b}(U^{(1)}_rx)=\hat{\tau}(U^{(1)}_{r_0}x)-\tau(U^{(1)}_{r_0}x),\quad \forall r_0\in\R;
        \end{gathered}
    \end{equation*}
    where $\hat{\tau}:\Zold{Z:BehaviorOfpsiW_t:FubiniConsequence}\to\R_{>0}$ is a measurable positive function, continuous along $U^{(1)}_t$-orbits, satisfying
    \[
    \norm{\hat{\tau}-1}_{\infty}<10^{-20}\eold{epsilon:Preliminaries:CoveringLemma},\qquad \int_{G_1/\Gamma_1}\hat{\tau}\,d\mu_1=1,
    \]
    and determined for every $x\in\Zold{Z:BehaviorOfpsiW_t:FubiniConsequence}$ by
    \begin{equation}\label{eq:BehaviorOfpsiW_t:DefOftau^(j)}
        \hat{\tau}(x)=\frac{\int_{F_\Wbf}\phi\left(\hat{b}(x)-f(x,g)\right)\tau(gx)dg}{\int_{F_\Wbf}\phi\left(\hat{b}(x)-f(x,g)\right)dg}=\frac{\int_{F_\Wbf}\phi\left(\hat{b}(gx)\right)\tau(gx)dg}{\int_{F_\Wbf}\phi\left(\hat{b}(gx)\right)dg}.
    \end{equation}
    
\end{Lemma}
\begin{proof}
First, we check that $\hat{b}:\Zold{Z:BehaviorOfpsiW_t:DomainOfSolForCohomologicalEq}\to\R$ given in (\ref{eq:BehaviorOfpsiW_t:DefOfSolForCohomologicalEq}) is a well-defined measurable function. For $x=hy$, $y\in\Zold{Z:BehaviorOfpsiW_t:FubiniConsequence}$, $h\in\Wbf_y$, 
\[
\Wbf\to\R:g\mapsto f(y,gh)-f(y,h)
\]
is a well-defined measurable function, and so
\[
\Wbf\to\R:\alpha\mapsto\int_{F_\Wbf}\varphi(\alpha-(f(y,gh)-f(y,h)))dg
\]
is well-defined continuous function, strictly increasing from $0$ to $1$. Hence, there exists a unique $\hat{b}(x)\in\R$ satisfying \eqref{eq:BehaviorOfpsiW_t:DefOfSolForCohomologicalEq}.

If $x=h_1y_1=h_2y_2$ for some $y_1,y_2\in\Zold{Z:BehaviorOfpsiW_t:FubiniConsequence},h_1\in\Wbf_{y_1},h_2\in\Wbf_{y_2}$, it follows from \textit{Proposition~\ref{Prop:BehaviorOfpsiW_t:psi(W_tx)=...psi(x)}} that for $g\in (\Wbf_{y_1}h_1^{-1})\cap(\Wbf_{y_2}h_2^{-1})$,
\begin{align*}
    &U^{(2)}_{f(y_1,gh_1)}\Tilde{\varpi}(\Tilde{g}\Tilde{h}_1)\psi(y_1)=\psi(gh_1y_1)=\psi(gh_2y_2)\\
    =\;&U^{(2)}_{f(y_2,gh_2)}\Tilde{\varpi}(\Tilde{g}\Tilde{h}_2)\psi(y_2)=U^{(2)}_{f(y_2,gh_2)-f(y_2,h_2)}\Tilde{\varpi}(\Tilde{g})\psi(h_2y_2)\\
    =\;&U^{(2)}_{f(y_2,gh_2)-f(y_2,h_2)}\Tilde{\varpi}(\Tilde{g})\psi(h_1y_1)=U^{(2)}_{f(y_2,gh_2)-f(y_2,h_2)+f(y_1,h_1)}\Tilde{\varpi}(\Tilde{g}\Tilde{h}_1)\psi(y_1).
\end{align*}
This and \ref{item:BehaviorOfpsiA_t:1} imply that
\[
f(y_1,gh_1)-f(y_1,h_1)=f(y_2,gh_2)-f(y_2,h_2),\quad \forall g\in (\Wbf_{y_1}h_1^{-1})\cap(\Wbf_{y_2}h_2^{-1}).
\]
Hence, $\hat{b}(x)$ given in (\ref{eq:BehaviorOfpsiW_t:DefOfSolForCohomologicalEq}) is independent of the choice of $y\in\Zold{Z:BehaviorOfpsiW_t:FubiniConsequence},h\in\Wbf_y$ satisfying $x=hy$. This implies that $\hat{b}(x)$ given in \eqref{eq:BehaviorOfpsiW_t:DefOfSolForCohomologicalEq} for $x\in\Zold{Z:BehaviorOfpsiW_t:DomainOfSolForCohomologicalEq}$ is well-defined, and by the formula in \eqref{eq:BehaviorOfpsiW_t:DefOfSolForCohomologicalEq}, we also deduce that $\hat{b}:\Zold{Z:BehaviorOfpsiW_t:DomainOfSolForCohomologicalEq}\to\R$ is measurable.

\medskip

Next, we show that $\hat{b}(x)-\hat{b}(hx)=f(x,h)$ for all $ x\in\Zold{Z:BehaviorOfpsiW_t:FubiniConsequence}$, $h\in\Wbf_x$.
For $x\in\Zold{Z:BehaviorOfpsiW_t:FubiniConsequence},h\in\Wbf_x$, by \eqref{eq:BehaviorOfpsiW_t:DefOfSolForCohomologicalEq}, it follows that
\[
\int_{F_\Wbf}\varphi\left(\hat{b}(x)-f(x,g)\right)dg=\frac{1}{2}=\int_{F_\Wbf}\varphi\left(\hat{b}(hx)-(f(x,gh)-f(x,h))\right)dg.
\]
On the other hand, \textit{Proposition~\ref{Prop:BehaviorOfpsiW_t:psi(W_tx)=...psi(x)}} implies that 
\[
f(x,g\gamma)=f(x,g),\qquad\forall g\in\Wbf_x,\quad\gamma\in\Gamma_{1,\nc},
\]
and hence, since $F_\Wbf\subseteq\Wbf$ is a fundamental domain for $\Gamma_{1,\nc}$, by \textit{Change of variables}, we have
\begin{align*}
    &\int_{F_\Wbf}\varphi\left(\hat{b}(hx)-(f(x,gh)-f(x,h))\right)dg\\
    =&\int_{F_\Wbf}\varphi\left((\hat{b}(hx)+f(x,h))-f(x,g)\right)dg.
\end{align*}
Thus, $\hat{b}(x)=\hat{b}(hx)+f(x,h)$, i.e. $\hat{b}(x)-\hat{b}(hx)=f(x,h)$.

\medskip

Now, we compute $\displaystyle{\left.\frac{d}{dr}\right\vert_{r=r_0}\hat{b}(U^{(1)}_rx)}$. For $x\in\Zold{Z:BehaviorOfpsiW_t:FubiniConsequence}$, by \textit{Proposition~\ref{Prop:BehaviorOfpsiW_t:psi(W_tx)=...psi(x)}},
\[
f(U^{(1)}_rx,g)-f(x,g)=\omega(gx,r)-\omega(x,r),\qquad\forall g\in\Wbf_x,\quad r\in\R.
\]
The above equation indeed follows from the fact that $\Zold{Z:BehaviorOfpsiW_t:psi(W_tx)=U_(f_j(x,s))g(s)psi(x)}(g)$ is $U^{(1)}_t$-invariant and then $\Wbf_x=\Wbf_{U^{(1)}_rx}$ for all $r\in\R$. This and (\ref{eq:BehaviorOfpsi:rho(x,s)andomega(x,s)}) imply that
\begin{equation}\label{eq:BehaviorOfpsiW_t:d/dr(f_j(U_rx,s))}
    \left.\frac{d}{dr}\right\vert_{r=r_0}f(U^{(1)}_rx,h)=\tau(hU^{(1)}_{r_0}x)-\tau(U^{(1)}_{r_0}x),\qquad \forall h\in\Wbf_x,\quad r_0\in\R.
\end{equation}
For fixed $x\in\Zold{Z:BehaviorOfpsiW_t:FubiniConsequence}$, consider
\[
F_x(r,a)=\int_{F_\Wbf}\varphi\left(a-f(U_r^{(1)}x,g)\right)dg.
\]
Since $\absolute{\partial_r f(U_r^{(1)}x,g)}\leq 2\cdot10^{-10}\epsilon_2$ (cf. \eqref{eq:BehaviorOfpsi:|partial_svarpi_i(x,s)-1|IsSmall} and \eqref{eq:BehaviorOfpsiW_t:d/dr(f_j(U_rx,s))}) and $\varphi$ is bounded, it follows from the \textit{Dominated convergence theorem} that $F_x$ is $C^1$ and
\[
\partial_aF_x(r,a)=\int_{F_\Wbf}\phi\left(a-f(U_r^{(1)}x,g)\right)dg>0.
\]
The \textit{Implicit function theorem} then implies that for $x\in\Zold{Z:BehaviorOfpsiW_t:FubiniConsequence}$, $\hat{b}(U^{(1)}_rx)$ is $C^1$ in $r\in\R$, and for $r_0\in\R$,
\[
\int_{F_\Wbf}\phi\left(\hat{b}(U^{(1)}_{r_0}x)-f(U^{(1)}_{r_0}x,g)\right)\left(\left.\frac{d}{dr}\right\vert_{r=r_0}\hat{b}(U^{(1)}_rx)-\tau(gU^{(1)}_{r_0}x)+\tau(U^{(1)}_{r_0}x)\right)dg=0.
\]
Thus, 
\[
\displaystyle{\left.\frac{d}{dr}\right\vert_{r=r_0}\hat{b}(U^{(1)}_rx)=\hat{\tau}(U^{(1)}_{r_0}x)-\tau(U^{(1)}_{r_0}x)},
\]
where $\hat{\tau}$ is determined as in \eqref{eq:BehaviorOfpsiW_t:DefOftau^(j)}.

\medskip

Lastly, we verify the remaining properties of $\hat{\tau}$. Since $\hat{b},f,\tau$ are measurable and $\phi$ is $C^{\infty}$, the formula \eqref{eq:BehaviorOfpsiW_t:DefOftau^(j)} implies that $\hat{\tau}$ is measurable. It follows from \eqref{eq:BehaviorOfpsi:|partial_svarpi_i(x,s)-1|IsSmall} and the constructions of $\Zold{Z:BehaviorOfpsiW_t:FubiniConsequence}$ and $\Zold{Z:BehaviorOfpsiW_t:psi(W_tx)=U_(f_j(x,s))g(s)psi(x)}(\cdot)$ that for every $x\in\Zold{Z:BehaviorOfpsiW_t:FubiniConsequence}$, $g\in\Wbf_x$, $r\in\R$,
\[
[\R\to\R:r\mapsto\tau(gU^{(1)}_rx)]\in C^\infty(\R),\quad\absolute{\tau(gU^{(1)}_rx)-1}<10^{-20}\eold{epsilon:Preliminaries:CoveringLemma}.
\]
This, \eqref{eq:BehaviorOfpsiW_t:d/dr(f_j(U_rx,s))}, the fact $\hat{b}(U^{(1)}_rx)$ is $C^1$ in $r\in\R$ for $x\in\Zold{Z:BehaviorOfpsiW_t:FubiniConsequence}$ (shown above), the \textit{Dominated convergence theorem}, and the formula of $\hat{\tau}$ in \eqref{eq:BehaviorOfpsiW_t:DefOftau^(j)}, imply that for every $x\in\Zold{Z:BehaviorOfpsiW_t:FubiniConsequence}$ and $r\in\R$,
\begin{equation}\label{eq:BehaviorOfpsiW_t:|tau^(j)-1|<epsilon}
    \left[\R\to\R:r\mapsto\hat{\tau}(U^{(1)}_rx)\right]\in C^0(\R),\quad\absolute{\hat{\tau}(U^{(1)}_rx)-1}<10^{-20}\eold{epsilon:Preliminaries:CoveringLemma},.
\end{equation}
Since $\displaystyle{\left.\frac{d}{dr}\right\vert_{r=r_0}\hat{b}(U^{(1)}_rx)}=\hat{\tau}(U^{(1)}_{r_0}x)-\tau(U^{(1)}_{r_0}x)$ for all $r_0\in\R,x\in\Zold{Z:BehaviorOfpsiW_t:FubiniConsequence}$ (shown above), $\hat{\tau},\tau\in L^\infty(\Gbf_1/\Gamma_1,\mu_1)$ (cf. \eqref{eq:BehaviorOfpsiW_t:|tau^(j)-1|<epsilon}, \eqref{eq:BehaviorOfpsi:|partial_svarpi_i(x,s)-1|IsSmall}), and the flow $U^{(1)}_t$ on $(\Gbf_1/\Gamma_1,\mu_1)$ is ergodic, it follows from the \textit{Pointwise ergodic theorem} that for $\mu_1$-a.e. $x\in\Zold{Z:BehaviorOfpsiW_t:FubiniConsequence}$,
\begin{multline*}
    \frac{1}{r}\left(\hat{b}(U^{(1)}_rx)-\hat{b}(x)\right)=\frac{1}{r}\int_0^r\hat{\tau}(U^{(1)}_sx)ds-\frac{1}{r}\int_0^r\tau(U^{(1)}_sx)ds\\
    \xrightarrow[r\to\infty]{}\int_{\Gbf_1/\Gamma_1}\hat{\tau}d\mu_1-\int_{\Gbf_1/\Gamma_1}\tau d\mu_1=\int_{\Gbf_1/\Gamma_1}\hat{\tau}d\mu_1-1.
\end{multline*}
On the other hand, we can choose one such $x\in\Zold{Z:BehaviorOfpsiW_t:FubiniConsequence}$ and $r_n\xrightarrow[n\to\infty]{}\infty$ such that $\frac{1}{r_n}\hat{b}(U^{(1)}_{r_n}x)\xrightarrow[n\to\infty]{}0$, since $\hat{b}:\Zold{Z:BehaviorOfpsiW_t:DomainOfSolForCohomologicalEq}\to\R$ is measurable (shown above) and the flow $U^{(1)}_t$ on $(\Gbf_1/\Gamma_1,\mu_1)$ is ergodic. Hence, it follows that $\int_{\Gbf_1/\Gamma_1}\hat{\tau}d\mu_1=1$.

This finishes the proof of \textit{Lemma \ref{Lemma:BehaviorOfpsiW_t:SolveCohomologicalEquation}}.
\end{proof}
With the solution $\hat{b}$ for (\ref{eq:BehaviorOfpsiW_t:CohomologicalEquation}) constructed, we prove \textit{Proposition~\ref{Prop:BehaviorOfpsiW_t}}:
\begin{proof}[Proof of \textit{Proposition~\ref{Prop:BehaviorOfpsiW_t}}]
Recall $\hat{b}$, $\hat{\tau}$ in \textit{Lemma~\ref{Lemma:BehaviorOfpsiW_t:SolveCohomologicalEquation}}. Define the sets
\[
\Zold{Z:BehaviorOfpsiW_t:ResultOfInduction}=\textstyle{\bigcap_{n\in\Z}}A_n^{(1)}(\Zold{Z:BehaviorOfpsiW_t:FubiniConsequence}\cap\Zold{Z:BehaviorOfpsiA_t:GoodSetForpsi(A_px)=U_()exp(p(a+X_0))psi(x)}),\qquad Z=\{hx:x\in\Zold{Z:BehaviorOfpsiW_t:ResultOfInduction},h\in\Wbf\}.
\]
Recall $\varpi$ in \textit{Proposition~\ref{Prop:BehaviorOfpsiW_t:psi(W_tx)=...psi(x)}}. Define the map $\hat{\psi}:Z\to\Gbf_2/\Gamma_2$ by
\begin{equation}\label{eq:BehaviorOfpsiW_t:DefOfpsi^(j)}
    \hat{\psi}(hx)=\Tilde{\varpi}(\Tilde{h})U^{(2)}_{\hat{b}(x)}\psi(x),\qquad\forall x\in\Zold{Z:BehaviorOfpsiW_t:ResultOfInduction},\quad h\in\Wbf.
\end{equation}

We claim $\hat{\psi}$ is well-defined. Indeed, if $x=h_1x_1=h_2x_2$ for $x_1,x_2\in\Zold{Z:BehaviorOfpsiW_t:ResultOfInduction}$, $h_1,h_2\in\Wbf$, then as $\Wbf_{x_1},\Wbf_{x_2}\subseteq\Wbf$ are of full measure, we can choose some $g\in(\Wbf_{x_1}h_1^{-1})\cap(\Wbf_{x_2}h_2^{-1})$. By \textit{Lemma~\ref{Lemma:BehaviorOfpsiW_t:SolveCohomologicalEquation}}, it follows that
\begin{align*}
    &\psi(gx)=\psi(gh_1x_1)=U^{(2)}_{f(x_1,gh_1)-\hat{b}(x_1)}\Tilde{\varpi}(\Tilde{g})\left(\Tilde{\varpi}(\Tilde{h}_1)U^{(2)}_{\hat{b}(x_1)}\psi(x_1)\right),\\
    &\psi(gx)=\psi(gh_2x_2)=U^{(2)}_{f(x_2,gh_2)-\hat{b}(x_2)}\Tilde{\varpi}(\Tilde{g})\left(\Tilde{\varpi}(\Tilde{h}_2)U^{(2)}_{\hat{b}(x_2)}\psi(x_2)\right),\\
    &f(x_1,gh_1)-\hat{b}(x_1)=-\hat{b}(gh_1x_1)=-\hat{b}(gh_2x_2)=f(x_2,gh_2)-\hat{b}(x_2).
\end{align*}
Hence, $\Tilde{\varpi}(\Tilde{h}_1)U^{(2)}_{\hat{b}(x_1)}\psi(x_1)=\Tilde{\varpi}(\Tilde{h}_2)U^{(2)}_{\hat{b}(x_2)}\psi(x_2)$, and we conclude that $\hat{\psi}$ is well-defined.

\medskip

From the construction of $\hat{\psi}$, $\Zold{Z:BehaviorOfpsiW_t:ResultOfInduction}$, and \ref{item:BehaviorOfpsiW_t:1}, we deduce that $\hat{\psi}$ is measurable, and for $x\in\Zold{Z:BehaviorOfpsiW_t:ResultOfInduction}$,
\begin{equation}\label{eq:BehaviorOfpsiW_t:psi^(j)(W^j_sx)}
    \begin{aligned}
        &\hat{\psi}(gx)=\tilde{\varpi}(\tilde{g})\hat{\psi}(x),\qquad \forall g\in\Wbf;\\
        &\hat{\psi}(x)=U^{(2)}_{\hat{b}(x)}\psi(x);\qquad\hat{\psi}(A^{(1)}_1x)=A^{(2)}_1U^{(2)}_{\hat{f}(x)}\hat{\psi}(x);
    \end{aligned}
\end{equation}
where $\Tilde{\varpi}$ is defined in \eqref{eq:tildeVarpiDef} and $\hat{f}:\Zold{Z:BehaviorOfpsiW_t:ResultOfInduction}\to\R$ is a measurable function determined by
\[
\hat{f}(x)=f(x)+e^{-2}\hat{b}(A^{(1)}_1x)-\hat{b}(x).
\]

Next, we show that $\hat{\psi}$ satisfies
\begin{equation}\label{eq:BehaviorOfpsiW_t:psi^(j)(U_rx)=U_omega(x,r)psi^(j)(x)}
    \hat{\psi}(U^{(1)}_rx)=U^{(2)}_{\hat{\omega}(x,r)}\hat{\psi}(x),\qquad\forall x\in\Zold{Z:BehaviorOfpsiW_t:ResultOfInduction},\quad r\in\R;
\end{equation}
\begin{equation}\label{eq:BehaviorOfpsiW_t:psi^(j)PreservesCorrespondingMeasures}
    \hat{\psi}_*\mu_{\hat{\tau}}=\mu_2;
\end{equation}
where $d\mu_{\hat{\tau}}=\hat{\tau}d\mu_1$ and
\begin{equation}\label{eq:BehaviorOfpsiW_t:DefOfomega^(j)}
    \hat{\omega}(x,r)=\int_0^r\hat{\tau}(U^{(1)}_sx)ds.
\end{equation}

Indeed, for $x\in\Zold{Z:BehaviorOfpsiW_t:ResultOfInduction}$, $r\in\R$, since $\Zold{Z:BehaviorOfpsiW_t:ResultOfInduction}$ is $U^{(1)}_t$-invariant, it follows from \eqref{eq:BehaviorOfpsiW_t:DefOfpsi^(j)} and (\ref{eq:BehaviorOfpsi:ConjugationRelation}) that
\[
\begin{aligned}
\hat{\psi}(U^{(1)}_rx)&=U^{(2)}_{\hat{b}(U^{(1)}_rx)}\psi(U^{(1)}_rx)=U^{(2)}_{\hat{b}(U^{(1)}_rx)+\omega(x,r)}\psi(x)\\
&=U^{(2)}_{\hat{b}(U^{(1)}_rx)-\hat{b}(x)+\omega(x,r)}\hat{\psi}(x).
\end{aligned}
\]
From this and the fact that
\begin{equation}\label{eq:BehaviorOfpsiW_t:intOfhat(tau)dmu_1=1}
    \hat{b}(U^{(1)}_rx)-\hat{b}(x)+\omega(x,r)=\hat{\omega}(x,r),
\end{equation}
which follows from \eqref{eq:BehaviorOfpsiW_t:DefOfomega^(j)} and (cf. \textit{Lemma~\ref{Lemma:BehaviorOfpsiW_t:SolveCohomologicalEquation}})
\[
\left.\frac{d}{dr}\right\vert_{r=r_0}\hat{b}(U^{(1)}_rx)=\hat{\tau}(U^{(1)}_{r_0}x)-\tau(U^{(1)}_{r_0}x),
\]
we obtain \eqref{eq:BehaviorOfpsiW_t:psi^(j)(U_rx)=U_omega(x,r)psi^(j)(x)}.

\medskip

For $f\in L^\infty(\Gbf_2/\Gamma_2,\mu_2)$, by \eqref{eq:BehaviorOfpsiW_t:DefOfpsi^(j)}, the definitions of $\mu_\tau$, $\mu_{\hat{\tau}}$, the facts that $\psi_*\mu_\tau=\mu_2$, $\mu_{\hat{\tau}}$ is a Borel probability measure on $\Gbf_1/\Gamma_1$, and the \textit{Pointwise ergodic theorem}, it follows that the subset of $x\in\Zold{Z:BehaviorOfpsiW_t:ResultOfInduction}$ satisfying
\begin{itemize}
    \item $x$ is generic for $f\circ\hat{\psi}$ w.r.t. the ergodic flow $U^{\hat{\tau}}_t$ on $(\Gbf_1/\Gamma_1,\mu_{\hat{\tau}})$,
    \item and $\hat{\psi}(x)=U^{(2)}_{\hat{b}(x)}\psi(x)$ is generic for $f$ w.r.t. the ergodic flow $U^{(2)}_t$ on $(\Gbf_2/\Gamma_2,\mu_2)$,
\end{itemize}
has full $\mu_1$-measure. Hence, with some $x$ in this subset, since \eqref{eq:BehaviorOfpsiW_t:psi^(j)(U_rx)=U_omega(x,r)psi^(j)(x)} implies $\hat{\psi}(U^{\hat{\tau}}_sx)=U^{(2)}_s\hat{\psi}(x)$ for all $s\in\R$, we deduce that
\begin{equation*}
    \begin{aligned}
        \int_{\Gbf_1/\Gamma_1}f\circ\hat{\psi}d\mu_{\hat{\tau}}&=\lim_{r\to\infty}\frac{1}{r}\int_0^rf\left(\hat{\psi}(U^{\hat{\tau}}_sx)\right)ds\\
        &=\lim_{r\to\infty}\frac{1}{r}\int_0^rf\left(U^{(2)}_s\hat{\psi}(x)\right)ds=\int_{\Gbf_2/\Gamma_2}fd\mu_2.
    \end{aligned}
\end{equation*}
Thus, we obtain \eqref{eq:BehaviorOfpsiW_t:psi^(j)PreservesCorrespondingMeasures}.

\medskip

Now, we see that $\hat{\psi}:\Gbf_1/\Gamma_1\to\Gbf_2/\Gamma_2$ is a Kakutani factor map between $(\Gbf_i/\Gamma_i,\mu_i,U^{(i)}_t)$, $i=1,2$, with the corresponding time change $\hat{\tau}$ satisfies all the desired conditions in \textit{Proposition \ref{Prop:BehaviorOfpsiW_t}}, except for the fact that $\hat{\tau}$ may not be $10^{-20}\eold{epsilon:Preliminaries:CoveringLemma}$-good, in particular, for $x\in\Zold{Z:BehaviorOfpsiW_t:ResultOfInduction}$, $[r\mapsto\hat{\tau}(U^{(1)}_rx)]$ may not be $C^\infty$. We address this by another convolution as follows.

Choose $\chi\in C^\infty_c(\R)$ such that $\chi\geq0$ and $\int_\R\chi(t)dt=1$. By the formula of $\hat{\tau}$ in \eqref{eq:BehaviorOfpsiW_t:DefOftau^(j)}, we can define $\hat{\tau}$ on $Z=\{hx:x\in\Zold{Z:BehaviorOfpsiW_t:ResultOfInduction},h\in\Wbf\}$ with the same formula. Recall that $\hat{\psi}$ is also defined on $Z$ (cf. \eqref{eq:BehaviorOfpsiW_t:DefOfpsi^(j)}), so for $x\in Z$, we can define
\begin{equation}\label{eq:BehaviorOfpsiW_t:DefinitionsInSmoothingStep}
\begin{gathered}
     \Tilde{\tau}(x)=\int_\R\chi(t)\hat{\tau}(U^{(1)}_tx)dt,\\
     \Tilde{k}(x)=\int_\R\chi(t)\hat{\omega}(x,t)dt,\qquad\Tilde{\psi}(x)=U^{(2)}_{\Tilde{k}(x)}\hat{\psi}(x).
\end{gathered}
\end{equation}

From the convolution construction of $\Tilde{\tau}$, (\ref{eq:BehaviorOfpsiW_t:|tau^(j)-1|<epsilon}), \textit{Fubini's theorem}, \textit{Change of variables}, and \textit{Dominated convergence theorem}, we deduce that 
\begin{equation*}
    \begin{gathered}
        \norm{\Tilde{\tau}-1}_{\infty}<10^{-20}\eold{epsilon:Preliminaries:CoveringLemma},\qquad\int_{\Gbf_1/\Gamma_1}\Tilde{\tau}(x)d\mu_1(x)=1,\\
        [r\mapsto\Tilde{\tau}(U^{(1)}_rx)]\in C^{\infty}(\R).
    \end{gathered}
\end{equation*}
Hence, $\Tilde{\tau}$ is a $10^{-20}\eold{epsilon:Preliminaries:CoveringLemma}$-good time change on $(\Gbf_1/\Gamma_1,\mu_1)$. Moreover, by the fact $\hat{\omega}$ is a cocycle with respect to $U^{(1)}_t$ (cf. \eqref{eq:BehaviorOfpsiW_t:DefOfomega^(j)}), \textit{Fubini's theorem}, and \eqref{eq:BehaviorOfpsiW_t:DefinitionsInSmoothingStep}, it follows that for $x\in\Zold{Z:BehaviorOfpsiW_t:ResultOfInduction}$,
\begin{align*}
    &\int_0^r\Tilde{\tau}(U^{(1)}_sx)ds=\int_0^r\int_\R\chi(t)\hat{\tau}(U^{(1)}_tU^{(1)}_sx)dt\,ds=\int_\R\chi(t)\int_0^r\hat{\tau}(U^{(1)}_{t+s}x)ds\,dt\\
    &=\int_\R\chi(t)\left(\hat{\omega}(U^{(1)}_rx,t)+\hat{\omega}(x,r)-\hat{\omega}(x,t)\right)dt=\Tilde{k}(U^{(1)}_rx)+\hat{\omega}(x,r)-\Tilde{k}(x).
\end{align*}
And since \eqref{eq:BehaviorOfpsiW_t:DefinitionsInSmoothingStep} and \eqref{eq:BehaviorOfpsiW_t:psi^(j)(U_rx)=U_omega(x,r)psi^(j)(x)} imply that for $x\in\Zold{Z:BehaviorOfpsiW_t:ResultOfInduction}$, $r\in\R$,
\[
\Tilde{\psi}(U^{(1)}_rx)=U^{(2)}_{\Tilde{k}(U^{(1)}_rx)+\hat{\omega}(x,r)-\Tilde{k}(x)}\Tilde{\psi}(x),
\]
we deduce that $\Tilde{\psi}:\Gbf_1/\Gamma_1\to\Gbf_2/\Gamma_2$ is a Kakutani factor map between $(\Gbf_i/\Gamma_i,\mu_i,U^{(i)}_t)$, $i=1,2$, with a $10^{-20}\eold{epsilon:Preliminaries:CoveringLemma}$-good time change $\Tilde{\tau}:\Zold{Z:BehaviorOfpsiW_t:ResultOfInduction}\to\R$ (the fact that $\Tilde{\psi}_*\mu_{\Tilde{\tau}}=\mu_2$ follows similarly as in the proof of \eqref{eq:BehaviorOfpsiW_t:psi^(j)PreservesCorrespondingMeasures}).

From \eqref{eq:BehaviorOfpsiW_t:DefOftau^(j)}, we note that our extension of $\hat{\tau}$ on $Z$ is $\Wbf$-invariant, hence $\hat{\omega}$ in \eqref{eq:BehaviorOfpsiW_t:DefOfomega^(j)} is defined on $Z\times\R$ and is also $\Wbf$-invariant (in the first variable). This implies via \eqref{eq:BehaviorOfpsiW_t:DefinitionsInSmoothingStep} that $\Tilde{k}$ is also $\Wbf$-invariant, which together with \eqref{eq:BehaviorOfpsiW_t:psi^(j)(W^j_sx)}, implies that for $x\in\Zold{Z:BehaviorOfpsiW_t:ResultOfInduction}$,
\begin{equation*}
    \begin{gathered}
         \Tilde{\psi}(gx)=\Tilde{\varpi}(\Tilde{g})\Tilde{\psi}(x),\qquad \forall g\in\Wbf;\\
    \Tilde{\psi}(x)=U^{(2)}_{\Tilde{b}(x)}\psi(x);\qquad\Tilde{\psi}(A^{(1)}_1x)=A^{(2)}_1U^{(2)}_{\Tilde{f}(x)}\Tilde{\psi}(x);
    \end{gathered}
\end{equation*}
where $\Tilde{b}=\hat{b}+\Tilde{k}$ and $\Tilde{f}=\hat{f}+e^{-2}\Tilde{k}\circ A^{(1)}_1-\Tilde{k}$ are measurable functions defined on $\Zold{Z:BehaviorOfpsiW_t:ResultOfInduction}$.

Thus, we finish the proof of \textit{Proposition \ref{Prop:BehaviorOfpsiW_t}}.
\end{proof}
\begin{Remark}\label{Remark:BehaviorOfpsiW_t:Assumepsi=Tilde(psi)}
    Notice that in \textit{Proposition~\ref{Prop:BehaviorOfpsiW_t}}, if $\psi:\Gbf_1/\Gamma_1\to\Gbf_2/\Gamma_2$ is a Kakutani equivalence, then so is $\Tilde{\psi}$. One way to see this is by a direct computation that $\Tilde{\psi}^{-1}(y)=U^{(1)}_{b(y)}\psi^{-1}(y)$, where $b(y)\in\R$ satisfies $\int_0^{b(y)}\Tilde{\tau}(U^{(1)}_s\psi^{-1}(y))ds=-\Tilde{b}(\psi^{-1}(y))$, is the inverse map of $\Tilde{\psi}$.  
\end{Remark}

\section{Renormalization}\label{Section:Renormalization}
In this section, we continue with the settings in \textit{Section~\ref{Section:BehaviorOfpsiW_t}}, and assume in addition that the ergodic $\slfrak_2$-induced unipotent flow $U^{(1)}_t$ on $(\Gbf_1/\Gamma_1,\mu_1)$ is nonsplitting. Our goal is to show $\phi(x)=\lim_{n\to\infty}\psi_n(x)$ exist $\mu_1$-a.e. as a point on the $U^{(2)}_t$-orbit of $\psi(x)$ and to show that $\phi$ conjugates the unipotent flows $U^{(1)}_t$ and $U^{(2)}_t$, i.e. $\phi(U^{(1)}_tx)=U^{(2)}_t\phi(x)$. Here, we denote $\psi_n(x)=A^{(2)}_{-n}\psi(A^{(1)}_nx)$ and $n$ is taken in an increasing sequence of natural numbers of large proportion depending on $x$. After stating our goal in detail, we briefly describe the renormalization procedure in \cite{lindenstrauss2025timechangerigidityunipotent}*{Section 11} and explain the technical changes needed for this procedure in our settings.

\medskip

Before stating this section's main result, as mentioned previously, we fix the notations and assumptions \ref{item:BehaviorOfpsis:1}-\ref{item:BehaviorOfpsis:9} in \textit{Section~\ref{Section:BehaviorOfpsiOnClosePointsAndAlongA_t}}, $\Kold{K:BehaviorOfpsi:GoodSetInLemmapsi(Kak)=Kak}$, $\deold{delta:BehaviorOfpsi:BoundOndeltaInTwoSideMatching}$, $\Rold{R:BehaviorOfpsi:BoundOnRInTwoSideMatching}$, $\eold{epsilon:BehaviorOfpsi:distanceBoundOnG_(d1)/Gamma_1}(\cdot)$ in \textit{Proposition~\ref{Proposition:BehaviorOfpsi:psi(Kak)=Kak}} for $\eta=10^{-6}$, \ref{item:BehaviorOfpsiA_t:1} in p.~\pageref{item:BehaviorOfpsiA_t:1}, \ref{item:BehaviorOfpsiW_t:1}-\ref{item:BehaviorOfpsiW_t:9} in \textit{Section~\ref{Section:BehaviorOfpsiW_t}}. And based on \textit{Proposition~\ref{Prop:BehaviorOfpsiW_t}}, we add:
\begin{enumerate}[label=(\ref{Section:Renormalization}:\roman*)]
    \item\label{item:Renormalization:1} There exists a full measure $U^{(1)}_t$ and $A_1^{(1)}$-invariant set $\Zold{Z:BehaviorOfpsiW_t:ResultOfInduction}\subseteq\Gbf_1/\Gamma_1$ such that for $x\in\Zold{Z:BehaviorOfpsiW_t:ResultOfInduction}$,
    \[
    \textstyle{\psi\left(gx\right)=\Tilde{\varpi}(\Tilde{g})\psi(x),\qquad \forall g\in\Wbf=\Sbf_{1,\operatorname{c}}\Qbf_1,\quad\Tilde{g}\in\pi_\Wbf^{-1}(g),}
    \]
    where $\Tilde{\varpi}:\widetilde{\Wbf}\to\Zbf_2=C_{\Gbf_2}(\Hbf_2)$ is a continuous homomorphism such that
    \[
    \Tilde{\varpi}(\Tilde{\gamma})\in\Gamma_2\cap Z(\Gbf_2),\qquad\forall\Tilde{\gamma}\in\Gamma_{\widetilde{\Wbf}}=\pi_\Wbf^{-1}(\Gamma_{1,\nc}).
    \]
\end{enumerate}

In the above settings, the main result of this section is as follows:
\begin{Theorem}\label{Theorem:Renormalization:phi=lim_(n)psi_n}

    For any $\epsilon\in(0,10^{-6})$, there exist $\Omega\subseteq \mathbf{G}_{1}/\Gamma_1$ a full measure set and measurable functions $n_i:\Omega\to\N,i\in\N,$ such that
    \begin{enumerate}[label=(\alph*)]
        \item\label{item:phi=lim_(n)psi_n:1} For $x\in\Omega$, $\{n_i(x)\}_{i\in\N}$ is increasing and
        \[
        \varliminf_{n\to\infty}\frac{\absolute{\{n_i(x)\}_{i\in\N}\cap[1,n]}}{n}>1-\epsilon.
        \]
        \item\label{item:phi=lim_(n)psi_n:2} For $\displaystyle{x\in\Omega,\phi(x)=\lim_{i\to\infty}\psi_{n_i(x)}(x)}$ exists and is on the $U^{(2)}_t$-orbit of $\psi(x).$
        \item\label{item:phi=lim_(n)psi_n:3}$\phi:\mathbf{G}_{1}/\Gamma_1\to \mathbf{G}_{2}/\Gamma_2$ is measurable and for $x\in\Omega,t\in\Q,$
        $$\phi(U^{(1)}_tx)=U^{(2)}_t\phi(x).$$
    \end{enumerate}
    
\end{Theorem}

\subsection{Set and functions satisfying \ref{item:phi=lim_(n)psi_n:1} and \ref{item:phi=lim_(n)psi_n:2}}\label{subSection:Renormalization:Step1}
We begin the proof of \textit{Theorem~\ref{Theorem:Renormalization:phi=lim_(n)psi_n}} with showing the existence of a set $\Omega'\subseteq \mathbf{G}_{1}/\Gamma_1$ having large measure and measurable functions $\tilde{n}_i:\Omega'\to\N,i\in\N,$ such that conditions \ref{item:phi=lim_(n)psi_n:1} and \ref{item:phi=lim_(n)psi_n:2} are satisfied. This is \cite{lindenstrauss2025timechangerigidityunipotent}*{Proposition 11.5}, and the precise statement is as follows:
\begin{Proposition}\label{Proposition:Renormalization:SetAndFunctionsSatisfyingConditions1,2}

    For any $\epsilon\in(0,10^{-6}),$ there exist $\Omega'\subset \mathbf{G}_{1}/\Gamma_1$ a compact set with $\mu_1(\Omega')>1-\epsilon$ and measurable functions $\tilde{n}_i:\Omega'\to\N,$ $i\in\N,$ such that for $x\in\Omega'$,
    \begin{enumerate}[label=(\~{\alph*})]
        \item\label{item:SetAndFunctionsSatisfyingConditions1,2:1} $\{\tilde{n}_i(x)\}_{i\in\N}$ is increasing and $\displaystyle{\varliminf_{n\to\infty}\frac{\absolute{\{\tilde{n}_i(x)\}_{i\in\N}\cap[1,n]}}{n}>1-\epsilon}$.
        \item\label{item:SetAndFunctionsSatisfyingConditions1,2:2} $\displaystyle{\tilde{\phi}(x)=\lim_{i\to\infty}\psi_{\tilde{n}_i(x)}(x)}$ exists uniformly and is on the $U^{(2)}_t$-orbit of $\psi(x).$
    \end{enumerate}
    
\end{Proposition}
Let us go over the main idea of the proof, outlined in \cite{lindenstrauss2025timechangerigidityunipotent}*{p.~87-90}.

For $\mu_1$-a.e. $x\in\Gbf_1/\Gamma_1$, \ref{item:BehaviorOfpsiW_t:1} gives us
\[
\psi_{n+1}(x)=U^{(2)}_{f(A^{(1)}_nx)e^{-2n}}\psi_n(x),\qquad \forall n\in\N,
\]
where $f:\Gbf_1/\Gamma_1\to\R$ is measurable. Since we have no other condition on $f$, we can only control $\absolute{f(\cdot)}$ on a large subset $K_1$ of $\Gbf_1/\Gamma_1$ via \textit{Lusin's theorem}. But for $\mu_1$-a.e. $x\in\Gbf_1/\Gamma_1$, the \textit{Pointwise ergodic theorem} only allows a subsequence $\{\tilde{n}_i(x)\}_i\subset\N$, which is of proportion equal to $\mu_1(K_1)$, to satisfy $A^{(1)}_{\tilde{n}_i(x)}x\in K_1$. Hence, it does not seem possible for us to ``compare" $\psi_{\tilde{n}_{i+1}(x)}(x)$ and $\psi_{\tilde{n}_i(x)}(x)$ by just applying \ref{item:BehaviorOfpsiW_t:1}, since it only gives us
\begin{equation}\label{eq:Renormalization:MainIdea:psi_n2(x)=U_tpsi_n1(x)}
    \psi_{\tilde{n}_{i+1}(x)}(x)=U^{(2)}_{t_i(x)}\psi_{\tilde{n}_i(x)}(x),\qquad \textstyle{t_i(x)=\sum_{n=\tilde{n}_i(x)}^{\tilde{n}_{i+1}(x)-1}}f(A^{(1)}_nx)e^{-2n},
\end{equation}
and we do not have control over $\absolute{f(A^{(1)}_nx)}$ for every $\tilde{n}_i(x)<n<\tilde{n}_{i+1}(x)$.

To address the above issue, we choose a small $h_n\in \mathbf{H}_1$ such that
\begin{equation}\label{eq:Renormalization:MainIdea:A_nh_nxinK_1}
    A^{(1)}_nh_nx\in K_1,
\end{equation}
then apply \ref{item:BehaviorOfpsiA_t:1} to deduce that
\[
\psi(A^{(1)}_{n+1}h_nx)=A_{1}^{(2)}U^{(2)}_{f(A^{(1)}_nh_nx)}\psi(A^{(1)}_nh_nx),
\]
where $f(A^{(1)}_nh_nx)$ is uniformly bounded. To connect these together and obtain a relation between $\psi_{\tilde{n}_i(x)}(x)$ and $\psi_{\tilde{n}_{i+1}(x)}(x)$, we need to estimate the difference between $\psi(A^{(1)}_nh_nx)$ and $\psi(A^{(1)}_nh_{n-1}x)$. Using \textit{Proposition~\ref{Proposition:BehaviorOfpsi:psi(Kak)=Kak}}, we establish this for $x\in\Gbf_1/\Gamma_1$ and large $n$ such that
\begin{equation}\label{eq:Renormalization:MainIdea:A_nh_nx,A_nh_(n-1)xinK_2}
    A^{(1)}_nh_nx,\;A^{(1)}_nh_{n-1}x\in K_2.
\end{equation}
Then, by \textit{Lemma~\ref{Lemma:Preliminaries:a_(-p)ga_pForgInKak(R,delta)}} and an induction argument, we show 
\[
\psi_{\tilde{n}_{i+1}(x)}(x)=g^{\hfrak}_ig^{\zfrak}_ig^{\tr}_i\psi_{\tilde{n}_i(x)}(x),
\]
where for some small $\delta>0$,
\[
g^\hfrak\in B^{\Hbf_2}_{10(\tilde{n}_{i+1}(x)-\tilde{n}_i(x)+1)\delta},\quad g^\zfrak\in B^{\Zbf_2}_{2(\tilde{n}_{i+1}(x)-\tilde{n}_i(x)+1)\delta},\quad g^\tr\in B^{\Gbf_2}_{e^{-0.9\tilde{n}_i(x)}}.
\]

From here, by a construction based on a \textit{Measure estimate lemma} (cf. \cite{lindenstrauss2025timechangerigidityunipotent}*{Corollary 6.7}), we show that for $x$ in a smaller ``good'' set $K_3$, if $\absolute{t_i(x)}$ in \eqref{eq:Renormalization:MainIdea:psi_n2(x)=U_tpsi_n1(x)} is bounded exponentially w.r.t. $\tilde{n}_i(x)$, then $\absolute{t_i(x)}$ is bounded linearly w.r.t. $\tilde{n}_{i+1}(x)-\tilde{n}_i(x)$.

With the previous estimate, we further show that by shrinking $K_3$ some more, in \eqref{eq:Renormalization:MainIdea:psi_n2(x)=U_tpsi_n1(x)}, $\absolute{t_i(x)}\leq e^{-\tilde{n}_i(x)}$. To do so, we assume for contradiction that $\absolute{t_i(x)}>e^{-\tilde{n}_i(x)}$. This let us find a suitable whole number $k$ such that $\absolute{t_i(x)}e^{-2k}$ is bounded above and below exponentially w.r.t. $\tilde{n}_i(x)$. Applying the previous step provides the exponent has a smaller upper bound for the coefficient of $\tilde{n}_i(x)$ than the previous bound, which is a contradiction. Hence, $\absolute{t_i(x)}\leq e^{-\tilde{n}_i(x)}$, and we finish the proof of \textit{Proposition~\ref{Proposition:Renormalization:SetAndFunctionsSatisfyingConditions1,2}}.

\medskip

For our settings, in the above procedure, we need some technical adjustments for the idea of choosing a small $h_n\in\Hbf_1$ such that, for large $n$, \eqref{eq:Renormalization:MainIdea:A_nh_nxinK_1} and \eqref{eq:Renormalization:MainIdea:A_nh_nx,A_nh_(n-1)xinK_2} hold. In \cite{lindenstrauss2025timechangerigidityunipotent}, the authors carry this out in p.~95, based on item (ii) - p.~93. In it, the sets $\mathcal{K}_{16}$ and $\mathcal{K}_{13}$ correspond to the sets $K_1$ and $K_2$ we describe above, respectively, and the authors use an \textit{$\SL_2(\R)$-ergodic theorem}, namely \cite{lindenstrauss2025timechangerigidityunipotent}*{Corollary 3.15}, on $\mathbf{1}_{\mathcal{K}_{13}\cap\mathcal{K}_{16}}$ to guarantee that we can choose our desired $h_n\in\Hbf_1$. In our settings, we have the corresponding \textit{Proposition~\ref{Prop:SpectralGapCondition:SL_2-ergodicTheorem}}, which (only) applies to $\Sbf_{1,\operatorname{c}}\Qbf_1$-invariant functions, i.e. $\Wbf$-invariant functions. Hence, we need to construct $K_1$ and $K_2$ to be $\Wbf$-invariant. This is where \ref{item:Renormalization:1}, i.e. \textit{Proposition~\ref{Prop:BehaviorOfpsiW_t}}, plays an important role in making sure that if $x$ is ``good", then $gx$ is also ``good" for all $g\in\Wbf$. 

\medskip

In the remainder of this subsection, with $\epsilon\in(0,10^{-6})$ fixed in \textit{Proposition~\ref{Proposition:Renormalization:SetAndFunctionsSatisfyingConditions1,2}}, we construct these two sets $K_1$ and $K_2$. We remark that we do not assume these sets to be contained in $\Zold{Z:BehaviorOfpsi:GoodSetWhereTheConjugationRelationHold}$ (cf. last assumption in \ref{item:BehaviorOfpsis:9}). Once we have $K_1$ and $K_2$, by substituting them and \textit{Lemma~\ref{Lemma:Renormalization:TechnicalLemma}} below in places of $\mathcal{K}_{16}$, $\mathcal{K}_{13}$, and \textit{Lemma~11.8} in \cite{lindenstrauss2025timechangerigidityunipotent}*{Subsection 11.3}, the arguments there work in our settings verbatim, giving us \textit{Proposition~\ref{Proposition:Renormalization:SetAndFunctionsSatisfyingConditions1,2}}.

\medskip

We first construct $K_1$. By \ref{item:BehaviorOfpsiW_t:1} and \textit{Lusin's theorem}, there exist a compact set $K'_1\subset\Zold{Z:BehaviorOfpsiA_t:GoodSetForpsi(A_px)=U_()exp(p(a+X_0))psi(x)}\cap\Zold{Z:BehaviorOfpsiUnderA_t:GenericPointsOfGoodSetInLemmapsi(Kak)=Kak}\cap \Zold{Z:BehaviorOfpsiW_t:ResultOfInduction}$ with $\mu_1(K'_1)>1-10^{-100}\epsilon$ such that $f|_{K'_1}$ is continuous and for some $N_1\geq1$, 
\begin{equation}\label{eq:Renormalization:MainIdea:K_1'}
    \psi(A^{(1)}_1x)=A^{(2)}_1U^{(2)}_{f(x)}\psi(x),\quad\absolute{f(x)}\leq N_1,\qquad\forall x\in K'_1.
\end{equation}
Then the compact $\Wbf$-invariant set
\[
K_1=\Wbf\cdot K_1'\subseteq\Gbf_1/\Gamma_1
\]
satisfies 
\[
\mu_1(K_1)>1-10^{-100}\epsilon.
\]
Given $x_0\in K'_1$ and $g_0\in\Wbf$, it follows from the facts that $x_0,A_1^{(1)}x_0\in\Zold{Z:BehaviorOfpsiW_t:ResultOfInduction}$,  $\Wbf$ commutes with $A_1^{(1)}$, \ref{item:Renormalization:1}, and  \eqref{eq:Renormalization:MainIdea:K_1'} that
\begin{equation}\label{eq:Renormalization:MainIdea:K_1'b}
    \psi(A_1^{(1)}g_0x_0)=\psi(g_0A_1^{(1)}x_0)=\tilde{\varpi}(\tilde{g}_0)\psi(A_1^{(1)}x_0)=A_1^{(2)}U_{f(x_0)}^{(2)}\psi(g_0x_0).
\end{equation}
If $g_1x_1=g_2x_2$ holds for some $x_1,x_2\in K'_1$ and $g_1,g_2\in\Wbf$, then  \ref{item:BehaviorOfpsiA_t:1}, \eqref{eq:Renormalization:MainIdea:K_1'b},  $\psi(g_1x_1)=\tilde{\varpi}(\tilde{g}_1)\psi(x_1)$, and $\tilde{\varpi}(\tilde{g}_1)\in C_{\Gbf_2}(\Hbf_2)$ imply that
\begin{equation}\label{eq:sameWFoliation}
    f(x_1)=f(x_2).
\end{equation}
Hence given $y=gx\in K_1$ for some $g\in\Wbf$ and $x\in K'_1$, we extend the definition of $f$ as
\[
f(y)=f(x),
\]
which is well defined due to \eqref{eq:sameWFoliation}. It is clear that this definition is compatible with the original definition of $f$. Notice that $\Wbf/\Gamma_{1,nc}$ is compact (cf. \ref{item:BehaviorOfpsiW_t:8}), this together with the fact that $f|_{K'_1}$ is continuous gives the measurability of the extension of $f$. Hence
\[
\psi(A^{(1)}_1x)=A^{(2)}_1U^{(2)}_{f(x)}\psi(x),\quad\absolute{f(x)}\leq N_1,\qquad\forall x\in K_1.
\]

Thus, we obtain our first ``good" set $K_1$.

\medskip

To construct $K_2$, we need the following technical lemma, corresponding to \cite{lindenstrauss2025timechangerigidityunipotent}*{Lemma 11.8}:
\begin{Lemma}\label{Lemma:Renormalization:TechnicalLemma}

    There exist $\denew\label{delta:Renormalization:TechnicalLemma}>0$ and $\enew\label{epsilon:Renormalization:TechnicalLemma}(\delta)\in(0,\frac{\delta}{100})$ for $\delta\in(0,\deold{delta:Renormalization:TechnicalLemma}]$, satisfying the following property. For $\epsilon'\in(0,10^{-6}),$ there exist $\Rnew\label{R:Renormalization:TechnicalLemma:A_tDeviation}(\epsilon'),$ $\Rnew\label{R:Renormalization:TechnicalLemma:U_tDeviation}(\epsilon')\geq1,$ and a compact $\Wbf$-invariant set $\Knew\label{K:Renormalization:TechnicalLemma:GoodSet}(\epsilon')\subset \mathbf{G}_{1}/\Gamma_1$ with $\mu_1(\Kold{K:Renormalization:TechnicalLemma:GoodSet}(\epsilon'))>1-\epsilon'$ such that if $x\in \mathbf{G}_{1}/\Gamma_1,h_1,h_2\in B^{\mathbf{H}_1}_{\eold{epsilon:Renormalization:TechnicalLemma}(\delta)}$, and $t\geq\Rold{R:Renormalization:TechnicalLemma:A_tDeviation}(\epsilon')$ satisfy
    \[
    A^{(1)}_th_1x, \;A^{(1)}_th_2x\in\Kold{K:Renormalization:TechnicalLemma:GoodSet}(\epsilon'),
    \]
    then there exist $s,c_{\bu},c_{\ba},c_{\bou}\in\R$, $g^{\zfrak}\in\Zbf_2$, and $g^{\tr}\in\Gbf_2$ such that
    \[
    \psi(A^{(1)}_th_2x)=U^{(2)}_{c_{\bu}}\Uh^{(2)}_{c_{\bou}}A^{(2)}_{c_{\ba}}g^{\zfrak}g^{\tr}U^{(2)}_s\psi(A^{(1)}_th_1x)
    \]
    holds with the following coefficients' bounds
    \begin{equation*}
        \begin{gathered}
        \absolute{s}\leq6(\Rold{R:Renormalization:TechnicalLemma:U_tDeviation}(\epsilon')+\eold{epsilon:Renormalization:TechnicalLemma}(\delta)e^{2t}),\quad\absolute{c_{\bu}}\leq\delta,\quad\absolute{c_{\ba}}\leq\delta,\quad\absolute{c_{\bou}}\leq\delta e^{-2t},\quad g^{\zfrak}\in B^{\Zbf_2}_{2\delta},\\
        d_{\Gbf_2}(U^{(2)}_rg^{\tr}U^{(2)}_{-r},e_{\Gbf_2})\leq2\delta,\qquad \forall r\in[-e^{2t},e^{2t}].
        \end{gathered}
    \end{equation*} 
\end{Lemma}

\medskip

With this lemma, we construct $K_2$ as follows. Let $\deold{delta:Renormalization:TechnicalLemma},$ $\eold{epsilon:Renormalization:TechnicalLemma}(\cdot)$ $\Kold{K:Renormalization:TechnicalLemma:GoodSet}(\cdot)$ be given in \textit{Lemma~\ref{Lemma:Renormalization:TechnicalLemma}}, and choose
\[
K_2=\Kold{K:Renormalization:TechnicalLemma:GoodSet}(10^{-100}\epsilon),
\]
then $K_2\subset\Gbf_1/\Gamma_1$ is compact, $\Wbf$-invariant, and $\mu_1(K_2)>1-10^{-100}\epsilon$.

Fix some $\delta\in(0,\frac{1}{2}\deold{delta:Renormalization:TechnicalLemma})$ small enough (depending only on $\Gbf_2$) and choose $\eta=\eta(\delta)>0$ small enough such that
\[
B^{\Hbf_1,\norm{\cdot}}_\eta\subseteq B^{\Hbf_1}_{\eold{epsilon:Renormalization:TechnicalLemma}(\delta)}.
\]
By construction of $K_1$, $K_2$, their intersection $K_2\cap K_1$ is a compact $\Wbf$-invariant set in $\Gbf_1/\Gamma_1$. Hence, we can apply \textit{Proposition~\ref{Prop:SpectralGapCondition:SL_2-ergodicTheorem}} to $\mathbf{1}_{K_2\cap K_1}$ and obtain
\begin{align*}
    \lim_{n\to\infty}\frac{1}{m_{\mathbf{H}_1}\left(B^{\mathbf{H}_1,\norm{\cdot}}_{\eta}\right)}\int_{B^{\mathbf{H}_1,\norm{\cdot}}_{\eta}}\mathbf{1}_{K_2\cap K_1}(A^{(1)}_nhx)dm_{\mathbf{H}_1}(h)\\
    =\mu(K_2\cap K_1)>1-10^{-99}\epsilon,\qquad\mu_1\text{-a.e. }x\in \mathbf{G}_{1}/\Gamma_1.
\end{align*}

Thus, we obtain our second ``good" set $K_2$.

\medskip

Now, we return to the proof of \textit{Lemma~\ref{Lemma:Renormalization:TechnicalLemma}}. It is similar to \cite{lindenstrauss2025timechangerigidityunipotent}*{Lemma~11.8}, with adjustments for the $\Wbf$-invariance of ``good" set.

\begin{proof}[Proof of \textit{Lemma~\ref{Lemma:Renormalization:TechnicalLemma}}]
We fix some constants and sets as follows:
\begin{enumerate}[label=(\roman*)]
    \item\label{item:Renormalization:TechnicalLemma:1} Recall $\deold{delta:BehaviorOfpsi:BoundOndeltaInTwoSideMatching},\Rold{R:BehaviorOfpsi:BoundOnRInTwoSideMatching},\Kold{K:BehaviorOfpsi:GoodSetInLemmapsi(Kak)=Kak},\eold{epsilon:BehaviorOfpsi:distanceBoundOnG_(d1)/Gamma_1}(\cdot)$ given in \textit{Proposition~\ref{Proposition:BehaviorOfpsi:psi(Kak)=Kak}} for $\eta=10^{-6}$ and $\Zold{Z:BehaviorOfpsiW_t:ResultOfInduction}$, $\Tilde{\varpi}$ given in \ref{item:Renormalization:1}.
    \item\label{item:Renormalization:TechnicalLemma:2} For $\epsilon\in(0,10^{-6})$, by the \textit{Pointwise ergodic theorem}, \textit{Egorov's Theorem}, and $\mu_1$'s regularity, we can choose a compact set $K_3(\epsilon)\subset \mathbf{G}_{1}/\Gamma_1$ with $\mu_1(K_3(\epsilon))>1-\epsilon^2$ and $\Rold{R:Renormalization:TechnicalLemma:U_tDeviation}(\epsilon)\geq1$ such that for $x\in K_3(\epsilon),R\geq\Rold{R:Renormalization:TechnicalLemma:U_tDeviation}(\epsilon),$
    \begin{equation*}
        \absolute{\left\{s\in[-R,R]:U^{(1)}_sx\in\Kold{K:BehaviorOfpsi:GoodSetInLemmapsi(Kak)=Kak}\right\}}>(1-6\cdot10^{-6})\cdot2R.
    \end{equation*}
    \item\label{item:Renormalization:TechnicalLemma:3} Since $\mu_1\left(K_3(\epsilon)\cap\Zold{Z:BehaviorOfpsiW_t:ResultOfInduction}\right)>1-\epsilon^2$, by \textit{Fubini's theorem}, we have
    \begin{align*}
        1-\epsilon^2&<\int_{F_\Wbf}\mu_1\left(g^{-1}\cdot\left(K_3(\epsilon)\cap\Zold{Z:BehaviorOfpsiW_t:ResultOfInduction}\right)\right)dm_\Wbf(g)\\
        &=\int_{(\Gbf_1/\Gamma_1)\times F_\Wbf}\mathbbm{1}_{\{(x,g):gx\in K_3(\epsilon)\cap\Zold{Z:BehaviorOfpsiW_t:ResultOfInduction}\}}(x,g)d(\mu_1\times m_\Wbf)(x,g)\\
        &=\int_{\Gbf_1/\Gamma_1}m_\Wbf\left(\left\{g\in F_\Wbf:gx\in K_3(\epsilon)\cap\Zold{Z:BehaviorOfpsiW_t:ResultOfInduction}\right\}\right)d\mu_1(x).
    \end{align*}
     
    Hence, there is a compact set $K_4(\epsilon)\subset\Gbf_1/\Gamma_1$ with $\mu_1(K_4(\epsilon))>1-\epsilon$ such that for $x\in K_4(\epsilon)$,
    \[
    m_\Wbf\left(\left\{g\in F_\Wbf:gx\in K_3(\epsilon)\cap\Zold{Z:BehaviorOfpsiW_t:ResultOfInduction}\right\}\right)>1-\epsilon.
    \]
    \item\label{item:Renormalization:TechnicalLemma:4} Choose $\Kold{K:Renormalization:TechnicalLemma:GoodSet}(\epsilon)=\Wbf\cdot K_4(\epsilon)\subseteq\Gbf_1/\Gamma_1$.
    \item\label{item:Renormalization:TechnicalLemma:5} Recall from \ref{item:BehaviorOfpsiW_t:9} (and \ref{item:BehaviorOfpsiW_t:8}) that $\overline{F_{\widetilde{\Wbf}}}$ is compact. Choose
    \[
    \kappa>10\left(1+\max\left\{\norm{\Ad\left(\Tilde{\varpi}(\Tilde{g})\right)}:\Tilde{g}\in\overline{F_{\widetilde{\Wbf}}}\right\}\right)
    \]
    sufficiently large, where $\norm{\Ad(\cdot)}$ is the operator norm on $\GL(\gfrak_2)$. We also define
    \[
    \deold{delta:Renormalization:TechnicalLemma}=10^{-2}\deold{delta:BehaviorOfpsi:BoundOndeltaInTwoSideMatching}.
    \]
    Then for all $\delta\in(0,\deold{delta:Renormalization:TechnicalLemma}]$, we define
    \[
    \eold{epsilon:Renormalization:TechnicalLemma}(\delta)=10^{-1}\eold{epsilon:BehaviorOfpsi:distanceBoundOnG_(d1)/Gamma_1}(\delta/\kappa)\in(0,10^{-2}\delta).
    \]
    \item\label{item:Renormalization:TechnicalLemma:6} Choose $\Rold{R:Renormalization:TechnicalLemma:A_tDeviation}(\epsilon)=\frac{1}{2}\log\left(\max\{200\Rold{R:Renormalization:TechnicalLemma:U_tDeviation}(\epsilon),\Rold{R:BehaviorOfpsi:BoundOnRInTwoSideMatching}\}\right)\geq1.$
\end{enumerate}

\medskip

By \ref{item:Renormalization:TechnicalLemma:3} and \ref{item:Renormalization:TechnicalLemma:4}, it follows that $\Kold{K:Renormalization:TechnicalLemma:GoodSet}(\epsilon)\subset\Gbf_1/\Gamma_1$ is a compact $\Wbf$-invariant set with $\mu_1(\Kold{K:Renormalization:TechnicalLemma:GoodSet}(\epsilon))>1-\epsilon$.

For $\delta\in(0,\deold{delta:Renormalization:TechnicalLemma}],$ consider $x\in \mathbf{G}_{1}/\Gamma_1,h_1,h_2\in B^{\mathbf{H}_1}_{\eold{epsilon:Renormalization:TechnicalLemma}(\delta)},$ and $t\geq\Rold{R:Renormalization:TechnicalLemma:A_tDeviation}(\epsilon)$ satisfying $A^{(1)}_th_1x,A^{(1)}_th_2x\in\Kold{K:Renormalization:TechnicalLemma:GoodSet}(\epsilon)$, and write $h_2h_1^{-1}=\Uh^{(1)}_{t_{\bou}}A^{(1)}_{t_{\ba}}U^{(1)}_{t_{\bu}}$ for some $\absolute{t_{\bou}},\absolute{t_{\ba}},\absolute{t_{\bu}}<2\eold{epsilon:Renormalization:TechnicalLemma}(\delta)$, we have
\[A^{(1)}_th_2x=\Uh^{(1)}_{t_{\bou}e^{-2t}}A^{(1)}_{t_{\ba}}U^{(1)}_{t_{\bu}e^{2t}}A^{(1)}_{t}h_1x.
\]
Letting $\phi(s)=\frac{se^{t_{\ba}}}{e^{-t_{\ba}}-t_{\bou}e^{-2t+t_{\ba}}s},$ it follows from \textit{Lemma~\ref{Lemma:Preliminaries:SL_2Computation}} that
\begin{equation}\label{eq:Renormalization:TechnicalLemma:SL_2Computation}
    \begin{alignedat}{2}
        &\text{for }\absolute{s}\leq\textstyle{\frac{10^{-2}e^{2t}}{2\eold{epsilon:Renormalization:TechnicalLemma}(\delta)}},\quad&&U^{(1)}_{\phi(s)}A^{(1)}_th_2x=\Uh^{(1)}_{t_{\bou}(s)}A^{(1)}_{t_{\ba}(s)}U^{(1)}_{s+t_{\bu}e^{2t}}A^{(1)}_th_1x,\\
        & &&\absolute{t_{\bou}(s)}<4\eold{epsilon:Renormalization:TechnicalLemma}(\delta)e^{-2t},\absolute{t_{\ba}(s)}<2\eold{epsilon:Renormalization:TechnicalLemma}(\delta)(1+2e^{-2t}\absolute{s});\\
        &\text{for }\absolute{s}\leq0.5e^{2t},&&\absolute{\phi'(s)-1}<5\eold{epsilon:Renormalization:TechnicalLemma}(\delta).
    \end{alignedat}
\end{equation}

Since $A^{(1)}_th_1x, A^{(1)}_th_2x\in\Kold{K:Renormalization:TechnicalLemma:GoodSet}(\epsilon)$, by \ref{item:Renormalization:TechnicalLemma:3} and \ref{item:Renormalization:TechnicalLemma:4}, it follows that there exists $w\in F_\Wbf$ such that 
\[
wA^{(1)}_th_1x, \;wA^{(1)}_th_2x\in K_3(\epsilon)\cap\Zold{Z:BehaviorOfpsiW_t:ResultOfInduction}.
\]

With this, \ref{item:Renormalization:TechnicalLemma:6}, and \ref{item:Renormalization:TechnicalLemma:2}, we can show that there exists $s_0\in\R$ with $\absolute{s_0}\leq2\Rold{R:Renormalization:TechnicalLemma:U_tDeviation}(\epsilon)+\eold{epsilon:Renormalization:TechnicalLemma}(\delta)e^{2t}$ such that
\[
U^{(1)}_{\phi(s_0)}wA^{(1)}_th_2x,\;U^{(1)}_{s_0+t_{\bu}e^{2t}}wA^{(1)}_th_1x\in\Kold{K:BehaviorOfpsi:GoodSetInLemmapsi(Kak)=Kak}.
\]
And as (\ref{eq:Renormalization:TechnicalLemma:SL_2Computation}) and \ref{item:Renormalization:TechnicalLemma:5} imply
\[
U^{(1)}_{\phi(s_0)}wA^{(1)}_th_2x\in\Kak\left(e^{2t},\eold{epsilon:BehaviorOfpsi:distanceBoundOnG_(d1)/Gamma_1}(\textstyle{\frac{\delta}{\kappa}})\right)U^{(1)}_{s_0+t_{\bu}e^{2t}}wA^{(1)}_th_1x,
\]
it follows from \textit{Proposition~\ref{Proposition:BehaviorOfpsi:psi(Kak)=Kak}} that 
\begin{equation}\label{eq:Renormalization:TechnicalLemma:psi(A_th_1x)andpsi(A_th_2x)}
    \begin{aligned}
        \psi(U^{(1)}_{\phi(s_0)}wA^{(1)}_th_2x)=g\psi(U^{(1)}_{s_0+t_{\bu}e^{2t}}wA^{(1)}_th_1x),\quad g\in\Kak\left(1.9e^{2t},\textstyle{\frac{\delta}{\kappa}}\right).
    \end{aligned}
\end{equation}

By (\ref{eq:BehaviorOfpsi:|partial_svarpi_i(x,s)-1|IsSmall}), (\ref{eq:BehaviorOfpsi:ConjugationRelation}), (\ref{eq:Renormalization:TechnicalLemma:SL_2Computation}), \ref{item:Renormalization:1}, and the fact that $\absolute{s_0}\leq2\Rold{R:Renormalization:TechnicalLemma:U_tDeviation}(\epsilon)+\eold{epsilon:Renormalization:TechnicalLemma}(\delta)e^{2t}$, let $\Tilde{w}\in F_{\widetilde{\Wbf}}$ be the unique lift of $w\in F_\Wbf$ (cf. \ref{item:BehaviorOfpsiW_t:8}, \ref{item:BehaviorOfpsiW_t:9}), we have
\[
\psi(U^{(1)}_{\phi(s_0)}wA^{(1)}_th_2x)=U^{(2)}_{r_2}\Tilde{\varpi}(\Tilde{w})\psi(A^{(1)}_th_2x),
\]
with $r_2=\omega(wA^{(1)}_th_2x,\phi(s_0))$ satisfying
\[
\absolute{r_2}\leq(1+10^{-9})\left(1+5\eold{epsilon:Renormalization:TechnicalLemma}(\delta)\right)(2\Rold{R:Renormalization:TechnicalLemma:U_tDeviation}(\epsilon)+\eold{epsilon:Renormalization:TechnicalLemma}(\delta)e^{2t});
\]
and
\[
\psi(U^{(1)}_{s_0+t_{\bu}e^{2t}}wA^{(1)}_th_1x)=U^{(2)}_{r_1}\Tilde{\varpi}(\Tilde{w})\psi(A^{(1)}_th_1x),
\]
with $r_1=\omega(wA^{(1)}_th_1x,s_0+t_{\bu}e^{2t})$ satisfying
\[
\absolute{r_1}\leq(1+10^{-9})(2\Rold{R:Renormalization:TechnicalLemma:U_tDeviation}(\epsilon)+3\eold{epsilon:Renormalization:TechnicalLemma}(\delta)e^{2t}).
\]

Using this, \eqref{eq:Renormalization:TechnicalLemma:psi(A_th_1x)andpsi(A_th_2x)}, and \ref{item:Renormalization:TechnicalLemma:5}, we finish the proof of \textit{Lemma~\ref{Lemma:Renormalization:TechnicalLemma}}.
\end{proof}

With the sets $K_1,K_2$ constructed and \textit{Lemma~\ref{Lemma:Renormalization:TechnicalLemma}} proved, as we explained previously, the arguments in \cite{lindenstrauss2025timechangerigidityunipotent}*{Subsection 11.3}, more specifically, in p.~91-101, apply to our settings and prove \textit{Proposition~\ref{Proposition:Renormalization:SetAndFunctionsSatisfyingConditions1,2}}.

\subsection{Finishing the proof of \textit{Theorem~\ref{Theorem:Renormalization:phi=lim_(n)psi_n}}}\label{subSection:Renormalization:FinishProofOfLimitTheorem}
In this subsection, we conclude the proof of \textit{Theorem~\ref{Theorem:Renormalization:phi=lim_(n)psi_n}}. The arguments we need are almost identical to those appear in \cite{lindenstrauss2025timechangerigidityunipotent}*{Subsection 11.2}. Hence, we refer the reader there for the details, and only mention at certain points what arguments are used.

\medskip

First, we construct $\Omega''\subset\Gbf_1/\Gamma_1$ and measurable functions $\overline{n}_i$, $i\in\N$, such that conditions \ref{item:phi=lim_(n)psi_n:1}, \ref{item:phi=lim_(n)psi_n:2} and a condition similar to \ref{item:phi=lim_(n)psi_n:3} are satisfied. The precise statement is as follows:
\begin{Proposition}\label{Proposition:Renormalization:SetAndFunctionsSatisfyingConditions1,2,3}

    For any $\epsilon\in(0,10^{-6}),$ there exist $\Omega''\subset\mathbf{G}_{1}/\Gamma_1$ compact with $\mu_1(\Omega'')>1-\epsilon$ and measurable functions $\overline{n}_i:\Omega''\to\N,$ $i\in\N,$ such that
    \begin{enumerate}[label=(${\overline{\mbox{\alph*}}}$)]
        \item\label{item:SetAndFunctionsSatisfyingConditions1,2,3:1} For $x\in\Omega''$, $\{\overline{n}_i(x)\}_{i\in\N}$ is increasing and
        \[
        \varliminf_{n\to\infty}\frac{\absolute{\{\overline{n}_i(x)\}_{i\in\N}\cap[1,n]}}{n}>1-\epsilon.
        \]
        \item\label{item:SetAndFunctionsSatisfyingConditions1,2,3:2} For $x\in\Omega'',\overline{\phi}(x)=\lim_{i\to\infty}\psi_{\overline{n}_i(x)}(x)$ exists uniformly and lies on the $U^{(2)}_t$-orbit of $\psi(x).$
        \item\label{item:SetAndFunctionsSatisfyingConditions1,2,3:3} For $t\in\R,\lim_{i\to\infty}\psi_{\overline{n}_i(x)}(U^{(1)}_tx)=U^{(2)}_t\overline{\phi}(x)$ uniformly in $x\in\Omega''.$
    \end{enumerate}
    
\end{Proposition}
We apply \textit{Proposition~\ref{Proposition:Renormalization:SetAndFunctionsSatisfyingConditions1,2}} to obtain $\Omega'^{(j)}$, $\{\Tilde{n}_i^{(j)}\}_{i\in\N}$ corresponding to $10^{-j}\epsilon$, $j\in\N$. From here, we apply \cite{lindenstrauss2025timechangerigidityunipotent}*{Lemma~11.6} to the sequence $\left(\Omega'^{(j)},\{\Tilde{n}_i^{(j)}\}_{i\in\N},10^{-j}\epsilon\right)_{j\in\N}$, and then apply \cite{lindenstrauss2025timechangerigidityunipotent}*{Lemma~11.7} to obtain \textit{Proposition~\ref{Proposition:Renormalization:SetAndFunctionsSatisfyingConditions1,2,3}}.

\medskip

With $\Omega''$, $\overline{n}_i$, $i\in\N$, given in \textit{Proposition~\ref{Proposition:Renormalization:SetAndFunctionsSatisfyingConditions1,2,3}}, we define
\[
\textstyle{\Omega=\bigcup_{q\in\Q}U^{(1)}_q\Omega''}.
\]
Let $\lessdot$ be a well-ordering of $\Q$, we can then define 
\[
\textstyle{\alpha:\Omega\to\Q:x\mapsto\min_{\lessdot}\{q\in\Q:U^{(1)}_qx\in\Omega''\},}
\]
\[
n_i:\Omega\to\N:x\mapsto\overline{n}_i(U^{(1)}_{\alpha(x)}x),\qquad\forall i\in\N.
\]
Now, using the arguments in \cite{lindenstrauss2025timechangerigidityunipotent}*{p.~87}, we obtain \textit{Theorem~\ref{Theorem:Renormalization:phi=lim_(n)psi_n}}.\qed

\section{Proof of \textit{Theorem~\ref{thm:main}}}\label{Section:ProofOfTheorem:MainTheorem}
In this section, we prove \textit{Theorem~\ref{thm:main}}.

\medskip

Let $(\Gbf_i/\Gamma_i,\mu_i,U^{(i)}_t)$ be an ergodic nonsplitting $\slfrak_2$-induced unipotent flow, $i=1,2$, and $\psi_{(3)}:\Gbf_1/\Gamma_1\to\Gbf_2/\Gamma_2$ a Kakutani factor map between them. We apply \textit{Proposition~\ref{Prop:epsilon-goodTimeChange:epsilon-goodTimeChange}}, \textit{Proposition~\ref{Proposition:BehaviorOfpsiA_t}}, \textit{Proposition~\ref{Prop:BehaviorOfpsiW_t}} consecutively to obtain the following Kakutani factor maps $\psi_{(2)},\psi_{(1)},\psi:\Gbf_1/\Gamma_1\to\Gbf_2/\Gamma_2$ between $(\Gbf_i/\Gamma_i,\mu_i,U^{(i)}_t)$, $i=1,2$, whose corresponding time changes are $10^{-20}\eold{epsilon:Preliminaries:CoveringLemma}$-good ($\eold{epsilon:Preliminaries:CoveringLemma}$ given on p.~\pageref{epsilon:Preliminaries:CoveringLemma} for $(\Gbf_2/\Gamma_2,\mu_2,U^{(2)}_t)$), such that
\begin{equation}\label{eq:ProofOfTheorem:epsilon-goodTimeChangeReduction}
    \psi_{(3)}(\cdot)=\psi_{(2)}(U^{(1)}_{\sigma_2(\cdot)}A^{(1)}_{p_0}(\cdot));
\end{equation}
\begin{equation}\label{eq:ProofOfTheorem:A_tDirecttionReduction}
    \psi_{(1)}=g_{\psi_{(2)}}\psi_{(2)};
\end{equation}
\begin{equation}\label{eq:ProofOfTheorem:CompactExtensionFactorReduction}
    \psi(\cdot)=U^{(2)}_{\sigma_1(\cdot)}\psi_{(1)}(\cdot).   
\end{equation}
Here, $p_0\in\R$, 
\[
\textstyle{g_{\psi_{(2)}}\in\exp\left(\R\bu^{(2)}\oplus\bigoplus_{j\notin J_2}\R\bx^{(2)}_{0,j}\right)},\qquad \text{$J_2$ given in \ref{item:BehaviorOfpsis:6}},
\]
and $\sigma_1,\sigma_2:\Gbf_1/\Gamma_1\to\R$ are measurable. Note that the required assumptions of each step are guaranteed by its preceding step. And in the end, $\psi$ also satisfies the required assumptions of \textit{Theorem~\ref{Theorem:Renormalization:phi=lim_(n)psi_n}}, so we can apply it and obtain a full measure set $\Omega\subseteq \mathbf{G}_{1}/\Gamma_1$ and measurable functions $n_i:\Omega\to\N$, $i\in\N$, such that
\begin{enumerate}[label=(\alph*)]
    \item\label{item:ProofOfMainTheorem:1} For $x\in\Omega$, $\{n_i(x)\}_{i\in\N}$ is increasing and
    \[
    \varliminf_{n\to\infty}\frac{\absolute{\{n_i(x):i\in\N\}\cap[1,n]}}{n}\geq1-10^{-9}.
    \]
    \item\label{item:ProofOfMainTheorem:2} For $x\in\Omega,\phi(x)=\lim_{i\to\infty}\psi_{n_i(x)}(x)=\lim_{i\to\infty}A^{(2)}_{-n_i(x)}\psi(A^{(1)}_{n_i(x)}x)$ exists and is on the $U^{(2)}_t$-orbit of $\psi(x).$
    \item\label{item:ProofOfMainTheorem:3} $\phi:\mathbf{G}_{1}/\Gamma_1\to \mathbf{G}_{2}/\Gamma_2$ is measurable and for $x\in\Omega,t\in\Q,$
    \[
    \phi(U^{(1)}_tx)=U^{(2)}_t\phi(x).
    \]
\end{enumerate}

\medskip

From conditions \ref{item:ProofOfMainTheorem:2} and \ref{item:ProofOfMainTheorem:3}, an argument using \textit{Ratner's Measure Classification Theorem} gives us the following:
\begin{Claim}\label{Claim:ProofOfTheorem:phiIsMeasurePreserving}

    $\phi_*\mu_1=\mu_2$, i.e. $\phi:(\Gbf_1/\Gamma_1,\mu_1)\to(\Gbf_2/\Gamma_2,\mu_2)$ is measure preserving.
    
\end{Claim}
\begin{proof}
By condition \ref{item:ProofOfMainTheorem:3}, $\phi_*\mu_1$ is an ergodic $\{U^{(2)}_t:t\in\Q\}$-invariant Borel probability measure on $\Gbf_2/\Gamma_2$. This implies via \textit{Ratner's Measure Classification Theorem} (cf. \cite{MR1135878}*{Theorem 1}) that $\phi_*\mu_1$ is supported on an orbit $G\cdot x_0$ in $\mathbf{G}_{2}/\Gamma_2,$ where $G$ is a closed subgroup of $\mathbf{G}_{2}$ containing $\{U^{(2)}_t:t\in\Q\}.$

Assume that $G\neq\Gbf_2$. Since $\Gbf_2$ is connected, $\dim(G)<\dim(\Gbf_2)$, which implies that $G\cdot x_0$ is a submanifold of $\Gbf_2/\Gamma_2$ of dimension less than $\dim(\mathbf{G}_{2}/\Gamma_2)$. Hence, $\mu_2(G\cdot x_0)=0$.

Recall that $G$ is closed and contains $\{U^{(2)}_t:t\in\Q\}$ and hence
\[
\{U^{(2)}_t:t\in\R\}\subseteq G.
\]
This together with 
\[
\mu_1(\phi^{-1}(G\cdot x_0))=\phi_*\mu_1(G\cdot x_0)=1
\]
and condition \ref{item:ProofOfMainTheorem:2} implies that 
\[
\psi(x)\in G\cdot x_0,\qquad \mu_1\text{-a.e. }x\in \mathbf{G}_{1}/\Gamma_1,
\]
or $\mu_1(\psi^{-1}(G\cdot x_0))=1$, and so $\mu_{\tau}(\psi^{-1}(G\cdot x_0))=1.$ As $\psi_*\mu_\tau=\mu_2$, we get $\mu_2(G\cdot x_0)=1,$ contradicting $\mu_2(G\cdot x_0)=0$.

Thus $G=\mathbf{G}_{2},$ and we obtain $\phi_*\mu_1=\mu_2$.
\end{proof}
We recall from \textit{Section~\ref{Section:Introduction}} the notations 
\[
\Gbf_i'=\Gbf_i/\Core_{\Gbf_i}(\Gamma_i),\quad\Gamma_i'=\Gamma_i/\Core_{\Gbf_i}(\Gamma_i),\quad\vartheta_i(g\Gamma_i)=(g\Core_{\Gbf_i}(\Gamma_i))\Gamma_i',
\]
and note that $\vartheta_i:\Gbf_i/\Gamma_i\to\Gbf_i'/\Gamma_i'$ is an isomorphism of homogeneous spaces, preserving the corresponding Haar probability measures, i.e. $(\vartheta_i)_*\mu_i=\mu_i'$, where $\mu_i'$ is the Haar probability measure on $\Gbf_i'/\Gamma_i'$, $i=1,2$. By \textit{Claim \ref{Claim:ProofOfTheorem:phiIsMeasurePreserving}} and condition \ref{item:ProofOfMainTheorem:3}, \cite{MR910358}*{Theorem 2.1'} (or \cite{MR2611846}*{Main theorem 1.3}) implies that $\phi$ is affine: there exists $(g_0,\alpha)\in\Gbf_2'\times\Sur(\Gbf_1',\Gbf_2')$ such that
\begin{equation}\label{eq:Claim:ProofOfTheorem:FormOfphi}
\begin{gathered}
        \alpha(\Gamma_1')\subseteq g_0\Gamma_2'g_0^{-1};\\
        \phi'(g\Gamma_1')=(\vartheta_2\circ\phi\circ\vartheta_1^{-1})(g\Gamma_1')=\alpha(g)g_0\Gamma_2',\quad \mu_1'\text{-a.e. }g\Gamma_1'\in\Gbf_1'/\Gamma_1'.
\end{gathered}
\end{equation}
\begin{Claim}\label{Claim:ProofOfTheorem:phi(U_t,A_t,Uh_t)}

    Let $U'^{(i)}_t=U^{(i)}_t\Core_{\Gbf_i}(\Gamma_i)$, $i=1,2$. Then
    \[
    \alpha(U'^{(1)}_t)=U'^{(2)}_t,\qquad\forall t\in\R.
    \]
    
\end{Claim}
This claim directly follows from \eqref{eq:Claim:ProofOfTheorem:FormOfphi} and condition \ref{item:ProofOfMainTheorem:3}. Combine it, (\ref{eq:ProofOfTheorem:epsilon-goodTimeChangeReduction})-(\ref{eq:Claim:ProofOfTheorem:FormOfphi}), and the fact that $\phi(\cdot)=U^{(2)}_{\sigma_0(\cdot)}\psi(\cdot)$ for some $\sigma_0:\Gbf_1/\Gamma_1\to\R$ measurable (cf. condition \ref{item:ProofOfMainTheorem:2}) together, we deduce that
\begin{equation}\label{eq:ProofOfTheorem:affinePair}
    \begin{gathered}
        \left(\vartheta_2\circ\left(U^{(2)}_{\sigma(\cdot)}\psi_{(3)}(\cdot)\right)\circ\vartheta_1^{-1}\right)(g\Gamma_1')=\alpha_1(g)g_1\Gamma_2',\qquad \mu_1'\text{-a.e. }g\in\Gbf_1',\\
        g_1=\left(g_{\psi_{(2)}}^{-1}\Core_{\Gbf_2}(\Gamma_2)\right)\alpha(A^{(1)}_{p_0}\Core_{\Gbf_1}(\Gamma_1))g_0\in\Gbf_2',
    \end{gathered}
\end{equation}
where $\sigma:\Gbf_1/\Gamma_1\to\R$ and $\alpha_1\in\Sur(\Gbf_1',\Gbf_2')$ are defined as
\begin{equation*}
    \begin{aligned}
        \sigma(x)&=-\sigma_2(x)+\sigma_1\left(U^{(1)}_{\sigma_2(x)}A^{(1)}_{p_0}x\right)+\sigma_0\left(U^{(1)}_{\sigma_2(x)}A^{(1)}_{p_0}x\right),\\
        \alpha_1(\cdot)&=(g_1g_0^{-1})\alpha(\cdot)(g_1g_0^{-1})^{-1}.
    \end{aligned}
\end{equation*}
Thus, $\psi_{(3)}$ is affine up to $U^{(2)}_t$-translations, and $(g_1,\alpha_1)\in\Gbf_2'\times\Sur(\Gbf_1',\Gbf_2')$ is an affine pair associated with $\psi_{(3)}$ up to $U^{(2)}_t$-translations.

\medskip

Next, assuming that $\psi_{(3)}$ is a Kakutani equivalence, we need to show that $\alpha_1(\Gamma_1')=g_1\Gamma_2'g_1^{-1}$ and $\alpha_1$ is an isomorphism. By the constructions of $g_1$ and $\alpha_1$ in (\ref{eq:ProofOfTheorem:affinePair}), it suffices to show that in this case, $\alpha(\Gamma_1')=g_0\Gamma_2'g_0^{-1}$ and $\alpha$ is an isomorphism.

Since $\psi_{(3)}$ is a Kakutani equivalence, it follows from \textit{Proposition~\ref{Prop:epsilon-goodTimeChange:epsilon-goodTimeChange}}, \textit{Remark \ref{Remark:BehaviorOfpsi:Assume(g_(psi))=I_(d_2+1)}}, and \textit{Remark \ref{Remark:BehaviorOfpsiW_t:Assumepsi=Tilde(psi)}} that $\psi_{(2)}$, $\psi_{(1)}$, and $\psi$ are also Kakutani equivalences. With this, using \textit{Claim~\ref{Claim:ProofOfTheorem:phiIsMeasurePreserving}}, \textit{Claim~\ref{Claim:ProofOfTheorem:phi(U_t,A_t,Uh_t)}}, \eqref{eq:BehaviorOfpsi:|partial_svarpi_i(x,s)-1|IsSmall}, and \eqref{eq:BehaviorOfpsi:ConjugationRelation}, we deduce that for some full measure subset $Z\subseteq\Gbf_1/\Gamma_1$, $\phi\vert_Z$ is injective. Hence, the following claim holds:
\begin{Claim}\label{Claim:ProofOfCoroll:phiIsInjective}
    $\phi:\Gbf_1/\Gamma_1\to\Gbf_2/\Gamma_2$ is a measurable isomorphism between Borel probability spaces $(\Gbf_i/\Gamma_i,\mu_i)$, $i=1,2$.\\
    Consequently, $\phi':\Gbf_1'/\Gamma_1'\to\Gbf_2'/\Gamma_2'$ is also a measurable isomorphism between Borel probability spaces $(\Gbf_i'/\Gamma_i',\mu_i')$, $i=1,2$.
\end{Claim}

With this claim and condition \ref{item:ProofOfMainTheorem:3} on p.~\pageref{item:ProofOfMainTheorem:3}, by the summary on the rigidity of unipotent flows under measurable isomorphism in \textit{Section~\ref{Section:Introduction}} and an argument similar to that of \cite{MR657240}*{Corollary 1}, it follows that $\alpha(\Gamma'_1)=g_0\Gamma_2'g_0^{-1}$ and $\alpha:\Gbf_1'\to\Gbf_2'$ is an isomorphism.

\medskip

Thus, we finish the proof of \textit{Theorem~\ref{thm:main}}.\qed

\bibliographystyle{plain} 
\bibliography{refs}

\end{document}